\documentclass[final, sort&compress]{siamart1116}

\usepackage{amsfonts}
\usepackage{amsmath}
\usepackage{amssymb}
\usepackage{bbm}
\usepackage{graphicx}
\usepackage{mathrsfs}
\usepackage{tikz-qtree}
\usepackage{rotating}
\usepackage{scalerel}
\usepackage{subcaption}
\usepackage{xspace}
\usepackage{comment}
\usepackage{array}
\usepackage{multirow}
\usepackage{changepage}

\usetikzlibrary{patterns}
\newcommand{\A}{\mathbf{A}}
\renewcommand{\b}{\mathbf{b}}

\renewcommand{\c}{\mathbf{c}}

\newcommand{\f}{\mathfrak{f}}
\newcommand{\s}{\mathfrak{s}}
\newcommand{\coup}{\mathfrak{c}}
\newcommand{\lte}{\textsc{lte}}
\newcommand{\one}{\text{\usefont{U}{bbm}{m}{n}1}} 
\renewcommand{\Re}{\mathbb{R}}
\newcommand{\norm}[1]{\bigl\lVert#1\bigr\rVert}
\newcommand{\normerr}[1]{\bigl\lVert#1\bigr\rVert_\textnormal{err}}
\newcommand{\abs}[1]{\lvert#1\rvert}
\newcommand{\bigo}[1]{\mathcal{O}{\left(#1\right)}}
\newcommand{\trsym}{^T}
\newcommand{\tr}[1]{#1\,\trsym}
\newcommand{\mgark}{MrGARK\xspace}
\newcommand{\mprod}{\,}
\newcommand{\sprod}{\times}
\newcommand{\explicit}{EX}
\newcommand{\implicit}{IM}
\newcommand{\method}[7]{#1#5-#2#6 #3(#4)[#7]}
\newcommand{\fitbox}[2]{%
	\noindent
	\resizebox{\linewidth}{!}{%
		\bgroup
		\renewcommand{\arraystretch}{3.5}
		\newcommand{\fullrow}[1]{%
			\multicolumn{#1}{>{\renewcommand{\arraystretch}{1.4}}l}{##1}%
		}
		\begin{tabular}{*#1{>{\renewcommand{\arraystretch}{1.4}}l}}
			#2
		\end{tabular}
		\egroup%
	}%
}
\newcommand{\fittboxrotate}[2]{
	\noindent
	\rotatebox{90}{\resizebox{0.9 \textheight}{!}{%
			\bgroup
			\renewcommand{\arraystretch}{5.5}
			\newcommand{\fullrow}[1]{%
				\multicolumn{#1}{>{\renewcommand{\arraystretch}{1.4}}l}{##1}%
			}
			\begin{tabular}{*#1{>{\renewcommand{\arraystretch}{1.4}}l}}
				#2
			\end{tabular}
			\egroup%
	}}
}

\newif\iftechreport
\techreporttrue

\allowdisplaybreaks
\newtheorem{remark}{Remark}
\numberwithin{theorem}{section}
\numberwithin{remark}{section}

\newcommand{\TheTitle}{Design of High-order Decoupled Multirate GARK Schemes} 
\newcommand{\TheAuthors}{Arash Sarshar \and Steven Roberts \and Adrian Sandu}

\headers{\TheTitle}{\TheAuthors}

\title{{\TheTitle}\thanks{
		\funding{The work of A. Sandu has been supported in part by NSF through awards NSF
			OCI--8670904397, NSF CCF--0916493, NSF DMS--0915047, NSF CMMI--1130667, 
			NSF CCF--1218454, AFOSR FA9550--12--1--0293--DEF, AFOSR 12-2640-06,
			and by the Computational Science Laboratory at Virginia Tech.}}}

\author{
	Arash Sarshar\thanks{
		Virginia Polytechnic Institute and State University, 
		Computational Science Laboratory, Department of Computer 
		Science, 2202 Kraft Drive, Blacksburg, VA 24060, USA
		(\email{sarshar@vt.edu}, \email{steven94@vt.edu}, \email{sandu@cs.vt.edu})
	}
	\and Steven Roberts\footnotemark[2]
	\and Adrian Sandu\footnotemark[2]
}

\begin{document}
\iftechreport

\thispagestyle{empty}
\setcounter{page}{0}
\begin{adjustwidth}{-2cm}{-1cm}

\makeatletter
\def\Year#1{%
  \def\yy@##1##2##3##4;{##3##4}%
  \expandafter\yy@#1;
}
\makeatother

\begin{Huge}
\begin{center}
\hspace{2cm} Computational Science Laboratory \\ \hspace{1cm} Technical Report CSL-TR-\Year{\the\year}-{4} \\
\today
\end{center}
\end{Huge}
\vfil
\begin{Large}
\begin{center}
Arash Sarshar, Steven Roberts and Adrian Sandu
\end{center}
\end{Large}

\vfil
\begin{Large}
\begin{it}
\begin{center}
\hspace*{2cm}``{Design of High-Order Decoupled Multirate GARK Schemes}''
\end{center}
\end{it}
\end{Large}
\vfil

\begin{large}
\begin{center}
Computational Science Laboratory \\
``Compute the Future!'' \\[12pt]
Computer Science Department \\
Virginia Polytechnic Institute and State University \\
Blacksburg, VA 24060 \\
Phone: (540)-231-2193 \\
Fax: (540)-231-6075 \\ 
Email: \url{sarshar@vt.edu} \\
Web: \url{http://csl.cs.vt.edu}
\end{center}
\end{large}
\vspace{1cm}

\resizebox{1 \textwidth}{!}{
	\hspace{50cm}
\begin{tabular}{cp{3in}c}
\raisebox{0pt}{\includegraphics[width=3\textwidth]{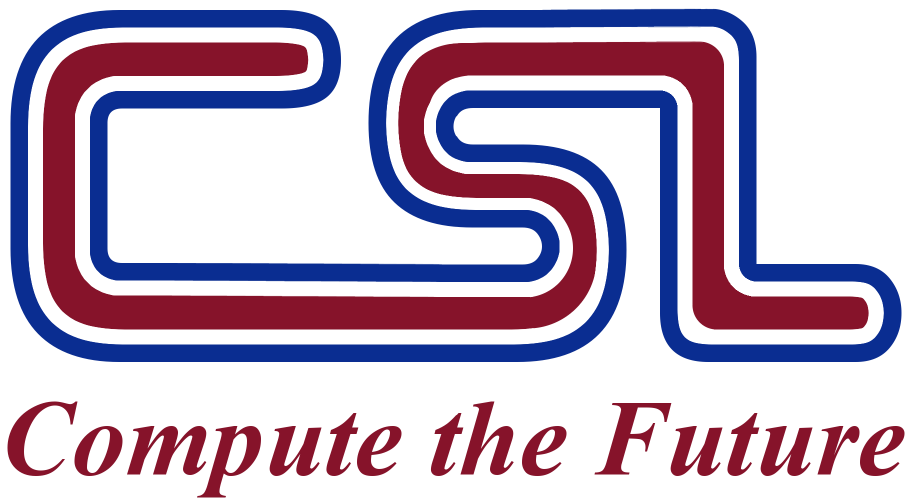}}
&&
\raisebox{10pt}{\includegraphics[width= 3.5 \textwidth]{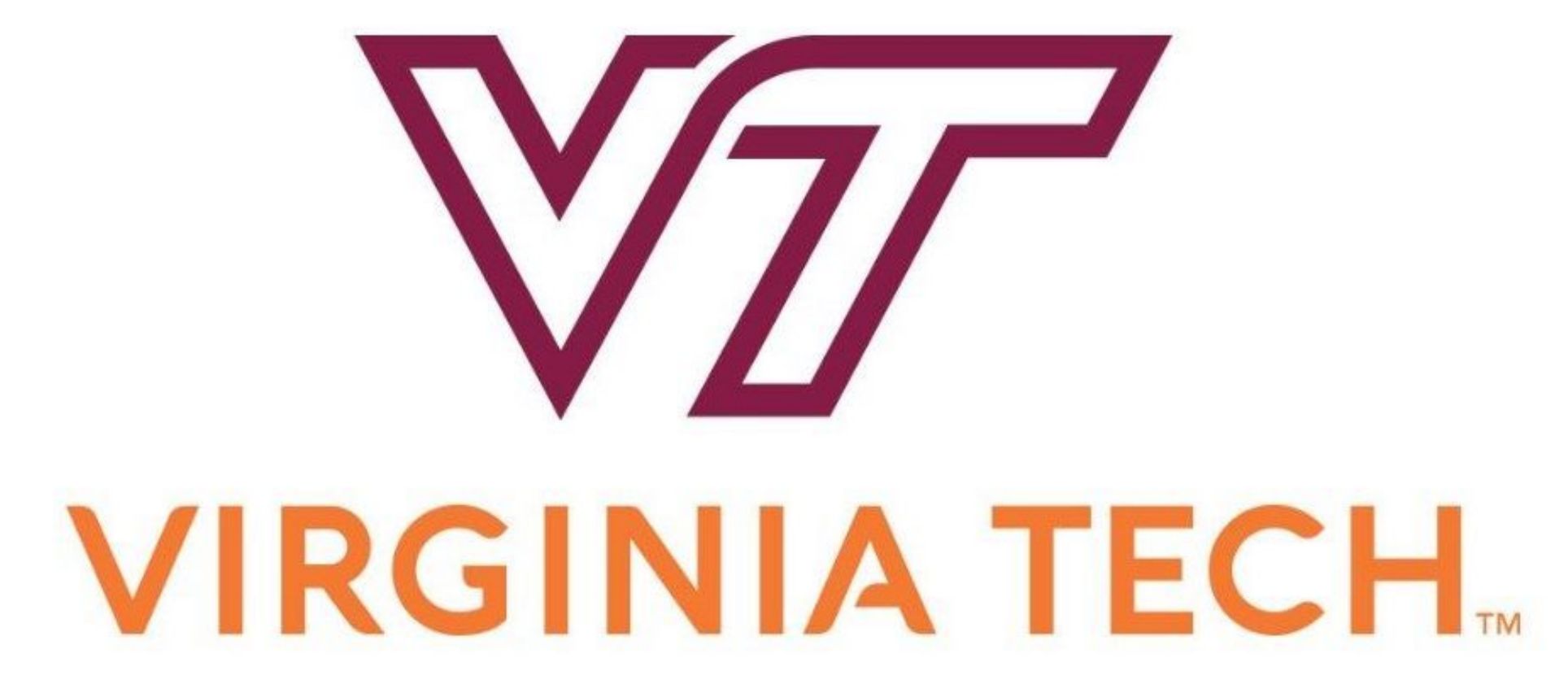}} 
\end{tabular}
}
\end{adjustwidth}

\clearpage
\fi

\maketitle

\begin{abstract}
Multirate time integration methods apply different step sizes to resolve different components of the system based on the local activity levels. This local selection of step sizes allows increased computational efficiency while achieving the desired solution accuracy. While the multirate idea is elegant and has been around for decades, multirate methods are not yet widely used in applications. This is due, in part, to the difficulties raised by the construction of high order multirate schemes.

Seeking to overcome these challenges, this work focuses on the design of practical high-order multirate methods using the theoretical framework of generalized additive Runge--Kutta (\mgark) methods \cite{Sandu_2016_GARK-MR}, which provides the generic order conditions and the linear and nonlinear stability analyses.   
A set of design criteria for practical multirate methods is defined herein: method coefficients should be generic in the step size ratio, 
but should not depend strongly on this ratio; unnecessary coupling between the fast and the slow components should be avoided; and the step size controllers should adjust both the micro- and the macro-steps.
Using these criteria, we develop \mgark schemes of up to order four that are explicit-explicit (both the fast and slow component are treated explicitly), implicit-explicit (implicit in the fast component and explicit in the slow one), and explicit-implicit (explicit in the fast component and implicit in the slow one). Numerical experiments illustrate the performance of these new schemes.
\end{abstract}



\clearpage

\section{Introduction}
\label{sec:introduction}

Many applications in science and engineering require the simulation of dynamical systems where different components evolve at different characteristic time scales. Clearly, these systems challenge traditional time discretizations that use a single time step for the entire system: either the fast components are resolved inaccurately, or the slow components are resolved with more accuracy than required, therefore increasing computational costs. Multirate methods exploit the partitioning of a system into components with different time scales, and use small step sizes to discretize the fast components, and large step sizes to discretize the slow components.

A first approach to constructing multirate methods is to employ traditional integrators with different time steps for different components, and to carefully orchestrate the coupling between these components. Early efforts to develop multirate Runge-Kutta methods are due to Rice \cite{Rice_1960} and Andrus \cite{Andrus_1979,Andrus_1993}.  In the first discussion of ``multirate methods" Gear and Wells \cite{Gear_1984_MR-LMM} propose pairing various linear multistep methods. This fundamental contribution already points to a number of challenges facing multirate methods such as coupling, automatic step size selection, and efficiency of the overall computational process. Other work to construct multirate linear multistep schemes includes \cite{Kato_1999}. G\"unther et al. \cite{Guenther_2001_MR-PRK,Guenther_1994_partition-circuits} developed multirate methods for partitioned Runge-Kutta schemes, as well as Rosenbrock-W methods  \cite{Guenther_1997_ROW} of order three that are well-suited for treatment of systems with both stiff and non-stiff variables. Similarly, Kv{\ae}rn{\o} and Rentrop  \cite{Kvaerno_2000_stability-MRK,Kvaerno_1999_MR-RK} constructed explicit multirate Runge-Kutta methods of order three. Bartel et al. \cite{Bartel_2002_MR-W} propose one-step methods where internal stages are used to provide the coupling between the fast and slow components. Constantinescu and Sandu developed strong stability preserving (SSP) multirate methods of Runge-Kutta  \cite{Sandu_2007_MR_RK2} and linear multistep \cite{Sandu_2009_MR_LMM} type that are suited for solving hyperbolic partial differential equations (PDEs).

Another approach, tracing back to Engstler et al. \cite{Engstler_1997_MR-extrapolation}, derives multirate methods using Richardson extrapolation. Here, the solution is recursively improved on partitions of the system until required tolerances are satisfied. An important advantage of such schemes is the naturally available dynamic partitioning of the system. Constantinescu and Sandu  \cite{Sandu_2008_MR-EXTRAP,Sandu_2013_extrapolatedMR,Sandu_2009_ICNAAM-Multirate} considered explicit and implicit base methods for multirate extrapolation methods and study their stability properties. Other multirate approaches include Galerkin \cite{Logg_2003_MAG1}, and combined multiscale \cite{Engquist_2005} methodologies.

G\"unther and Sandu  \cite{Sandu_2016_GARK-MR} built a class of multirate methods in based on the General Additive Runge--Kutta framework (GARK) \cite{Sandu_2015_GARK}. Bremicker-Tr\"{u}belhorn and Ortleb \cite{Bremicker_2017_MGARK} developed third order multirate GARK (\mgark) methods for fluid-structure interaction, and allowed for non-uniform fast steps in the order conditions of the methods.

This study develops a systematic design approach for constructing multirate methods in the \mgark framework of G\"unther and Sandu \cite{Sandu_2016_GARK-MR,Sandu_2015_GARK}. Several high-order schemes are constructed that combine implicit and explicit components for the fast and slow subsystems.

The paper is organized as follows. \Cref{sec:intro} reviews the GARK framework and multirate GARK family. We study order conditions for high order \mgark methods in \cref{sec:order-conditions,}, and their scalar linear stability in \cref{sec:linear-stability}. \Cref{sec:coupling} provides insight into fast-slow coupling requirements. \Cref{sec:design} discusses the design criteria for practical methods, and \cref{sec:step_control} studies error estimation and the adaptivity of both micro- and macro-steps. Newly constructed methods are listed in \cref{sec:new_mgark}, and method coefficients are detailed in \cref{sec:new_mr_schemes}. Numerical tests are performed on different test problems and the results are reported in \cref{sec:numerics}.

\section{Multirate generalized additive Runge--Kutta schemes (\mgark)}
\label{sec:intro} 

In this section we review some background on \mgark methods.

\subsection{GARK methods}
\label{subsec:gark-review}
The generalized additive Runge--Kutta (GARK) methods introduced in \cite{Sandu_2014_SCEE,Sandu_2015_GARK} allow derivation of advanced multi-methods for solving {\it additively} partitioned systems of ordinary differential equations: 
\begin{equation} 
	\label{eqn:additive-ode}
	 y'= f(y) = \sum_{m=1}^N f^{\{m\}}(y), \qquad y(t_0)=y_0,
\end{equation}
where the right-hand side function $f: \Re^d \rightarrow \Re^d$ is split into {$N$} different partitions  based on properties such as stiffness, non-linearity, dynamical behavior, and evaluation cost. We note that additive partitioning includes the case of \textit{component partitioning} in which the state vector is split into a number subsets.

A GARK method advances the numerical solution as follows \cite{Sandu_2015_GARK}:
\begin{subequations}
\label{eqn:gark-step}
\begin{align}
	\label{eqn:delta-stage2}
	Y_i^{\{q\}} &= y_{n} + h \sum_{m=1}^N \sum_{j=1}^{s^{\{m\}}} a_{i,j}^{\{q,m\}} \, f^{\{m\}}\left(Y_j^{\{m\}}\right)
	\qquad \substack{{q=1,\ldots,N,}\\ {i=1,\ldots,s^{\{q\}}}},  \\
	\label{eqn:delta-sol2}
	y_{n+1} &= y_n + h \sum_{q=1}^N \sum_{i=1}^{s^{\{q\}}} b_{i}^{\{q\}} \, f^{\{q\}}\left(Y_i^{\{q\}}\right).
\end{align}
\end{subequations}
%

\subsection{Multirate GARK methods}
\label{subsec: review}

In the case of a two-way partitioned system \cref{eqn:additive-ode} with slow component $\{\s\}$, and fast component $\{\f\}$ we have:
\begin{equation} 
	\label{eqn:multirate-ode}
	 y'= f(y) = f^{\{\s\}} (y) +  f^{\{\f\}} (y), \qquad y(t_0)=y_0.
\end{equation}
A multirate GARK method \cite{Sandu_2015_GARK} integrates the slow component with a Runge--Kutta method $(A^{\{\s,\s\}},b^{\{\s\}})$ and a large step size $H$, and the fast component with another Runge--Kutta method $\bigl(A^{\{\f,\f\}},b^{\{\f\}}\bigr)$ and a small step size $h = H/M$. Here $M \ge 1$ represents the (integer) number of fast steps that are executed for each of the slow steps.
\mgark methods are formally derived in \cite{Sandu_2016_GARK-MR}. One step of the method utilizes $ s^{\{\s\}}$ slow stages, denoted by $Y_i^{\{\s\}}$, and $M\, s^{\{\f\}}$ fast stages, denoted by $Y_i^{\{\f,\lambda\}}$:
\begin{subequations}
	\label{eqn:GARK-MR}
	\begin{align}
	        \label{eqn:GARK-MR-slow-stage}
		\begin{split}
			Y_i^{\{\s\}} &= y_n + H \, \sum_{j=1}^{s^{\{\s\}}} a_{i,j}^{\{\s,\s\}} f^{\{\s\}}\left(Y_j^{\{\s\}}\right) \\
			& \quad + h \, \sum_{\lambda=1}^M \,\sum_{j=1}^{s^{\{\f\}}} a_{i,j}^{\{\s,\f,\lambda\}} f^{\{\f\}}\left(Y_j^{\{\f,\lambda\}}\right),
			\qquad i=1,\dots, s^{\{\s\}},
		\end{split} \\
		\begin{split}
			\label{eqn:GARK-MR-fast-stage}
			Y_i^{\{\f,\lambda\}} & = \widetilde{y}_{n+\left(\lambda-1\right)/M} + H \, \sum_{j=1}^{s^{\{\s\}}} a_{i,j}^{\{\f,\s,\lambda\}} f^{\{\s\}}\bigl(Y_j^{\{\s\}}\bigr) \\
			& \quad + h \, \sum_{j=1}^{s^{\{\f\}}} a_{i,j}^{\{\f,\f\}} f^{\{\f\}}\bigl(Y_j^{\{\f,\lambda\}}\bigr),\qquad \substack{i=1,\dots, s^{\{\f\}}\\{\lambda=1,\ldots,M }},
		\end{split} \\
	        \label{eqn:GARK-MR-final-solution}
		\widetilde{y}_{n+\lambda/M} &= \widetilde{y}_{n+\left(\lambda-1\right)/M} + h \sum_{i=1}^{s^{\{\f\}}} b_{i}^{\{\f\}} f^{\{\f\}}\bigl(Y_i^{\{\f,\lambda\}}\bigr).
	\end{align}
	The full step is then formed with:
	\begin{equation}
	y_{n+1} = \widetilde{y}_{n+M/M} + H \, \sum_{i=1}^{s^{\{\s\}}} b_{i}^{\{\s\}} f^{\{\s\}}\bigl(Y_i^{\{\s\}}\bigr).
	\end{equation}
\end{subequations}
Let the coupling between the two methods be described by $A^{\{\s,\f,\lambda\}}$ and $A^{\{\f,\s,\lambda\}}$ for $\lambda \in \{ 1, 2, \cdots, M \}$. The GARK Butcher tableau for the method \eqref{eqn:GARK-MR} is \cite{Sandu_2016_GARK-MR}:
\begin{equation}
	\label{eqn:mrRK-butcher}
	\bgroup
	\def\arraystretch{1.75}
	\begin{array}{c|c}
	\A^{\{\f,\f\}} & \A^{\{\f,\s\}}  \\ \hline
	\A^{\{\s,\f\}} & \A^{\{\s,\s\}} \\ \hline 
	\tr{{\b}^{\{\f\}}} & \tr{{\b}^{\{\s\}}}
	\end{array} ~~ :=~~ 
	\begin{array}{cccc|cccc}  
	\frac{1}{M} A^{\{\f,\f\}}      &          0                   & \cdots & 0 & A^{\{\f,\s,1\}}  \\
	\frac{1}{M} \one \mprod  \tr{b^{\{\f\}}} & \frac{1}{M} A^{\{\f,\f\}}        & \cdots & 0 &  A^{\{\f,\s,2\}}  \\
	\vdots                     &                             & \ddots &   & \vdots  \\
	\frac{1}{M} \one \mprod \tr{b^{\{\f\}}} & \frac{1}{M} \one \mprod \tr{b^{\{\f\}}} & \cdots & \frac{1}{M} A^{\{\f,\f\}} &A^{\{\f,\s,M\}} \\
	\hline 
	\frac{1}{M} A^{\{\s,\f,1\}} & \frac{1}{M} A^{\{\s,\f,2\}} & \cdots & \frac{1}{M} A^{\{\s,\f,M\}} & A^{\{\s,\s\}}   \\   \hline 
	\frac{1}{M} \tr{b^{\{\f\}}} & \frac{1}{M} \tr{b^{\{\f\}}} & \cdots & \frac{1}{M} \tr{b^{\{\f\}}} & \tr{b^{\{\s\}}}
	\end{array}~.
	\egroup
\end{equation}
\begin{remark}[Slow and fast stage numbers]
The structure of the Butcher tableau implies that the fast stage $\ell$ of the fast micro-step $\lambda$ corresponds to row $(\lambda-1) s^{\{\f\}} + \ell$ in  \cref{eqn:mrRK-butcher}, and the slow stage $j$ corresponds to row $M s^{\{\f\}} + j$ in  \cref{eqn:mrRK-butcher}.
\end{remark}
\begin{definition}[Telescopic \mgark schemes]
A \textit{telescopic \mgark} method \cite{Sandu_2016_GARK-MR} uses the same base scheme for the slow and fast partitions:
\begin{equation}
	\label{eqn:mr-telescopic-conditions}
	A^{\{\f,\f\}} = A^{\{\s,\s\}} = A, \qquad b^{\{\f\}} = b^{\{\s\}} = b.
\end{equation}
\end{definition}
In this paper we will focus on schemes with telescopic property since it allows a simple extension of the \mgark method to an arbitrary number of partitions (time scales) using successively larger time steps, each an integer multiple of the previous one. 
%
\section{Order conditions for multirate GARK methods}
\label{sec:order-conditions}
We consider the following \textit{internal consistency} conditions \cite{Sandu_2016_GARK-MR,Sandu_2015_GARK} to ensure that fast and slow right-hand side evaluations are performed at the same points in time:
\begin{equation}
\label{eqn:internal-consistency}
\A^{\{\s,\f\}} \mprod {\one}^{\{ \s \}}  =  \A^{\{\s,\s\}}\mprod {\one}^{\{ \s\}} :={\c^{\{\s\}}}, \qquad
\A^{\{\f,\s\}} \mprod {\one}^{\{ \f \}}  =  \A^{\{\f,\f\}}\mprod {\one}^{\{ \f \}} :={\c^{\{\f\}}},
\end{equation}
where ${\one}^{\{ \sigma \}}  := [1, 1, \cdots, 1 ]\trsym \in \Re^{s^{\{\sigma\}}}$. 

Assume that the internal consistency conditions \cref{eqn:internal-consistency} hold, and further assume that each individual method $\bigl(A^{\{\s,\s\}},b^{\{\s\}}\bigr)$ and $\bigl(A^{\{\f,\f\}},b^{\{\f\}}\bigr)$ has at least order 4. Then the \mgark scheme \eqref{eqn:GARK-MR}--\eqref{eqn:mrRK-butcher} has order four if and only if the following coupling conditions hold \cite{Sandu_2016_GARK-MR}: 
\begin{subequations}
\label{eq:mgark-high-order-conditions}
\begin{align}
	  \frac{1}{6} &= \tr{\b^{\{\f\}}} \mprod \A^{\{\f,\s\}}  \mprod \c^{\{\s\}}
, & \qquad (\text{order}~ 3) \\
	  \frac{1}{6} &= \tr{\b^{\{\s\}}} \mprod \A^{\{\s,\f\}}  \mprod \c^{\{\f\}}
, & \qquad (\text{order}~ 3) \\
	  \frac{1}{8} &= \tr{\b^{\{\f\}}} \mprod \left(\c^{\{\f\}} \sprod \A^{\{\f,\s\}} \mprod \c^{\{\s\}}\right)
, & \qquad (\text{order}~ 4) \\
	  \frac{1}{8} &= \tr{\b^{\{\s\}}} \mprod \left(\c^{\{\s\}} \sprod \A^{\{\s,\f\}} \mprod \c^{\{\f\}}\right)
, & \qquad (\text{order}~ 4) \\
	  \frac{1}{12} &= \tr{\b^{\{\f\}}} \mprod \A^{\{\f,\s\}}  \mprod
	\left( \c^{\{\s\}} \sprod \c^{\{\s\}} \right)
,  & \qquad (\text{order}~ 4) \\
	  \frac{1}{12} &= \tr{\b^{\{\s\}}} \mprod \A^{\{\s,\f\}}  \mprod
	\left( \c^{\{\f\}} \sprod \c^{\{\f\}} \right)
,  & \qquad (\text{order}~ 4) \\
	 \frac{1}{24} &= \tr{\b^{\{\s\}}} \mprod \A^{\{\s,\s\}} \mprod \A^{\{\s,\f\}} \mprod \c^{\{\f\}}, & \qquad (\text{order}~ 4)\\
	 \frac{1}{24} &= \tr{\b^{\{\s\}}} \mprod \A^{\{\s,\f\}} \mprod \A^{\{\f,\s\}} \mprod \c^{\{\s\}}, & \qquad (\text{order}~ 4)\\
	 \frac{1}{24} &= \tr{\b^{\{\s\}}} \mprod \A^{\{\s,\f\}} \mprod \A^{\{\f,\f\}} \mprod \c^{\{\f\}}, & \qquad (\text{order}~ 4)\\
	 \frac{1}{24} &= \tr{\b^{\{\f\}}} \mprod \A^{\{\f,\f\}} \mprod \A^{\{\f,\s\}} \mprod \c^{\{\s\}}, & \qquad (\text{order}~ 4)\\
	 \frac{1}{24} &= \tr{\b^{\{\f\}}} \mprod \A^{\{\f,\s\}} \mprod \A^{\{\s,\s\}} \mprod \c^{\{\s\}}, & \qquad (\text{order}~ 4)\\
	 \frac{1}{24} &= \tr{\b^{\{\f\}}} \mprod \A^{\{\f,\s\}} \mprod \A^{\{\s,\f\}} \mprod \c^{\{\f\}}. & \qquad (\text{order}~ 4)
\end{align}
\end{subequations}
Here ``$\times$'' denotes component-wise multiplication of two vectors.

Without imposing any special structure on coupling matrices, we can rewrite \cref{eq:mgark-high-order-conditions} using the block structure shown in \cref{eqn:mrRK-butcher}. 
\iftechreport
Considering the following intermediate simplifications:
\begin{subequations}
	\label{eq:order3_simplifying_terms}
	\begin{align}
	\c^{\{\f\}} &= \frac{1}{M} \left[ c^{\{\f\}} + \left(\lambda-1\right) \one  \right]_{\lambda=1,\dots,M},\\
	\A^{\{\f,\f\}} \mprod \c^{\{\f\}} &= \left[ \frac{1}{M^2}\left(\frac{(\lambda-1)^2}{2} \one + (\lambda-1) c^{\{\f\}} + A^{\{\f,\f\}} \mprod  c^{\{\f\}} \right) \right]_{\lambda=1,\dots,M},\\
	\A^{\{\s,\f\}}\mprod \c^{\{\f\}}
	&= \frac{1}{M^2} \sum_{\lambda=1}^M A^{\{\s,\f,\lambda\}} \mprod \left((\lambda - 1) \one + c^{\{\f\}}\right), \\
	\A^{\{\f,\f\}} \mprod \A^{\{\f,\s\}} \mprod \c^{\{\s\}} &= \frac{1}{M}\left[ \sum_{k=1}^{\lambda -1} \one \mprod \tr{b^{\{\f\}}} \mprod  A^{\{\f,\s,k\}} \mprod c^{\{\s\}} +  A^{\{\f,\f\}} \mprod A^{\{\f,\s,\lambda\}} \mprod c^{\{\s\}}  \right]_{\lambda=1,\dots,M}.
	\end{align}
\end{subequations}
\fi 
The order conditions \cref{eq:mgark-high-order-conditions} can be written in terms of base methods and coupling coefficients as follows: 
\begin{subequations}
	\begin{small}
	\label{eq:mgark-high-order-conditions-blk}
	\begin{align}
	 \frac{M^{\phantom{.}}}{6} &= \sum_{\lambda=1}^M \tr{b^{\{\f\}}} \mprod A^{\{\f,\s,\lambda\}} \mprod c^{\{\s\}}, 
	& \qquad (\text{order}~ 3) \\
	\frac{M^2}{6} &= \sum_{\lambda=1}^M \tr{b^{\{\s\}}} A^{\{\s,\f,\lambda\}} \mprod \left((\lambda - 1) \one + c^{\{\f\}}\right),
	& \qquad (\text{order}~ 3) \\
	\begin{split}
  	\frac{M^2}{8} &= \sum_{\lambda=1}^M  (\lambda - 1)\;\tr{b^{\{\f\}}} \mprod A^{\{\f,\s,\lambda\}} \mprod c^{\{\s\}} \\
	&+\sum_{\lambda=1}^M \tr{b^{\{\f\}}} \mprod \left( c^{\{\f\}} \sprod A^{\{\f,\s,\lambda\}} \mprod c^{\{\s\}} \right),
	\end{split}
  	& \qquad (\text{order}~ 4) \\
	\frac{M^2}{8} &= \tr{b^{\{\s\}}} \mprod  \sum_{\lambda=1}^M \Bigg( c^{\{\s\}} \sprod  \Big( A^{\{\s,\f\, \lambda\}} \big( (\lambda-1)\; \one  + c^{\{\f\}}  \big ) \mprod   \Big) \Bigg),
	& \qquad (\text{order}~ 4) \\
	\frac{M^{\phantom{.}}}{12} &= \sum_{\lambda=1}^M  \tr{b^{\{\f\}}} \mprod A^{\{\f,\s,\lambda\}}  c^{\{\s\} \sprod 2},
	& \qquad (\text{order}~ 4) \\
	\begin{split}
	\frac{M^3}{12} &= \sum_{\lambda=1}^M  \tr{b^{\{\s\}}} \mprod A^{\{\s,\f,\lambda\}} \mprod c^{\{\f\} \sprod 2} 
 	+\sum_{\lambda=1}^M (\lambda -1 )^2 \tr{b^{\{\s\}}} \mprod A^{\{\s,\f,\lambda\}} \mprod \one \\
	&+2\sum_{\lambda=1}^M (\lambda -1 ) \tr{b^{\{\s\}}} \mprod A^{\{\s,\f,\lambda\}} \mprod c^{\{\f\}},
	\end{split} 	& \qquad (\text{order}~ 4) \\
	\frac{M^2}{24} &= \sum_{\lambda=1}^M  \tr{b^{\{\s\}}} \mprod A^{\{\s,\s\}} \mprod A^{\{\s,\f,\lambda\}} \left( (\lambda -1) \one  + \mprod c^{\{\f\}} \right),
	& \qquad (\text{order}~ 4)\\
	\frac{M^{\phantom{.}}}{24} &= \sum_{\lambda=1}^M  \tr{b^{\{\s\}}} \mprod A^{\{\s,\f,\lambda\}} \mprod A^{\{\f,\s,\lambda\}} \mprod c^{\{\s\}},
	& \qquad (\text{order}~ 4)\\
	\begin{split}
		\frac{M^3}{24} &= \sum_{\lambda=1}^M \frac{(\lambda -1)^2}{2}  \tr{b^{\{\s\}}} \mprod A^{\{\s,\f,\lambda\}}\mprod \one \\
		&+\sum_{\lambda=1}^M (\lambda -1) \tr{b^{\{\s\}}} \mprod A^{\{\s,\f,\lambda\}} \mprod c^{\{\f\}}
		+\sum_{\lambda=1}^M \tr{b^{\{\s\}}} \mprod A^{\{\s,\f,\lambda\}}\mprod A^{\{\f,\f\}}\mprod c^{\{\f\}},
	\end{split} 
	& \qquad (\text{order}~ 4) \\
	\frac{M^2}{24} &= \sum_{\lambda=1}^M  \sum_{k=1}^{{\lambda-1}}  \tr{b^{\{\f\}}} \mprod A^{\{\f,\s,k\}}\mprod c^{\{\s\}} + 
	\sum_{\lambda=1}^M  \tr{b^{\{\f\}}} \mprod A^{\{\f,\f\}} \mprod A^{\{\f,\s,\lambda\}}\mprod c^{\{\s\}},
	& \qquad (\text{order}~ 4)\\
	\frac{M^{\phantom{.}}}{24} &= \sum_{\lambda=1}^M  \tr{b^{\{\f\}}} \mprod A^{\{\f,\s,\lambda\}} \mprod A^{\{\s,\s\}} \mprod  c^{\{\s\}},
	& \qquad (\text{order}~ 4) \\
	\frac{M^3}{24} &= \sum_{\lambda=1}^M  \sum_{k=1}^M  \tr{b^{\{\f\}}} \mprod A^{\{\f,\s,\lambda\}} \mprod A^{\{\s,\f, k\}} \mprod \left( (k-1)\one + c^{\{\f\}} \right). 
	& \qquad (\text{order}~ 4) 
	\end{align}
\end{small}
\end{subequations}
\begin{remark}[Design process]
A practical design procedure for \mgark schemes is to first select the base slow and fast methods of desired order, and then solve the order conditions \cref{eq:mgark-high-order-conditions-blk} for the coupling coefficients $ \A^{\{\f,\s\}} $ and $ \A^{\{\s,\f\}} $.
\end{remark}
\section{Linear stability analysis}
\label{sec:linear-stability}
Following \cite{Sandu_2016_GARK-MR,Sandu_2015_GARK} we consider the following scalar model problem:
\begin{equation}
\label{eq:linear-test-problem}
y' = \lambda^{\{\f\}}\, y + \lambda^{\{\s\}}\, y.
\end{equation}
where the ratio of the fast to the slow variable is $M$, the step size ratio. 
Denote $z^{\{\f\}}=H\, \lambda^{\{\f\}}$, $z^{\{\s\}}=H\, \lambda^{\{\s\}}$, and
\begin{equation}
s = M s^{\{\f\}} + s^{\{\f\}} , \quad
Z = \begin{bmatrix}
z^{\{\f\}}\mathbf{I}_{ Ms^{\{\f\}} \times Ms^{\{\f\}}}  & \mathbf{0} \\
\mathbf{0} & z^{\{\s\}}\mathbf{I}_{ s^{\{\s\}} \times s^{\{\s\}}}
\end{bmatrix}.
\end{equation}
It was shown in \cite{Sandu_2016_GARK-MR,Sandu_2015_GARK} that application of \mgark method \cref{eqn:mrRK-butcher} to \cref{eq:linear-test-problem} leads to the solution
\begin{equation*}
y_{n+1} = R\left( z^{\{\f\}}, z^{\{\s\}} \right)  y_n,
\end{equation*}
with the stability function:
\begin{equation}
\label{eqn:stability-function}
R\left( z^{\{\f\}}, z^{\{\s\}} \right) = 1 + \mathbf{b}^T_\textsc{gark} \cdot Z   \cdot \left( \mathbf{I}_{s\times s} -   \mathbf{A}_\textsc{gark} \cdot Z \right)^{-1}\cdot \one_{s\times 1} .
\end{equation}
In order to visualize the stability region of a method, we choose:
\begin{equation}
	z^{\{\f\}} = M\, \rho\, e^{-i\, \theta^{\{\f\}}}, \qquad z^{\{\s\}} = \rho\, e^{-i\, \theta^{\{\s\}}}, 
	\qquad \frac{\pi}{2} \le \theta^{\{\f\}},\theta^{\{\s\}} \le \frac{3 \pi}{2},
\end{equation}
such that the ratio of the variable magnitudes matches the ratio of the step sizes.
This reduces \cref{eqn:stability-function} to a function of three real variables by assuming the fast eigenvalue is $M$ times larger in magnitude than the slow eigenvalue.   The stability region is the volume in the \{$\theta^{\{\f\}}$, $\theta^{\{\s\}}$, $\rho$\} space where the magnitude of the stability function \eqref{eqn:stability-function} is less than or equal to one. Stability regions for each method developed here are provided in \cref{sec:new_mr_schemes}. For some of the methods the region decreases with increasing $M$. The plots for different values of $M$ correspond to different test problems: when $M$ increases, the scale separation of the test problem also increases (we apply smaller micro-steps to faster problems).  The increasing stiffness of the fast component restricts the macro-step, a phenomenon known as stiffness leakage. Coupling coefficients whose magnitude increases rapidly with $M$ can lead to a degradation of stability.
\section{Decoupled \mgark methods}
\label{sec:coupling}

We now discuss in detail the structure of the coupling coefficient matrices, and how this structure defines the way the fast and slow stage computations \cref{eqn:GARK-MR} are carried out, and therefore determines the practicality of the \mgark method. In order to construct practical \mgark methods we need to avoid complex couplings between multiple fast and slow stages. To this end we define {\it decoupled \mgark methods}.

\begin{definition}[Decoupled \mgark methods]
An \mgark method is decoupled if the computation of its stages proceeds in sequence, such that each slow stage uses only information from other slow stages and the already computed fast stages, and vice-versa. There is no coupling that requires fast and slow stages to be solved together. Any form of implicitness is entirely within the fast or within the slow system.
\end{definition}

\subsection{Structure of the slow-fast coupling (including fast information into the slow stage calculations)}
\label{subsec:slow-to-fast-coup}

We first introduce a notation for the order of computation of the slow stages with respect to the fast stages. Consider the $j$-th slow stage \eqref{eqn:GARK-MR-slow-stage} at abscissa $c_j^{\{\s\}}$, i.e., the row $M\, s^{\{\f\}} + j$ in the Butcher tableau \cref{eqn:mrRK-butcher}. We denote its order of computation by $(L_j,I_j)$, i.e., the $j$-th slow stage is computed immediately after the $I_j$-th stage of the $L_j$-th micro-step is computed. This means that the first $L_j-1$ micro-steps have been completed, and the $L_j$-th micro-step has partially progressed to compute stage $I_j$, when the $j$-th slow stage is evaluated. 

When the slow stage $c_j^{\{\s\}}$ is evaluated after the last stage of micro-step $\lambda-1$, but before the before the first stage of the $\lambda$-th micro-step, we have $(L_j,I_j) = (\lambda-1,s^{\{\f\}})$. Equivalently, this situation can be represented by $(L_j,I_j) = (\lambda,0)$.

In order to construct practical \mgark methods we require that the evaluation of the $j$-th slow stage depends only on those fast stages that have been completed; in the GARK Butcher tableau \cref{eqn:mrRK-butcher}, the $j$-th slow stage depends only on the rows $1: (L_j-1) s^{\{\f\}} + I_j$. 

The $j$-th row of the slow-fast coupling matrix: 
\begin{subequations}
\label{eqn:structure-Asf}
\begin{equation}
	\A^{\{\s,\f\}} = \frac{1}{M}\,
	\begin{bmatrix} A^{\{\s,\f,1\}} & \dots & A^{\{\s,\f,\lambda\}} & \dots & A^{\{\s,\f,M\}} \end{bmatrix} \in \Re^{s^{\{\s\}} \times M s^{\{\f\}}}
\end{equation}
contains the coefficients that bring fast information into the computation of slow stage $j$. A decoupled \mgark scheme enjoys the following properties:
\begin{itemize}
	\item The coupling matrices $A^{\{\s,\f,\lambda\}}$ with the completed microsteps $\lambda=1:L_j-1$ can have full rows:
	\begin{equation*}
		a_{j,k}^{\{\s,\f,\lambda\}} \ne 0, \quad \text{ for } \quad \lambda=1:L_j-1.
	\end{equation*} 
	\item Rows of the coupling matrix $A^{\{\s,\f,L_j\}}$ can have non-zero entries only in the first $I_j$ positions:
	\begin{alignat*}{2}
		a_{j,k}^{\{\s,\f,L_j\}} &\ne 0, &\qquad k&=1,\cdots,I_j, \\
	        a_{j,k}^{\{\s,\f,L_j\}} &= 0,   &\qquad k&=I_j+1,\cdots,s^{\{\f\}}.
	\end{alignat*}
	\item The coupling matrices $A^{\{\s,\f,\lambda\}}$ with the future microsteps $\lambda=L_j+1:M$ are zero:
	\begin{equation*}
		 a_{j,k}^{\{\s,\f,\lambda\}} = 0, \quad \text{ for } \quad \lambda=L_j+1:M.
	\end{equation*}
\end{itemize}

Consequently, for decoupled \mgark schemes, the slow-fast coupling matrix is lower block triangular:
\begin{equation}
	\bigl(\A^{\{\s,\f\}}\bigr)_{j,k} = 0, \quad \textnormal{for} \quad  1 \le j \le s^{\{\s\}}, \quad (L_j-1) s^{\{\f\}} + I_j + 1 \le k \le M s^{\{\f\}},
\end{equation}
\end{subequations}
where we used the fact that the fast stage $I_j$  of micro-step $L_j$ has GARK stage number $(L_j-1) s^{\{\f\}} + I_j$ in the Butcher tableau \cref{eqn:mrRK-butcher}.
%

\subsection{Structure of the fast-slow coupling 
   (including slow information into the fast stage calculations)}
\label{subsec:fast-to-slow-coup}
\newcommand{\Jellam}{J_{\scaleto{L,i}{4pt}}}

We next introduce a notation for the order of computation of the fast stages with respect to the slow stages.  We denote by $\Jellam$ the index of the last slow stage computed before starting the fast stage $i$  of micro-step $L$. Specifically, the fast stage $i$ at micro-step $L$ is computed after the slow stage ${\Jellam}$, but before the slow stage ${\Jellam +1}$. 
%
For a decoupled \mgark this stage can only use information from slow stages $1$ to $\Jellam$. Consequently, the $i$-th row of the coupling matrix $A^{\{\f,\s,\lambda\}}$ can have nonzero entries only in the first $\Jellam$ columns:
\begin{align*}
	a^{\{\f,\s,\lambda\}}_{i,j} = 0, \quad \textnormal{for} \quad \Jellam+1 \le j \le s^{\{\s\}}.
\end{align*}
The coupling matrix:
\begin{subequations}
\label{eqn:structure-Afs}
\begin{equation}
	\A^{\{\f,\s\}} =
	\begin{bmatrix} A^{\{\f,\s,1\}}\,^T & \dots & A^{\{\f,\s,\lambda\}}\,^T & \dots & A^{\{\f,\s,M\}}\,^T \end{bmatrix}^T \in \Re^{M s^{\{\f\}}  \times s^{\{\s\}}}
\end{equation}
has a block lower triangular structure:
\begin{equation}
\A^{\{\f,\s\}}_{\scaleto{(L-1) s^{\{\f\}} + i,j}{6pt}} = 0 \quad \text{for} \quad \Jellam+1 \le j \le s^{\{\s\}},
\end{equation}
\end{subequations}
where we used the fact that the fast stage $i$  of micro-step $L$ has GARK stage number $(L-1) s^{\{\f\}} + i$ in the Butcher tableau \cref{eqn:mrRK-butcher}.
%
\subsection{Relation between the coupling matrices}

\begin{figure}[h]
\centering
\begin{subfigure}{0.475 \textwidth}
	\centering
	\includegraphics[width=\linewidth]{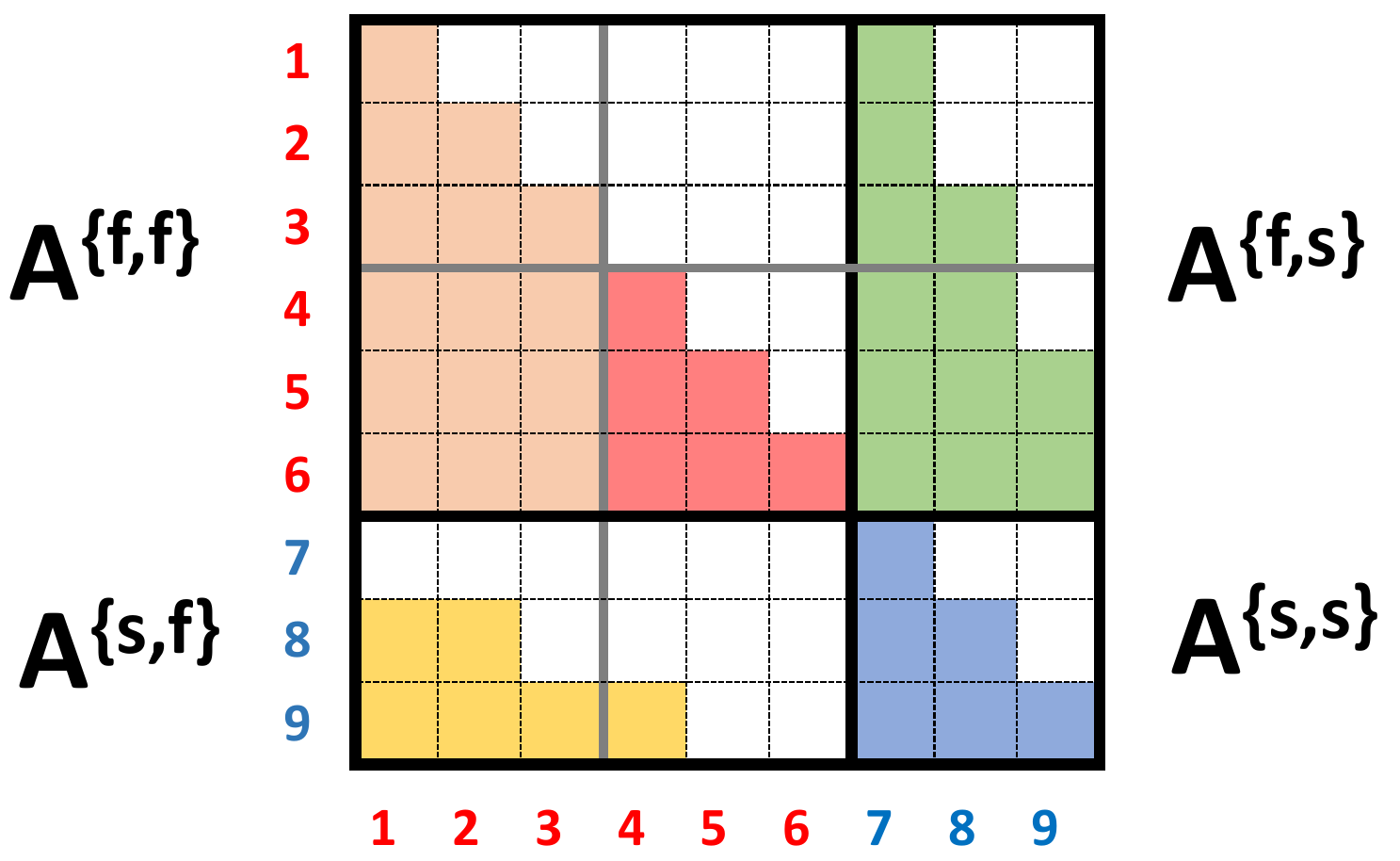}
	\caption{Butcher tableau of a decoupled \mgark. The coupling matrices have complementary sparsity patterns.}
	\label{fig:MGARK-butcher-decoupled-standard}
\end{subfigure}
\hfill
\begin{subfigure}{0.475 \textwidth}
	\centering
	\includegraphics[width=\linewidth]{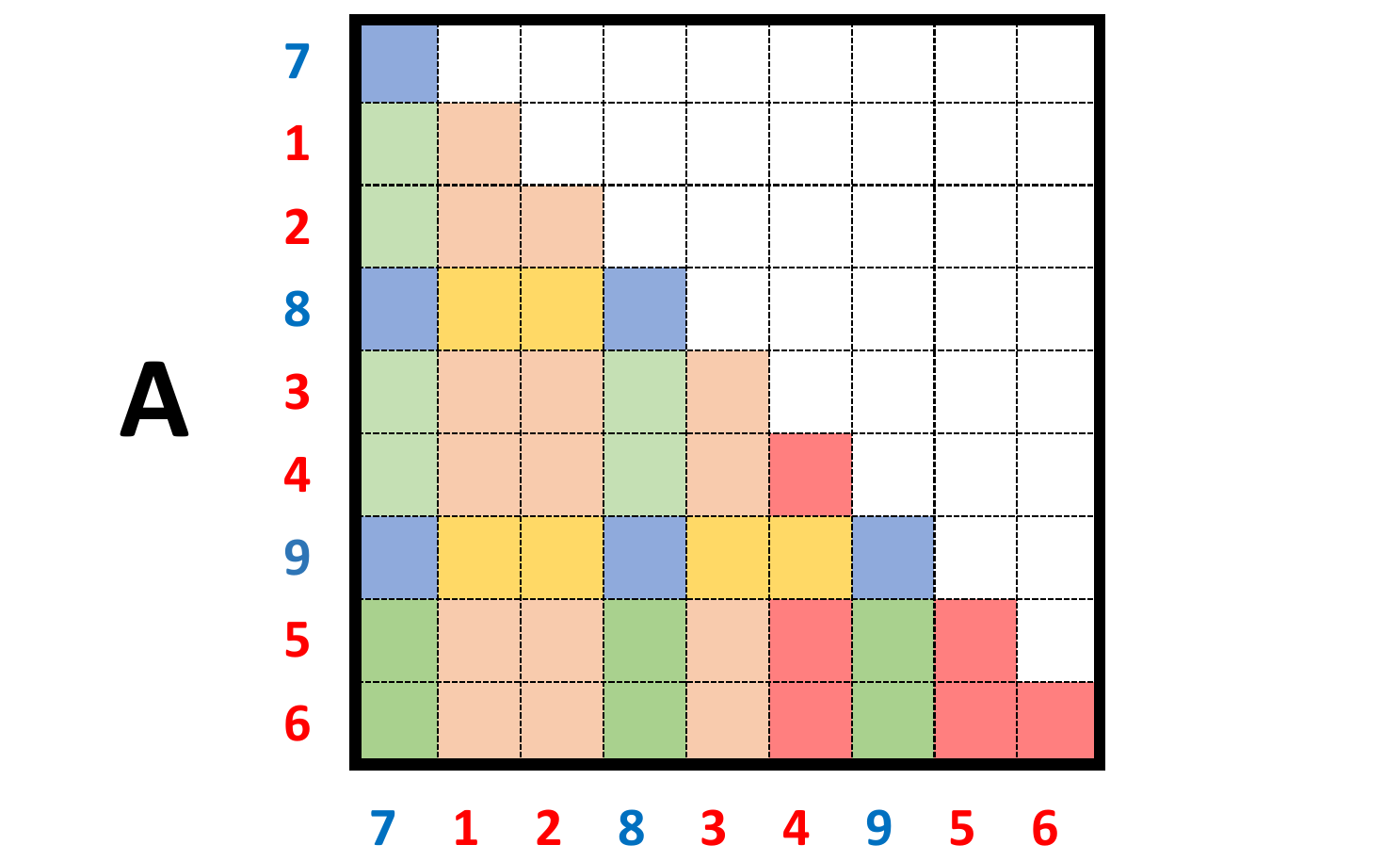}
	\caption{Permuted Butcher tableau of a decoupled \mgark shows a lower diagonal structure.}
	\label{fig:MGARK-butcher-decoupled-reordered}
\end{subfigure}

\begin{subfigure}{0.475 \textwidth}
	\centering
	\includegraphics[width=\linewidth]{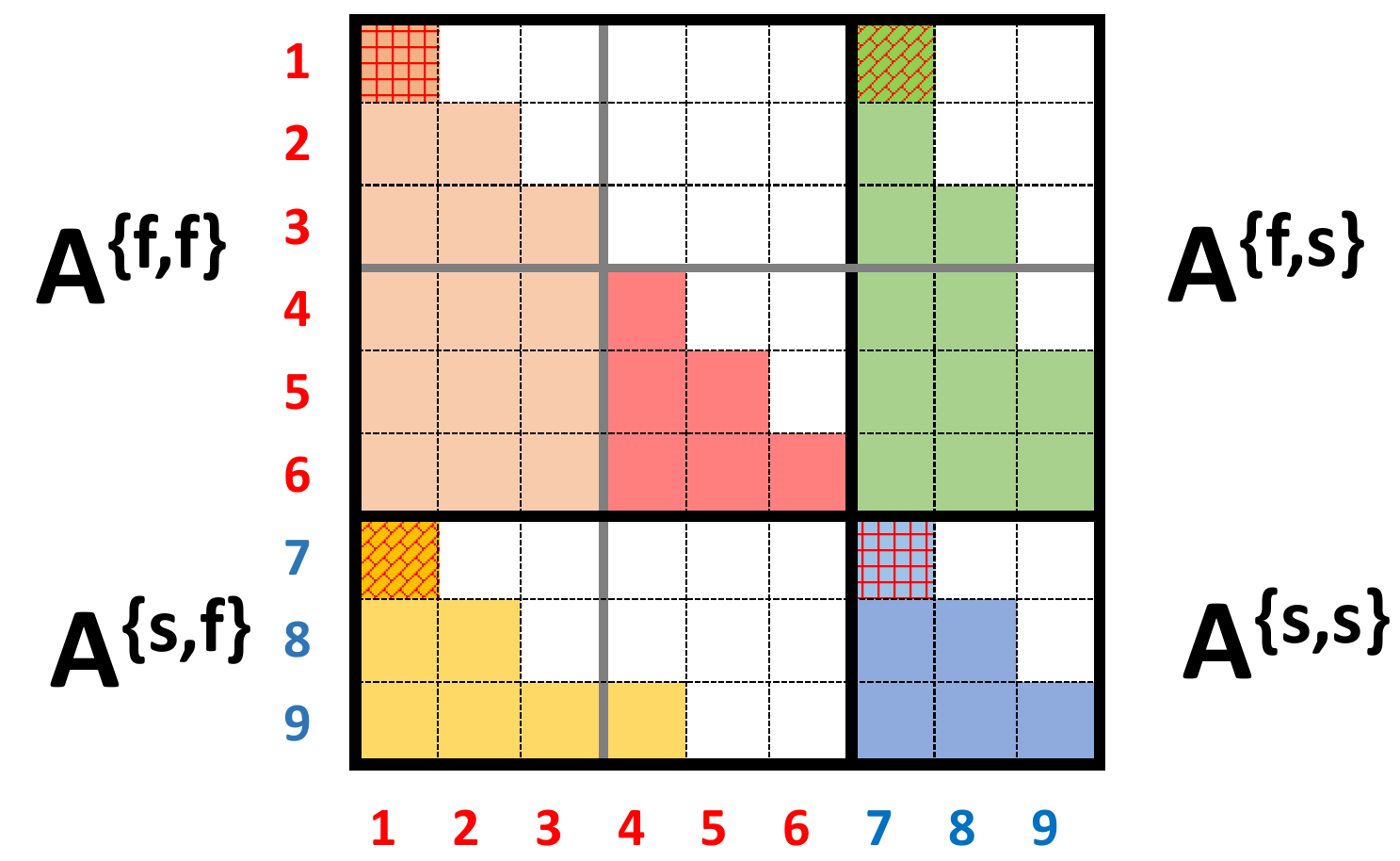}
	\caption{Butcher tableau of a coupled \mgark. The dashed entries of the coupling matrices violate the sparsity complementarity.}
	\label{fig:MGARK-butcher-coupled-standard}
\end{subfigure}
\hfill
\begin{subfigure}{0.475 \textwidth}
	\centering
	\includegraphics[width=\linewidth]{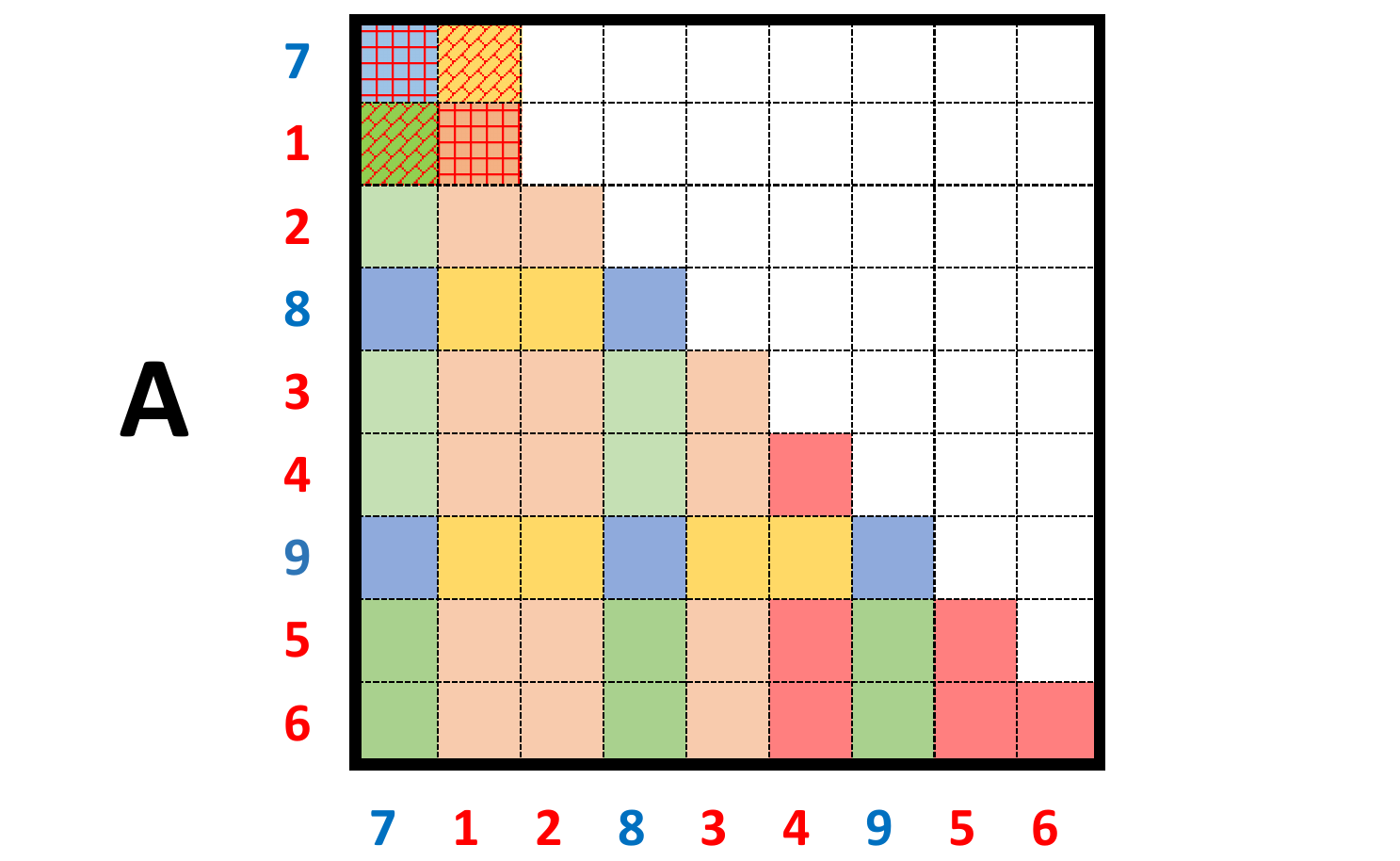}
	\caption{Permuted Butcher tableau of a coupled \mgark. The first slow and the first fast steps, dashed, need to computed together in a coupled manner.}
	\label{fig:MGARK-butcher-coupled-reordered}
\end{subfigure}

\caption{Example of decoupled and coupled \mgark with $M=2$ and $s=3$. Blue is the slow method, pink the first fast step, dark pink the second fast step, and green and yellow are the couplings.  The permuted versions of the tableaus reflect the sequential order of stage computations.  Note the entry above the diagonal for the coupled, permuted tableau due to the non-complementary coupling structure.}
\label{fig:MGARK-butcher-tableaus}
\end{figure}

From \cref{eqn:structure-Asf,eqn:structure-Afs} we see that the two coupling matrices $\A^{\{\f,\s\}}$ and $\A^{\{\s,\f\}}$ have related sparsity structures. For {\it decoupled \mgark methods} the sparsity structures are complementary, in the sense that one must have zeros in the entries where the other matrix can has nonzero elements:
\begin{equation}
\label{eq:decoupled-methods-equation}
	\A^{\{\s,\f\}} \sprod \tr{\A^{\{\f,\s\}}} = \boldsymbol{0}_{s^{\{\s\}} \times M s^{\{\f\}}},
\end{equation}
where $\sprod$ denotes element-by-element multiplication. This is schematically illustrated in \cref{fig:MGARK-butcher-decoupled-standard}.

For {\it coupled \mgark methods} the sparsity structures can overlap, in the sense that they both can have non-zeros entries in the same location:
\begin{equation}
	\A^{\{\s,\f\}} \sprod \tr{\A^{\{\f,\s\}}} \ne \boldsymbol{0}_{s^{\{\s\}} \times M s^{\{\f\}}}.
\end{equation}
This is illustrated in \cref{fig:MGARK-butcher-coupled-standard}. The overlapping non-zero coupling coefficients, indicated by dashed boxes in \cref{fig:MGARK-butcher-coupled-standard}, imply that the corresponding fast and slow stages need to be computed together, in a step that involves the entire non-partitioned system. Coupled methods are not pursued further in this paper.

\subsection{Order of the slow and fast stage evaluation}

The sparsity structure of the coupling matrices $\A^{\{\s,\f\}}$ and  $\A^{\{\f,\s\}}$ determines the order in which the fast and slow stages can be evaluated such as to respect the data dependencies between them. Using the notation defined in \cref{subsec:slow-to-fast-coup} and \cref{subsec:fast-to-slow-coup}, the general order in which stages are evaluated is follows:
\begin{enumerate}
\item Start by evaluating the fast stages $i$ for which $a^{\{\s,\f\}}_{1,i} \ne 0$; these are the fast stages needed by the computation of the first slow stage $c_1^{\{\s\}}$. If the entire first row $a^{\{\s,\f\}}_{1,:} = 0$ then proceed with evaluating the first slow stage $c_1^{\{\s\}}$.
\item After stage $I_j$ of the $L_j$-th micro-step is computed, the $j$-th stage $c_j^{\{\s\}}$ of the slow method is evaluated.
\item Continue with evaluating stages of the fast method, and stop after stage $I_{j+1}$ of the $L_{j+1}$-st micro-step to compute stage $c_{j+1}^{\{\s\}}$ of the slow method.
\item Continue until all fast and slow stages are evaluated.
\end{enumerate}
%
\begin{remark}[Time ordering]
Assuming that the slow  stage abscissae are increasing, and that fast stage abscissae are non-decreasing in time:
\begin{equation}
	c_1^{\{\s\}} < \cdots < c_{s^{\{\s\}}}^{\{\s\}}, \qquad
	c_1^{\{\f\}} \le \cdots \le c_{s^{\{\f\}}}^{\{\f\}},
\end{equation}
a natural order to evaluate the \mgark stages follows the time ordering of their abscissae. Note that the $j$-th slow stage approximates the solution at time  $T_{M s^{\{\f\}}+j} = t_n+c_{j}^{\{\s\}}\,H$, while the $j$-th fast stage of the $L$-th microstep approximates the solution at time  $T_{(L-1) s^{\{\f\}}+j} =t_n+(L -1 + c_{j}^{\{\f\}})\,(H/M)$.
The pairs $(L_j,I_j)$ are chosen such that the stages are evaluated in increasing order of their approximation times:
\begin{subequations}
\begin{align}
	\frac{L_j -1 + c_{I_j}^{\{\f\}}}{M} \le c_j^{\{\s\}} \le \frac{L_j -1 + c_{I_j+1}^{\{\f\}}}{M} \quad
	&\textnormal{for}\quad 1 \le I_j < s^{\{\f\}}, \\
	\frac{L_j + c_{s^{\{\f\}}}^{\{\f\}}}{M} \le c_j^{\{\s\}} \le \frac{L_j + 1 + c_{1}^{\{\f\}}}{M} \quad
	&\textnormal{for}\quad I_j = s^{\{\f\}}, \\
\nonumber
\textnormal{or}	\quad\frac{L_j -1 + c_{s^{\{\f\}}}^{\{\f\}}}{M} \le c_j^{\{\s\}} \le \frac{L_j + c_{1}^{\{\f\}}}{M} \quad 
	&\textnormal{for}\quad I_j =0.
\end{align}
\end{subequations}
\end{remark}
%

\begin{remark}[Simpler time ordering]
It is possible to simplify the time ordering such that the slow stages are evaluated at the end of full micro-steps.
The stage $c_j^{\{\s\}}$ is evaluated after the end of micro-step $\lambda-1$, but before micro-step $\lambda$, if
$(\lambda -1)/M \le c_j^{\{\s\}}  < \lambda/M$, in which case $L_j := \lambda$ and $I_j=0$. When $c_j^{\{\s\}} \ge 1$ we take $L_j=M, I_j=0$.
\end{remark}

\subsection{Reordering the GARK Butcher tableau}

In the standard form of the GARK Butcher tableau \cref{eqn:mrRK-butcher} the first $M s^{\{\f\}}$ stages are for the fast method, and stages $M s^{\{\f\}}+1$ to $M s^{\{\f\}} + s^{\{\s\}}$ are for the slow method. 
%
%
\begin{figure}
	\centering
	\includegraphics*[width = 0.7 \textwidth]{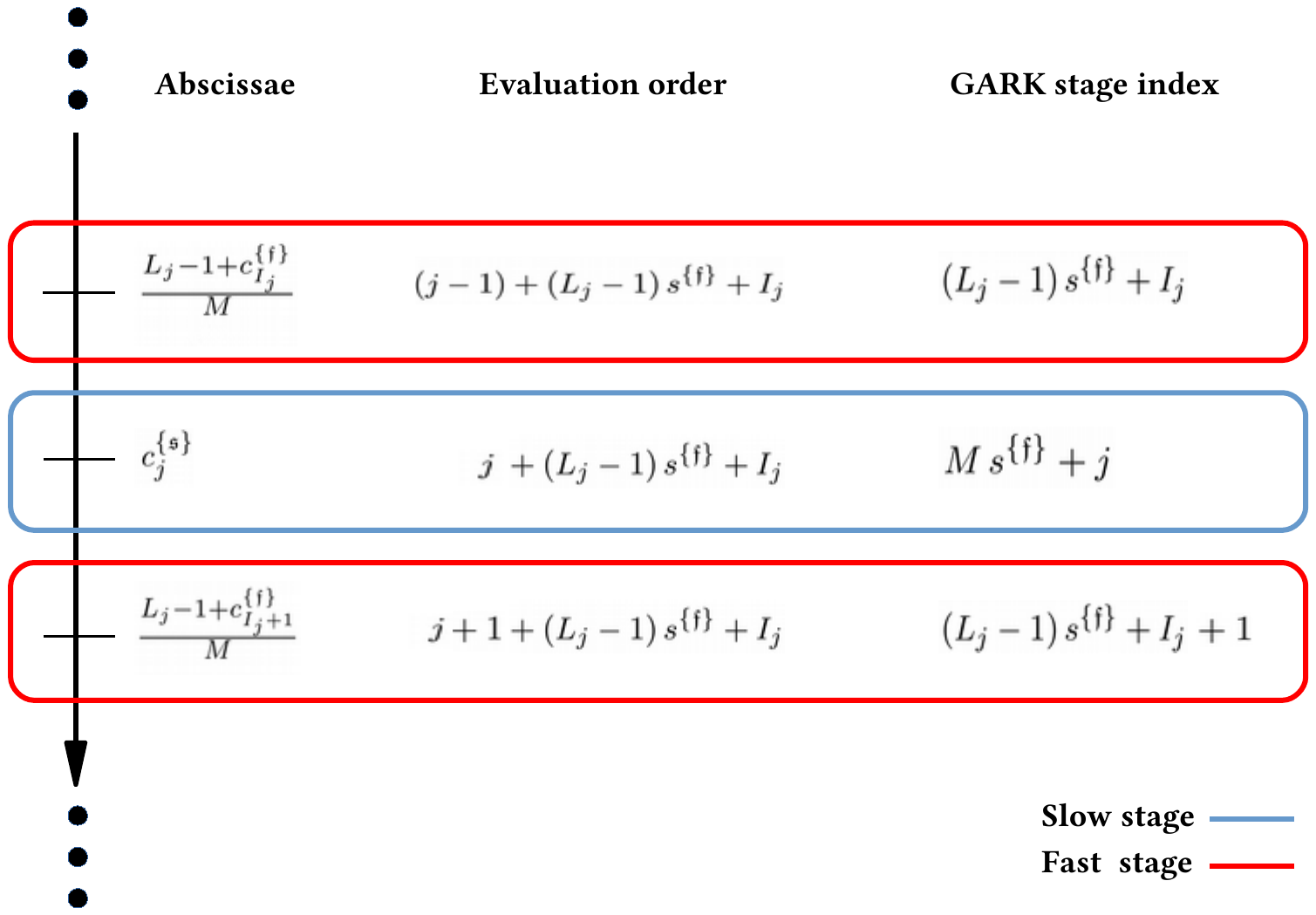}
	\caption{An example of the order of computing stages in \mgark with stage sequence $\{(L_j, I_j)\}$. }
	\label{fig:stage_orders}
\end{figure}
Let $\textsc{ic}$ be the vector of GARK tableau stage indices sorted in the order of computations, as summarized in  \cref{fig:stage_orders}. A renumbering of stages leads to a row and column permutation of the Butcher tableau \cref{eqn:mrRK-butcher}. The reordered Butcher matrix is $\A_\textsc{gark}(\textsc{ic},\textsc{ic})$. The reordering is illustrated in \cref{fig:MGARK-butcher-tableaus}. We distinguish the following cases:
\begin{itemize}
\item If the reordered Butcher tableau is strictly lower triangular then the \mgark method is explicit, and each stage uses only previously computed information. 
\item If the reordered Butcher tableau is lower triangular, with some non-zero diagonal entries, then the \mgark method is implicit, but decoupled. The non-zero diagonal entries correspond to implicit fast or slow stages. This case is illustrated in \cref{fig:MGARK-butcher-decoupled-reordered}. 
\item Finally, if the reordered Butcher tableau has block-diagonal entries then the \mgark method is coupled, in the sense that several fast and slow stages need to be solved together, in a coupled manner. This case is illustrated in \cref{fig:MGARK-butcher-coupled-reordered}.
\end{itemize}
\iftechreport
As an example let us consider the GARK Butcher matrix for the explicit method \method{\explicit}{\explicit}{2}{1}{2}{2}{A} introduced in \cref{subsec:method_EXEX2} for $M=3$:
%
\bgroup
\def\arraystretch{1.75}
\begin{tiny}
\begin{align*}
\A &= \left(\begin{array}{cccccc|cc} 0 & 0 & 0 & 0 & 0 & 0 & 0 & 0\\ \frac{2}{9} & 0 & 0 & 0 & 0 & 0 & \frac{2}{9} & 0\\ \frac{1}{4} & \frac{3}{4} & 0 & 0 & 0 & 0 & \frac{11}{60} & \frac{3}{20}\\ \frac{1}{4} & \frac{3}{4} & \frac{2}{9} & 0 & 0 & 0 & - \frac{19}{180} & \frac{9}{20}\\ \frac{1}{4} & \frac{3}{4} & \frac{1}{4} & \frac{3}{4} & 0 & 0 & \frac{31}{60} & \frac{3}{20}\\ \frac{1}{4} & \frac{3}{4} & \frac{1}{4} & \frac{3}{4} & \frac{2}{9} & 0 & - \frac{79}{180} & \frac{9}{20}\\ \hline
0 & 0 & 0 & 0 & 0 & 0 & 0 & 0\\ 2 & 0 & 0 & 0 & 0 & 0 & \frac{2}{3} & 0 \end{array}\right),
\quad
\A(\textsc{ic},\textsc{ic})=\left(\begin{array}{cccccccc} 0 & 0 & 0 & 0 & 0 & 0 & 0 & 0\\ 0 & 0 & 0 & 0 & 0 & 0 & 0 & 0\\ \frac{2}{9} & \frac{2}{9} & 0 & 0 & 0 & 0 & 0 & 0\\ \frac{2}{3} & 2 & 0 & 0 & 0 & 0 & 0 & 0\\ \frac{11}{60} & \frac{1}{4} & \frac{3}{4} & \frac{3}{20} & 0 & 0 & 0 & 0\\ - \frac{19}{180} & \frac{1}{4} & \frac{3}{4} & \frac{9}{20} & \frac{2}{9} & 0 & 0 & 0\\ \frac{31}{60} & \frac{1}{4} & \frac{3}{4} & \frac{3}{20} & \frac{1}{4} & \frac{3}{4} & 0 & 0\\ - \frac{79}{180} & \frac{1}{4} & \frac{3}{4} & \frac{9}{20} & \frac{1}{4} & \frac{3}{4} & \frac{2}{9} & 0 \end{array}\right),
\end{align*}
\end{tiny}
\egroup
where $\textsc{ic} = [7,1,2,8,3,4,5,6]$.
%
The permuted Butcher tableau verifies the explicit nature of the method, and also provides a progression of stage computations that is practically easy to implement.
\fi
\section{Design of practical decoupled \mgark methods}
\label{sec:design}

In this section we discuss several desirable properties of practical \mgark methods, and the design methodologies to incorporate them in the construction of new high order schemes.

\subsection{Design principles}
\label{sec:design-principles}
The following set of design principles incorporates properties that are important for ensuring the practicality of \mgark schemes:
\begin{enumerate}

\item Practical \mgark methods need to be {\it high order} (have an order of accuracy larger than two), and {\it generic} in the multirate ratio (have the same order of accuracy for any $M$.) This typically requires that the coupling coefficients are functions of $M$.  We have addressed these requirements via the order conditions derived in \cref{sec:order-conditions}.

\item Methods where the fast and slow base schemes are both either explicit or implicit should be {\it telescopic} \cref{eqn:mr-telescopic-conditions} in order to be easily applicable to multi-partitioned systems with multiple time scales. All explicit-explicit methods proposed in this paper are telescopic.

\item In order to maintain computational efficiency, a complex coupling between multiple fast and slow stages needs to be avoided. In this paper, we focus on {\it decoupled} multirate methods satisfying equation \cref{eq:decoupled-methods-equation}, as they are far less computationally demanding than coupled methods. However, the overall stability of the scheme may be affected by decoupling.

\item Methods should be optimized for general-purpose time integration, i.e. have small principal error terms and large stability regions.  Note that the principal error for \mgark methods is a function of $M$.  If not properly controlled, the coupling errors can grow with $M$, thereby reducing the overall efficiency of the method. All schemes derived herein have coefficients optimized for small principal error and large stability region.

\item Consider the case where one of the base methods $q \in \{\f,\s\}$ is implicit.  A favorable property of the entire \mgark method is stiff accuracy \cite{Sandu_2016_GARK-MR}:
\begin{equation}
\label{eqn:stiff-accuracy}
	\mathbf{e}_{s^{\{q\}}}\trsym \mprod \A^{\{q, \f\}} = \tr{\b^{\{\f\}}}, \qquad 
	\mathbf{e}_{s^{\{q\}}}\trsym \mprod \A^{\{q, \s\}} = \tr{\b^{\{\s\}}}.
\end{equation}
This property implies that the base implicit Runge--Kutta method is stiffly accurate, and that the \mgark stability function approaches zero as the $q$-th component of the system becomes infinitely stiff \cite{Sandu_2016_GARK-MR}.  All schemes derived herein having an implicit base method enjoy the property \cref{eqn:stiff-accuracy}.

\item Practical use of \mgark methods demands an error control mechanism that is capable of adapting both the step size $H$ and the multirate ratio $M$. Therefore the error control problem in multirate integration is fundamentally more complex than in traditional, single rate integration. We fully address the error control issue in \cref{sec:step_control}. 

\item It is desirable to derive multirate methods with reduced coupling errors by requiring that both the main and the embedded methods satisfy higher order coupling conditions. This strategy isolates the dominant local truncation errors to the solution of the slow and fast components and greatly simplifies the task of adaptively choosing $H$ and $M$, as discussed in \cref{sec:step_control}. 
\begin{definition}[Naturally adaptive schemes]
\label{def:naturally-adaptive}
A \mgark scheme of order $p$ is {\it naturally adaptive} if the fast and slow local truncation error terms (corresponding to N-trees with nodes of the same color) are $\mathcal{O}(h^{p+1})$, and the coupling local truncation error terms (corresponding to N-trees with nodes of different colors) are $\mathcal{O}(h^{p+2})$.
\end{definition}
We construct {\it naturally adaptive} second and third order schemes in \cref{sec:new_mr_schemes}.

\item Methods of type S are optimized for simplicity and stability. First, by design, one seeks to maximize the sparsity of the coupling coefficient matrices, such as to simplify the implementation and decrease data dependencies. Next, the coupling coefficients can potentially become large in magnitude when $M$ increases, and this can lead to large cancellation errors in the context of finite precision arithmetic, as well as to a degradation of numerical stability.  It is desirable that the coupling coefficients are optimized such that their magnitude remains bounded for large values $M$. 

\end{enumerate}

\subsection{Design process}
\label{sec:design-process}

The order conditions \cref{eq:mgark-high-order-conditions} and the additional constraints associated with the design principles of \cref{sec:design-principles} lead to large, nonlinear systems of equations that need to be solved for the method coefficients. The resulting coupling coefficients $\A^{\{\s,\f\}}$ and $\A^{\{\f,\s\}}$ are {\it functions} of the micro-step number $\lambda$ and the multi-rate step size ratio $M$.  Our proposed approach for designing \mgark methods of order $p$ is a three-step  process, as follows.

\begin{enumerate}

\item The first step is to construct optimized base slow $\bigl(A^{\{\s,\s\}},b^{\{\s\}}\bigr)$ and fast $\bigl(A^{\{\f,\f\}},b^{\{\f\}}\bigr)$ schemes of order $p$.  Implicit base methods are selected to be stiffly-accurate SDIRK schemes.  Explicit base methods are chosen to have small principal error components. When both the slow and the fast methods are of the same type (explicit or implicit) we choose the same discretization coefficients $A^{\{\s,\s\}}=A^{\{\f,\f\}}$, $b^{\{\s\}}=b^{\{\f\}}$ such as to obtain a telescopic \mgark method.

\item The second step is to define the sparsity patterns of the coupling matrices $\A^{\{\s,\f\}}$ and $\A^{\{\f,\s\}}$. These patterns determine the computational flow of the method and its implementation complexity, and influence the overall stability and accuracy properties. We manually test various sparsity patterns to balance all of these properties.

\item The third step computes the coupling coefficients $\A^{\{\s,\f\}}$ and $\A^{\{\f,\s\}}$ such as to satisfy the coupling conditions \cref{eq:mgark-high-order-conditions} up to order $p$.  Any free parameters in the family are used to minimize the Euclidean norm of the residuals of the order $p + 1$ coupling conditions; for naturally adaptive \mgark methods all these residuals are cancelled. In this work the solution of order conditions and the minimization of the error coefficients are carried out with Mathematica Version 11.2.

\end{enumerate}

An alternative to this derivation strategy is a monolithic constrained optimization procedure that minimizes the residuals of the order $p + 1$ coupling conditions, subject to solving the conditions \cref{eq:mgark-high-order-conditions} up to order $p$, and subject to structural constraints such as decoupling \eqref{eq:decoupled-methods-equation} and stiff accuracy \eqref{eqn:stiff-accuracy}.  

The proposed three-step procedure is preferable for two reasons.  First, we expect \mgark methods to be applied to problems with a rather weak coupling between partitions, where the primary sources of error are the base methods and not the coupling.  Our approach gives precedence to the base errors first.  Second, our procedure is practical as it reduces the number of nonlinear equations that need to be solved together during the design.
%


\section{Error estimation and adaptive \mgark methods}
\label{sec:step_control}
Adaptivity of traditional (single rate) methods adjusts the step size such as to ensure the desired accuracy of the solution at a minimal computational effort. In the context of Runge--Kutta schemes a second ``embedded'' method is used to provide an aposteriori estimate of the local truncation error \cite[CH II.4]{Hairer_book_I}. The step size is adjusted to bring the local error estimate to the user prescribed level using the asymptotic relation that  this error is proportional to $\propto H^{p+1}$.

Adaptivity of multirate methods is more complex, as there are two independent parameters that control the solution accuracy and efficiency: the macro-step size $H$ and the multirate ratio $M$ (or, equivalently, the macro-step $H$ and the micro-step $h$.) 
In order to achieve adaptivity of both $H$ and $M$ we propose to construct error estimates using not one, but several embedded methods, such as to obtain additional information about the structure of the local truncation error. Analytical asymptotic formulas are used to quantify the dependency of the local error terms on both $H$ and $M$. 
Armed with this information, we develop several criteria for adaptively selecting both the macro-step size and the multirate ratio.


\subsection{Local error structure for a second order \mgark method}
\label{sec:multirate_errors_3}

%
In order to understand the structure of the local truncation error we first consider a second order, internally consistent \mgark scheme. The principal error terms of order $\bigo{H^3}$ are associated with the third order conditions, and are provided in \cref{tab:imex_err}.

\begin{table}[ht!]
	\centering
	\begin{tabular}{|c|c|c|c|c|} \hline
		Error & Error & Elementary & Order condition & Residuals for   \\
		 type & term & differential &  & \method{\implicit}{\explicit}{2}{1}{2}{2}{A} \\ \hline
		Slow & $e^{\{\s\}}_{3,1}$ & $f^{\{\s\}}_{y,y}\left( f^{\{\s\}},f^{\{\s\}} \right)$ & $\tr{\b^{\{\s\}}}\mprod \c^{\{\s\} \times 2} = \frac{1}{3}$ & $0$ \\
		Slow & $e^{\{\s\}}_{3,2}$ & $f^{\{\s\}}_y f^{\{\s\}}_y f^{\{\s\}}$ & $\tr{\b^{\{\s\}}} \mprod \A^{\{\s,\s\}} \mprod \c^{\{\s\}} = \frac{1}{6}$ & $\frac{1}{6}$ \\
		Fast & $e^{\{\f\}}_{3,1}$ & $f^{\{\f\}}_{y,y}\left( f^{\{\f\}},f^{\{\f\}} \right)$ & $\tr{\b^{\{\f\}}}\mprod \c^{\{\f\} \times 2} = \frac{1}{3}$ & $\frac{4 - 3 \sqrt{2}}{12 M^2}$ \\
		Fast & $e^{\{\f\}}_{3,2}$ & $f^{\{\f\}}_y f^{\{\f\}}_y f^{\{\f\}}$ & $\tr{\b^{\{\f\}}} \mprod \A^{\{\f,\f\}} \mprod \c^{\{\f\}} = \frac{1}{6}$ & $\frac{4 - 3 \sqrt{2}}{6 M^2}$ \\
		Coupling & $e^{\{\coup\}}_{3,1}$ & $f^{\{\s\}}_y f^{\{\f\}}_y f^{\{\f\}}$ & $\tr{\b^{\{\s\}}} \mprod \A^{\{\s,\f\}} \mprod \c^{\{\f\}} = \frac{1}{6}$ & $\frac{3 \sqrt{2}-3-M}{12 M}$ \\
		Coupling & $e^{\{\coup\}}_{3,2}$ & $f^{\{\f\}}_y f^{\{\s\}}_y f^{\{\s\}}$ & $\tr{\b^{\{\f\}}} \mprod \A^{\{\f,\s\}} \mprod \c^{\{\s\}} = \frac{1}{6}$ & $\frac{1}{6}$ \\
		\hline
	\end{tabular}
	\caption{Principal error terms of $\bigo{H^3}$ for second order \mgark schemes.}
	\label{tab:imex_err}
\end{table}
As an example, consider the second order \mgark method \method{\implicit}{\explicit}{2}{1}{2}{2}{A} from \cref{subsec:method_IMEX2}.
We note from \cref{tab:imex_err} that the  third order residuals for this method are all asymptotically bounded as $M$ increases,  and the local truncation error behaves like:
\begin{align}
\begin{split}
\lte &= \Big( 0 \cdot f^{\{\s\}}_{y,y}\bigl( f^{\{\s\}},f^{\{\s\}} \bigr) + \frac{1}{6} \cdot f^{\{\s\}}_y f^{\{\s\}}_y f^{\{\s\}}  \Big)\,H^3 \\
&\quad + \Big(\frac{4 - 3 \sqrt{2}}{12 M^2} \cdot f^{\{\f\}}_{y,y}\bigl( f^{\{\f\}},f^{\{\f\}} \bigr) + \frac{4 - 3 \sqrt{2}}{6 M^2} \cdot f^{\{\f\}}_y f^{\{\f\}}_y f^{\{\f\}}  \Big)\,H^3 \\
&\quad + \Big( \frac{3 \sqrt{2}-3-M}{12 M} \cdot f^{\{\s\}}_y f^{\{\f\}}_y f^{\{\f\}} + \frac{1}{6} \cdot f^{\{\f\}}_y f^{\{\s\}}_y f^{\{\s\}} \Big)\,H^3 + \bigo{H^4}, \\
\norm{\lte} &\approx \left( K_1 + K_2 \frac{1}{M} + K_3 \frac{1}{M^2} \right)\,H^3 + \bigo{H^4}.
\end{split}
\end{align}

\begin{figure}[h]
	\begin{subfigure}[t]{0.45\textwidth} 
		\includegraphics[width=\linewidth]{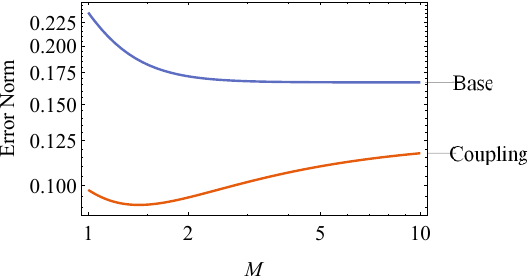}
		\caption{\method{\explicit}{\explicit}{2}{1}{2}{2}{A}.  The base error is $\bigo{H^3}$, but the coupling is $\bigo{H^4}$ since the method is naturally adaptive.}
		\label{subfig:EXEX21_error}
	\end{subfigure}%
	\hfill
	\begin{subfigure}[t]{0.45\textwidth} 
		\includegraphics[width=\linewidth]{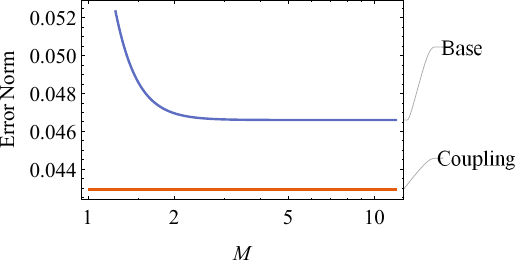}
		\caption{\method{\explicit}{\explicit}{3}{2}{3}{3}{A}.  Both the base and coupling errors are $\bigo{H^4}$.}
		\label{subfig:EXEX32_error}
	\end{subfigure}%
	\caption{Behavior of the local truncation error components for two \mgark schemes as the multirate ratio $M$ increases.}
	\label{fig:leading_error_coeff}
\end{figure}

In some cases, when there are sufficient degrees of freedom left after solving the order conditions, it is possible for a method of order $p$ to either cancel out coupling errors of order $\bigo{H^{p+1}}$ (\cref{subfig:EXEX21_error}) or to minimize them to assure better coupling behavior (\cref{subfig:EXEX32_error}).  \cref{fig:leading_error_coeff} shows how leading error coefficients change with $M$ for two optimized \mgark methods. 

\subsection{General structure of the \mgark local truncation error}
\label{sec:mgark-error-structure}
%
As the analysis in \cref{sec:multirate_errors_3} reveals, the local truncation error of the \mgark method  \cref{eqn:GARK-MR} has three components associated with the slow integration, the fast integration, and with the coupling:
\begin{equation}
	\lte = \lte^{\{\s\}} +  \lte^{\{\f\}} +  \lte^{\{\coup\}}.
\end{equation}
For a method of order $p$ the slow truncation error is that of applying one step with the base slow Runge--Kutta method with a step size $H$:  
\begin{subequations}
\label{eq:leading_error_terms}
\begin{equation}
\label{eq:lte_s}
		\norm{\lte^{\{\s\}}} \approx C^{\{\s\}} \cdot H^{p+1} + \bigo{H^{p+2}}.
\end{equation}
The fast truncation error is that of applying $M$ consecutive steps with the base fast Runge--Kutta method with a step size $H/M$. We use the global error estimate \cite[CH II.3]{Hairer_book_I}:
\begin{align*}
		\norm{\lte^{\{\f\}}} &\leq 
		\left(\frac{H}{M} \right)^p \frac{C'}{L} \big( \exp{\left(L\,H\right)}-1 \big)
		= \left(\frac{H}{M} \right)^p \frac{C'}{L} \big( L H + \bigo{H^2} \big),
\end{align*}
where $\norm{\partial f^{\{\f\}}/\partial y} \leq L$ and $C'$ is a constant, to obtain: 
\begin{equation}
\label{eq:lte_f}
\norm{\lte^{\{\f\}}} \le \frac{C^{\{\f\}}}{M^p} \cdot H^{p+1} + \bigo{H^{p+2}}.
\end{equation}
The principal terms of the coupling errors correspond to N-trees having $p+1$ nodes of two colors. A tree with $k$ fast nodes and $p+1-k$ slow nodes corresponds to a residual term constructed from multiplying $k$ fast matrix blocks in the Butcher tableau \eqref{eqn:mrRK-butcher} and $p+1-k$ slow matrix blocks. Note that each of the fast matrix blocks ($\tr{{\b}^{\{\f\}}},\A^{\{\s,\f\}},\A^{\{\f,\f\}}$) in \eqref{eqn:mrRK-butcher} carries a scaling factor of $1/M$. Assuming that the coupling coefficients remain bounded for large $M$, a product of $k$ fast blocks is $\bigo{1/M^k}$. Since we have coupling trees with $p \ge k \ge 1$ fast nodes, the corresponding products of fast blocks have scalings ranging from $1/M$ to $1/M^p$. Residuals contain sums of $M$ elementary coupling blocks as seen in \eqref{eq:mgark-high-order-conditions-blk}; sums of $M$ fast blocks have scaling factors ranging from $1$ to $1/M^{p-1}$, but some of the sums can have a small, fixed number of terms. In addition, the magnitude of the coupling coefficients, as resulted from the solution of order conditions, can scale as $M^\ell$ (for some $\ell \ge 0$). Consequently, a generic expression for the local coupling error is: 
\begin{equation}
\label{eq:lte_c}
		\norm{\lte^{\{\coup\}}} \approx \left(\sum_{i=-\ell}^{p} C^{\{\coup\}}_i \, M^{-i} \right) \cdot H^{p+1} + \bigo{H^{p+2}},
\end{equation}
\end{subequations}
While the slow and fast errors are simple functions of $H$ and $M$, the dependency of the coupling error on $M$  is more difficult to describe even with access to the residuals of each of the $p + 1$ order conditions.  These residuals can be positive or negative. Moreover, the evolution of the coupling error as a function of $M$ depends not only on residuals, but also on elementary differentials, since the constants $C^{\{\coup\}}_i$ in \eqref{eq:lte_c} are linear combinations of products of method residuals and norms of elementary differentials. 
Therefore, in an adaptive method, the effect of changing $M$ on the coupling error is more difficult to quantify.

\subsection{Use of multiple embedded methods to estimate the local truncation error}

Following the traditional Runge--Kutta strategy, we look to use embedded methods to obtain estimates of the local truncation error.  Specifically, we design pairs of main weight vectors $(b^{\{\f\}}, b^{\{\s\}})$ and embedded weight vectors $(\widehat{b}^{\{\f\}}, \widehat{b}^{\{\s\}})$ that produce solutions of different orders when paired with $(\A^{\{\f,\f\}}$,$\A^{\{\f,\s\}}$,$\A^{\{\s,\f\}}$, $\A^{\{\s,\s\}})$. The choice of the embedded weights is made such that additional function evaluations are avoided.


In the traditional Runge--Kutta approach, a single embedded method is used to provide a single estimate of the error, and a step size scaling factor is computed (for the next step) such that the corresponding scaled error estimate satisfies the accuracy requirements at hand. In multirate integration the solution accuracy depends on two parameters, $H$ and $M$, that can be adjusted independently. We have seen in \cref{sec:mgark-error-structure} that changes in these parameters affect differently the slow, fast, and coupling components of the local truncation error.

Our proposed strategy to select $H$ and $M$ adaptively relies on multiple embedded methods to independently estimate different parts of the error. Specifically, consider an \mgark scheme that produces a main solution $y_{n+1}$ of order $p$, i.e., cancels all residuals for two-trees with up to $p$ nodes. We seek to build an embedded solution $\widehat{y}^{\{\s\}}_{n + 1}$ that cancels all residuals up to order $p-1$, as well as all fast and all coupling residuals of order $p$. Similarly, we seek to build a second embedded solution $\widehat{y}^{\{\f\}}_{n + 1}$ that cancels all residuals up to order $p-1$, as well as all slow and all coupling residuals of order $p$. The three components of the local truncation error can then be estimated as follows:
\begin{equation}
\label{eqn:multiple-error-estimates}
\begin{split}
	\lte_{n+1} &\approx y_{n + 1} - \widehat{y}_{n + 1}, \quad
	\lte_{n+1}^{\{\s\}} \approx y_{n + 1} - \widehat{y}^{\{\s\}}_{n + 1}, \quad
	\lte_{n+1}^{\{\f\}} \approx y_{n + 1} - \widehat{y}^{\{\f\}}_{n + 1}, \\
	\lte_{n+1}^{\{\c\}} &= \lte_{n+1} - \lte_{n+1}^{\{\s\}} - \lte_{n+1}^{\{\f\}} \approx \widehat{y}^{\{\s\}}_{n + 1} 
	   +  \widehat{y}^{\{\f\}}_{n + 1}  - \widehat{y}_{n + 1} -y_{n + 1}.
\end{split}	   
\end{equation}

The construction of three different embedded methods that capture precisely some components of the error can be very difficult to achieve for high order schemes. For this reason we now discuss a simplified error estimation strategy that utilizes only two embedded methods. Consider an \mgark scheme with primary weights $(b^{\{\f\}}, b^{\{\s\}})$ that produces a main solution $y_{n+1}$ of order $p$. The residuals corresponding to two-trees with up to $p$ nodes are zero. We construct three different embedded solutions, as follows:
\begin{itemize}
\item A generic embedded solution $\widehat{y}_{n+1}$ of order $p-1$, obtained with weights  $(\widehat{b}^{\{\f\}}, \widehat{b}^{\{\s\}})$, that captures all the error terms and is used to approximate the overall local truncation error. The residuals associated with two-trees with up to $p-1$ nodes are zero, while the residuals of two-trees with $p$ nodes can be nonzero.
\item An embedded solution $\widehat{y}^{\{\s\}}_{n+1}$, generated with the weights $(b^{\{\f\}}, \widehat{b}^{\{\s\}})$, captures all of the slow error, part of the coupling error, and none of the fast error. The weight vector $b^{\{\f\}}$ cancels the residuals for all two-trees with up to $p$ nodes with a fast-colored root. The weight vector $\widehat{b}^{\{\s\}}$ cancels the residuals for all two-trees with up to $p-1$ nodes with a slow-colored root. However, the residuals of the two-trees with $p$ nodes and a slow-colored root, and the associated $\bigo{H^p}$ error terms, can be nonzero. Trees with a slow-colored root correspond to either slow or to coupling trees, and the solution difference $y_{n+1} - \widehat{y}^{\{\s\}}_{n+1}$ captures the sum of the corresponding errors terms.
\item Similarly, an embedded solution $\widehat{y}^{\{\f\}}_{n+1}$, generated with $(\widehat{b}^{\{\f\}}, b^{\{\s\}})$, captures all of the fast error, part of the coupling error, and none of the slow error.
\end{itemize}
In general, these mixed embeddings do not exactly isolate the slow, fast, and coupling error, but they can serve as approximations in \eqref{eqn:multiple-error-estimates}.  Note that the resulting approximate $\lte^{\{\s\}}$ and $\lte^{\{\f\}}$ contain the parts of the coupling error corresponding to trees of order $p$ with slow-colored roots and with fast colored-roots, respectively. Consequently the  $\lte^{\{\coup\}}$ estimated by the solution difference in \eqref{eqn:multiple-error-estimates} is zero.

However, for naturally adaptive methods (\cref{def:naturally-adaptive}), the simplified error estimation strategy does separate the main components of the fast and slow errors, as explained next. For a naturally adaptive method of order $p$ with weights $(b^{\{\f\}}, b^{\{\s\}})$ the residuals corresponding to all two-trees with up to $p$ nodes are zero, and the coupling residuals for two-trees with $p+1$ nodes are also zero.  The non-zero residuals of order $p$ correspond to either fast or low trees with $p$ nodes.

A naturally adaptive embedded solution $\widehat{y}_{n+1}$ of order $p-1$, obtained with weights  $(\widehat{b}^{\{\f\}}, \widehat{b}^{\{\s\}})$, cancels all residuals of order up to $p-1$; it also cancels the coupling residuals for two-trees of up to $p$ nodes. Consequently, the difference $y_{n + 1} - \widehat{y}^{\{\s\}}_{n + 1}$ contains only $\bigo{H^p}$ terms corresponding to slow-colored trees, and $y_{n + 1} - \widehat{y}^{\{\f\}}_{n + 1}$ contains only $\bigo{H^p}$ terms corresponding to fast-colored trees. Consequently, naturally adaptive embedded methods \textit{do} isolate the slow and the fast errors.  This property allows to construct more accurate adaptivity mechanisms.

\subsection{Controlling errors by adapting both the macro-step and the micro-step}
\label{subsec:adaptivity}

Based on our understanding of the structure of errors in the \mgark framework, we now look into practical approaches to adaptivity of multirate methods.  The following error estimates are available via the set of embedded methods:
\begin{align*}
	\widehat{\varepsilon}_{n+1} &:= \normerr{ y_{n+1} - \widehat{y}_{n+1} }, \quad
	\widehat{\varepsilon}_{n+1}^{\{\s\}} := \normerr{ y_{n+1} - \widehat{y}^{\{\s\}}_{n+1} }, \quad
	\widehat{\varepsilon}_{n+1}^{\{\f\}} := \normerr{ y_{n+1} - \widehat{y}^{\{\f\}}_{n+1} },
\end{align*}
where the error is measured by the following relative error norm \cite[CH II.4]{Hairer_book_I}:
\begin{align*}
	\normerr{ x - y }  := \sqrt{\frac{1}{d} \sum_{i = 1}^d \left(\frac{x_i - y_i}{\text{AbsTol}_i + \text{RelTol}_i \cdot \max{(|(x_i|,|y_i|)}}\right)^2}.
\end{align*}
Several heuristic strategies that use these error estimates to adapt both $H$ and $M$ are discussed below.

\subsubsection{Balancing error strategy}
\label{subsec:BES_adapt}
In this approach, we first use the estimated total truncation error $\widehat{\varepsilon}_{n+1}$ to control the macro-step size using the traditional mechanism based on the asymptotic error behavior \cite[CH II.4]{Hairer_book_I}: 
\begin{equation*}
\widehat{\varepsilon}_{n+2} \le 1 \qquad \Rightarrow \qquad
	H_{\text{new}} = \textnormal{fac} \cdot H \cdot \left(\widehat{\varepsilon}_{n+1}\right)^{-\frac{1}{p}},
\end{equation*}
where $\textnormal{fac} < 1$ is a safety factor.  

Next, the multirate ratio $M$ is adjusted such that the estimated slow and fast error components over the next step are equal to each other, i.e., the slow and fast contributions to the error are balanced. Assuming $q = \min{(p, \hat p)}$ and using the asymptotic formulas \eqref{eq:leading_error_terms} we have that:
\begin{align} 
\label{eq:simple_adaptivity_derivation}
	\begin{split}
		\widehat{\varepsilon}_{n+2}^{\{\s\}} &= \widehat{\varepsilon}_{n+2}^{\{\f\}}, \\
		C^{\{\s\}} \cdot H_{\text{new}}^{q+1} &\approx \frac{C^{\{\f\}}}{M_{\text{new}}^{q}} \cdot H_{\text{new}}^{q+1}, \\
		\widehat{\varepsilon}_{n+1}^{\{\s\}} \cdot \frac{H_{\text{new}}^{q+1}}{H^{{q+1}}} &\approx \widehat{\varepsilon}_{n+1}^{\{\f\}} \cdot \frac{M^{q}}{M_{\text{new}}^{q}} \cdot \frac{H_{\text{new}}^{q+1}}{H^{q+1}}, 
	\end{split}
\qquad \Rightarrow \qquad			
M_{\text{new}} \approx M \cdot \left( \frac{\widehat{\varepsilon}_{n+1}^{\{\f\}}}{\widehat{\varepsilon}_{n+1}^{\{\s\}}} \right)^{\frac{1}{q}}.
\end{align}
In this strategy, as well as the others, $M$ can be rounded up, down, or to the nearest integer.  Also in practice, the re-scaling of $H$ and $M$ for the next step can be bounded up and down in order to avoid large jumps and oscillations.  This adaptivity strategy serves as a very simple heuristic.  However, the choices of the new $H$ and of the new $M$ are made independently of each other, and the approach does not account for how their interaction impacts the error; the only mechanism for controlling the coupling error is the change in $H$.

\subsubsection{Efficiency optimization strategy}
\label{subsec:POS_adapt}

This approach focuses on the important aspect of the overall cost of multirate integration. Evaluation of the slow and fast partitions can have very different computational costs in some applications.  Moreover, an implicit method (e.g., applied to solve the fast component) is likely to be much more expensive than an explicit method (e.g., applied to the slow component).  We require the adaptive selection of $H$ and $M$ to satisfy the error tolerance criteria at a minimal overall computational cost. 

Let $t^{\{\s\}}$  and $t^{\{\f\}}$ represent the computational costs of a slow macro-step and a fast micro-step, respectively. We define the computational efficiency of a step as the progress made during the step ($H$) divided by the total cost of executing step ($t^{\{\s\}} + M_{\text{new}}  t^{\{\f\}}$).

The new values of $H$ and $M$ are selected such as to achieve the desired accuracy while maximizing the computational efficiency. This requires solving the following constrained optimization problem to minimize the inverse of efficiency:
\begin{equation}
\begin{split}
\min_{H_{\text{new}}, M_{\text{new}}} & \frac{t^{\{\s\}} + M_{\text{new}}\,  t^{\{\f\}}}{H_{\text{new}}}, \\
\text{subject to} & ~\widehat{\varepsilon}_{n+2} = 1.
\end{split}
\label{eq:HM-Adapt_Opt}
\end{equation}
Expanding the constraint yields: 
\begin{equation*}
	1
	= \widehat{\varepsilon}_{n+2}
	\approx \widehat{\varepsilon}_{n+1}^{\{\s\}} \cdot \frac{H_{\text{new}}^{q+1}}{H^{q+1}} + \widehat{\varepsilon}_{n+1}^{\{\f\}} \cdot \frac{M^q}{M_{\text{new}}^q} \cdot \frac{H_{\text{new}}^{q+1}}{H^{q+1}},
\end{equation*}
where we have used the fact that, for naturally adaptive methods, the coupling component of the local truncation error is negligible in the asymptotic regime. The constraint equation can be solved explicitly for $H_{\text{new}}$:
\begin{equation} \label{eq:opt_h_new}
	H_{\text{new}} = H \cdot \left(\widehat{\varepsilon}_{n+1}^{\{\s\}} + \widehat{\varepsilon}_{n+1}^{\{\f\}} \cdot \frac{M^q}{M_{\text{new}}^q}\right)^{-\frac{1}{q+1}}.
\end{equation}
After eliminating the constraint the optimization problem \cref{eq:HM-Adapt_Opt} simplifies to:
\begin{equation}
\min_{M_{\text{new}}} \frac{t^{\{\s\}} + M_{\text{new}}\, t^{\{\f\}}}{H} \left(\widehat{\varepsilon}_{n+1}^{\{\s\}} + \widehat{\varepsilon}_{n+1}^{\{\f\}} \cdot \frac{M^{q}}{M_{\text{new}}^{q}}\right)^{\frac{1}{q+1}}.
\label{eq:EOS}
\end{equation}
Note this is an integer minimization problem. One can solve it as a continuous optimization problem, then round the result to an integer to find the optimal $M_{\text{new}}$. Afterwards $H_{\text{new}}$ is computed from \cref{eq:opt_h_new}.

\begin{remark}[Timing]
The CPU times $t^{\{\s\}}$  and $t^{\{\f\}}$ can be evaluated online by timing the slow macro-steps and the fast micro-steps during their execution. The algorithm adjusts automatically if these compute times vary during the application lifetime.
\end{remark}
\section{New high-order \mgark methods}
\label{sec:new_mgark}

Using the design process outlined above we construct several \mgark methods for use in practical applications. We use the naming convention {\it \method{FAST}{SLOW}{$p$}{$\widehat{p}$}{$f$}{$s$}{type}}, where $p$ is the method order, $\widehat{p}$ is the embedded order, $f$ is the number of stages in the fast base method, and $s$ is the number of stages in the slow base method. Each component method is either explicit or implicit: FAST, SLOW $\in$ \{EX,IM\}. We distinguish between methods of type A (optimized for accuracy and for better step size control) and methods of type S (optimized for simplicity and for stability), therefore  $type \in$ \{A,S\}. The newly developed methods are as follows:
\begin{itemize} 
\item  \method{\explicit}{\explicit}{2}{1}{2}{2}{A}, \method{\explicit}{\explicit}{2}{1}{2}{2}{S},  \method{\explicit}{\explicit}{3}{2}{3}{3}{A}, \method{\explicit}{\explicit}{3}{2}{4}{4}{A}, \method{\explicit}{\explicit}{3}{2}{3}{3}{S}, and \method{\explicit}{\explicit}{4}{3}{5}{5}{A}  are explicit-explicit, schemes of order two, three, three and four respectively.
\item \method{\explicit}{\implicit}{2}{1}{2}{2}{A}, \method{\explicit}{\implicit}{3}{2}{3}{3}{A}, and \method{\explicit}{\implicit}{4}{3}{6}{5}{A} are methods with an explicit fast part and an implicit slow part of order two, three, and four respectively.
\item  \method{\implicit}{\explicit}{2}{1}{2}{2}{A}, \method{\implicit}{\explicit}{3}{2}{3}{3}{A}, and \method{\implicit}{\explicit}{4}{3}{6}{5}{A} are methods with an implicit fast part and an explicit slow part of order two, three, and four, respectively.
\end{itemize}
The coefficients of these methods are given in \cref{sec:new_mr_schemes}.
\section{Numerical experiments}
\label{sec:numerics}

In this section we carry out numerical experiments to validate and test the newly derived \mgark methods.
\subsection{Additive partitioning tests}

The first experiment is carried out using a two-dimensional unsteady convection-diffusion equation \cite[Ch. 3]{elman2006finite} in a square spatial domain $\Omega = \{ 0 \leq x,y \leq 1\}$ and with a circular wind profile:
\begin{subequations}
	\label{eq:convection-diffusion}
	\begin{equation}
	\begin{split}
	u_t -\varepsilon \nabla^2 u + w \cdot \nabla u = 0 ~ \text{in} ~ \Omega, \quad
	u =  0  ~ \text{on}  ~ \partial \Omega, \quad w =
        \begin{bmatrix}
        \phantom{-} 2y(1-x^2)\\
 - 2x(1-y^2)
\end{bmatrix}.
	\end{split}
	\end{equation}
A Streamline Upwind Petrov-Galerkin (SUPG) spatial discretization is used, which leads to a semi-discrete system of linear ODEs: 
%
%
	\begin{align}
	\label{eq:CD_ODEFull}
	\mathbf{M}^h u^h_t &= \A\, u^h  + (\vec{n} + \vec{n}^{stab})\, u^h ,
	\end{align}
\end{subequations}
where $\mathbf{M}^h $ and $\A$ are mass and stiffness matrices and $\vec{n} + \vec{n}^{stab}$ represent linear forms of the convective term and SUPG stabilization. 

In the multirate experiments we designate the first term in \cref{eq:CD_ODEFull} as the slow component, $f^{\{\s\}} = (\mathbf{M}^h)^{-1} \,\A\, u^h$, and the second term as the fast component, $f^{\{\f\}} = (\mathbf{M}^h)^{-1} \,(\vec{n} + \vec{n}^{stab})\, u^h $. In practice the splitting choice is informed by inspecting the spectral radius of the right hand side operators. The weak form of the PDE and the corresponding \mgark schemes are implemented in the FEniCS package \cite{AlnaesBlechta2015a}, which is used to carry out the convergence experiments. 

\begin{figure}[h]
	   \begin{subfigure}[t]{0.45\textwidth} 
			\includegraphics[height=0.9 \linewidth]{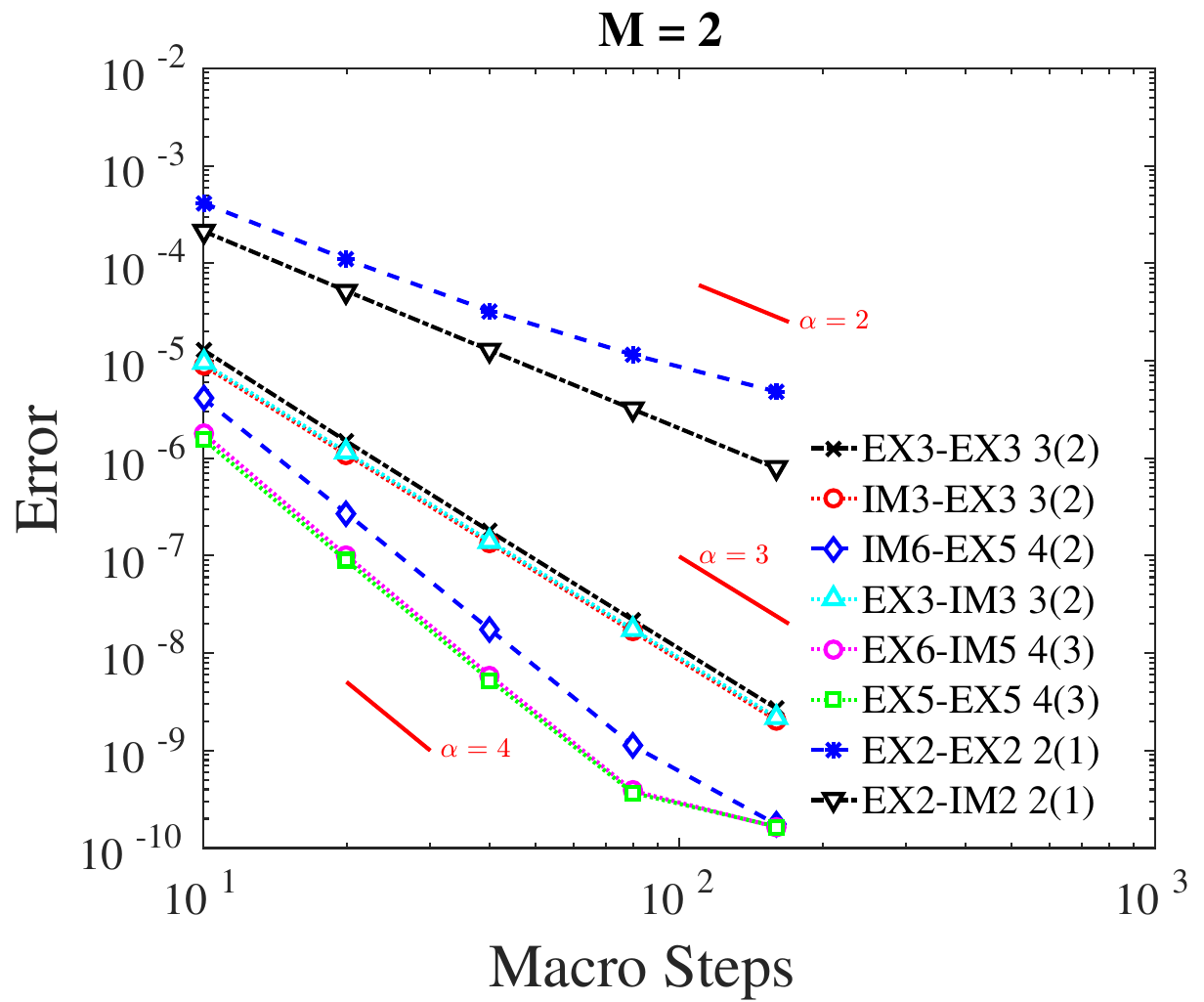}
	   \end{subfigure}%
      $\quad$
 	   \begin{subfigure}[t]{0.45\textwidth} 
      	\includegraphics[height=0.9 \linewidth]{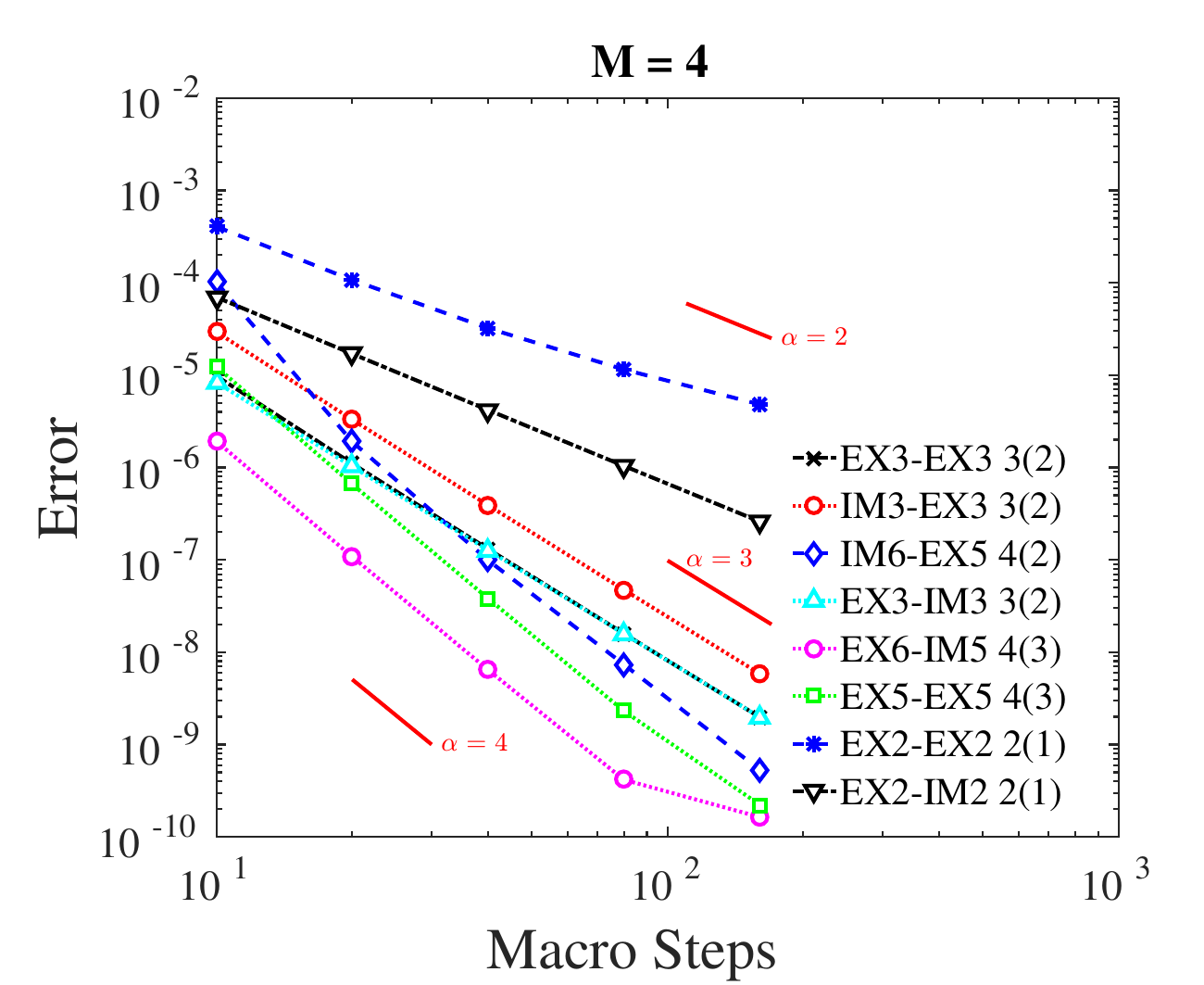}
       \end{subfigure}%
       \hfill
       \caption{Convergence plots for the unsteady convection-diffusion test \eqref{eq:convection-diffusion} over the time span $T = [0, 10]$ seconds. A fixed macro-step time integration is carried out with varying multirate step ratios $M$ using \mgark type A methods.}
       \label{fig:CD_convergence}
\end{figure}

The convergence diagrams for this test are shown in \cref{fig:CD_convergence}. The numerical orders of convergence for all schemes match their theoretical orders for the multirate step ratios tested. 
\iftechreport
\Cref{fig:CD_snapshots} shows the evolution of the model solution over time.
\begin{figure}[h]
	\begin{subfigure}[t]{0.33\textwidth} 
		\includegraphics[width=0.9\linewidth]{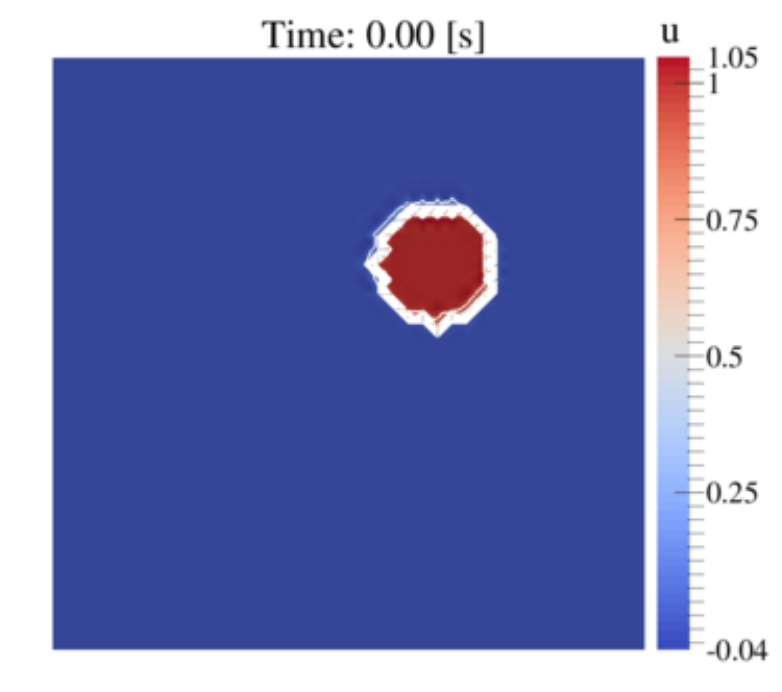}
	\end{subfigure}%
	\hfill
	\iftechreport
	\begin{subfigure}[t]{0.33\textwidth} 
		\includegraphics[width=0.9\linewidth]{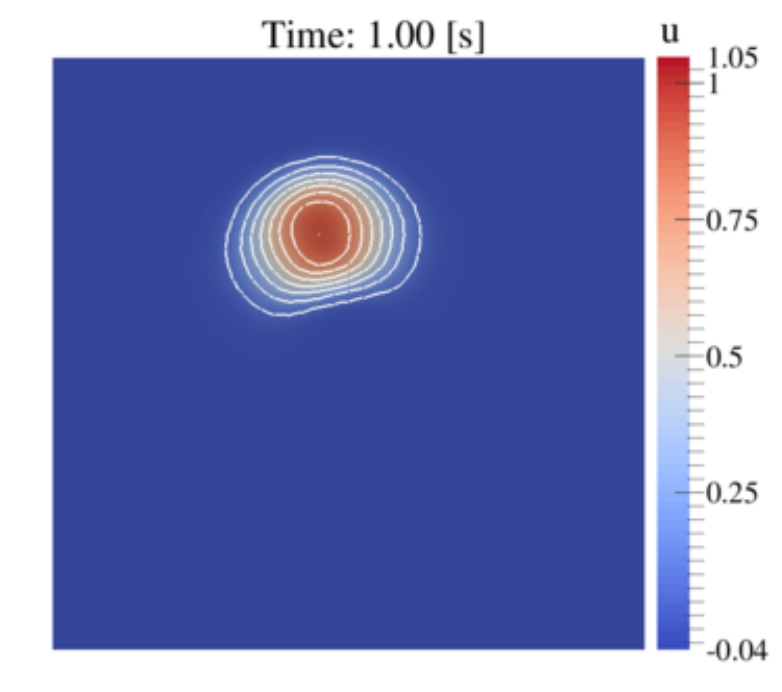}
	\end{subfigure}%
	\hfill
	\fi 
	\begin{subfigure}[t]{0.33\textwidth} 
		\includegraphics[width=0.9\linewidth]{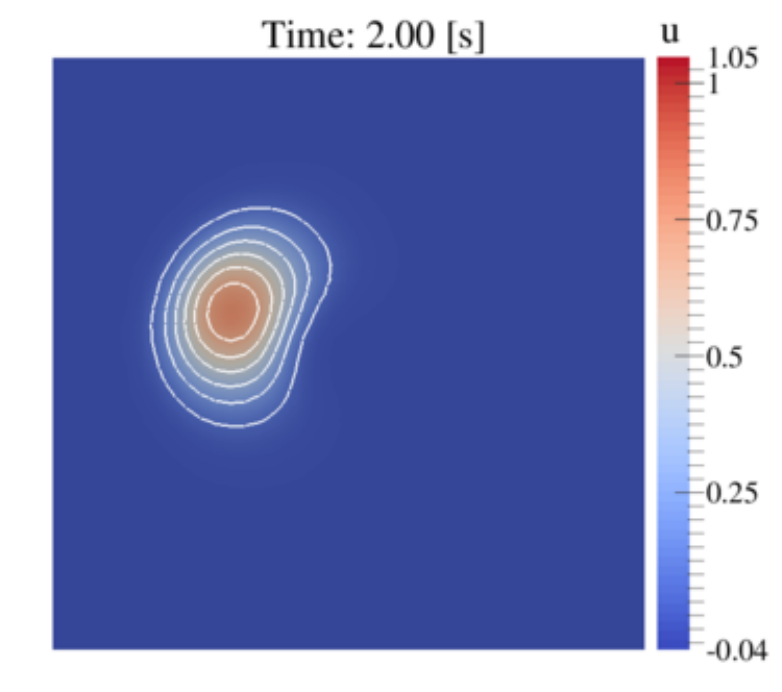}
	\end{subfigure}%

	\begin{subfigure}[t]{0.33\textwidth} 
		\includegraphics[width=0.9\linewidth]{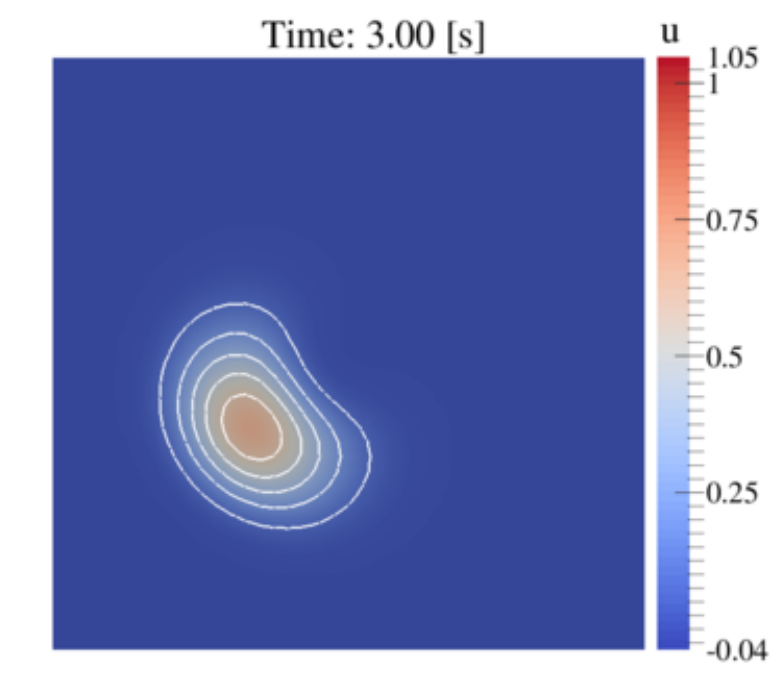}
	\end{subfigure}%
	\hfill  
		\begin{subfigure}[t]{0.33\textwidth} 
		\includegraphics[width=0.9\linewidth]{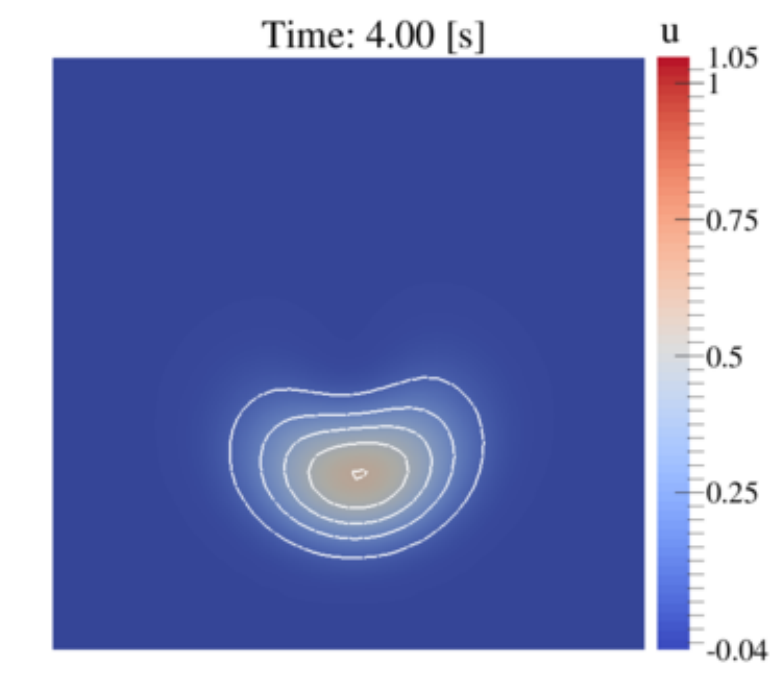}
	\end{subfigure}%
	\hfill
		\begin{subfigure}[t]{0.33\textwidth} 
		\includegraphics[width=0.9\linewidth]{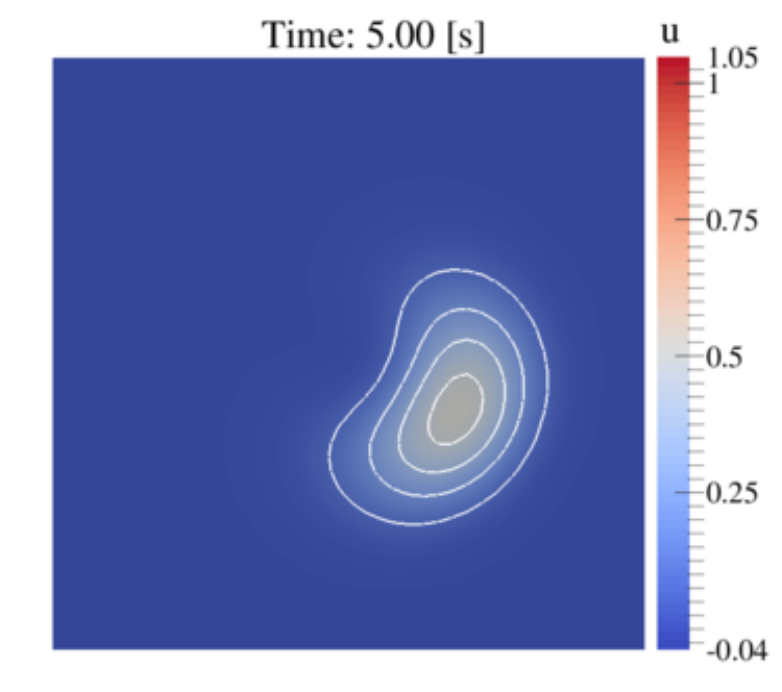}
	\end{subfigure}%
	\hfill
	\caption{Evolution of the unsteady convection-diffusion problem \eqref{eq:convection-diffusion} solution in time.}
	\label{fig:CD_snapshots}
\end{figure}
\fi 
%

\subsection{Timing experiments}
Timing experiments are performed in MATLAB using the Gray-Scott model \cite{lee1993pattern}. Here, we are interested in the application of \mgark methods to additive splittings of the right hand side into  linear and nonlinear terms. This reaction-diffusion PDE is:
\begin{equation}
\label{eq:Gray-Scott}
\left\{
\begin{split}
		u_t &= \nabla \cdot ( \varepsilon_u\, \nabla u ) - u\,v^2 + \mathfrak{f} (1-u), \\
		v_t &=  \nabla \cdot ( \varepsilon_v\, \nabla v ) + u\,v^2 - (\mathfrak{f} + \mathfrak{k})\, v.
\end{split}
\right.
\end{equation}
The domain is the unit square discretized with second order finite differences. The reaction parameters are $\mathfrak{k}= 0.0520$ and $\mathfrak{f}= 0.0180$. 

In the following experiments the nonlinear reaction terms on the right hand side are considered the fast system, and the diffusion terms are regarded as the slow one. 

\begin{figure}[h]
	\begin{subfigure}[t]{0.33\textwidth} 
		\includegraphics[width=0.9\linewidth]{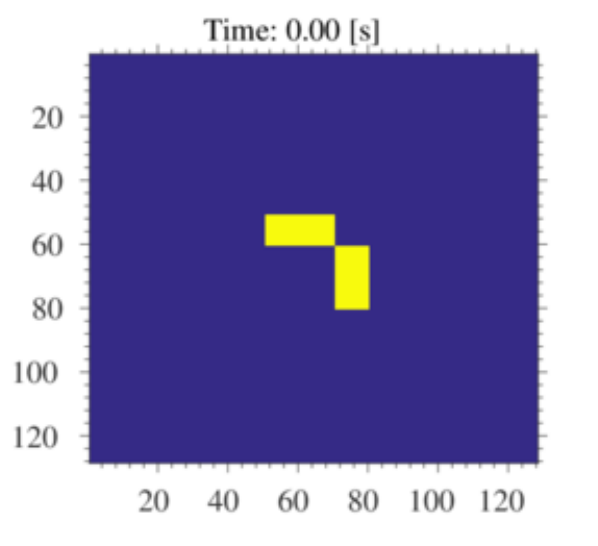}
	\end{subfigure}%
	\hfill
	\begin{subfigure}[t]{0.33\textwidth} 
		\includegraphics[width=0.9\linewidth]{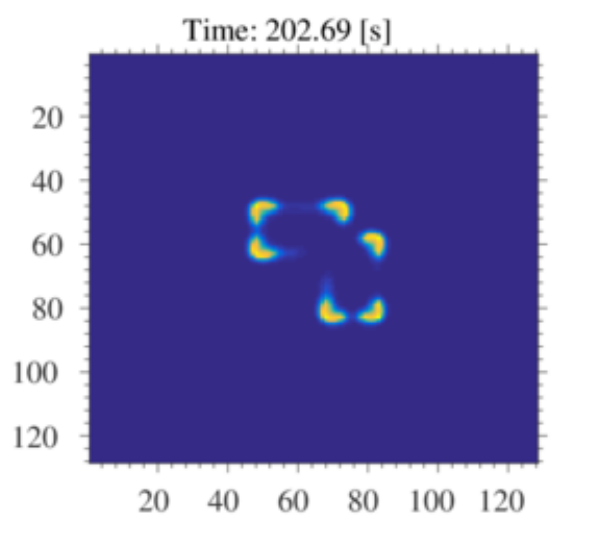}
	\end{subfigure}%
	\hfill
	\begin{subfigure}[t]{0.33\textwidth} 
		\includegraphics[width=0.9\linewidth]{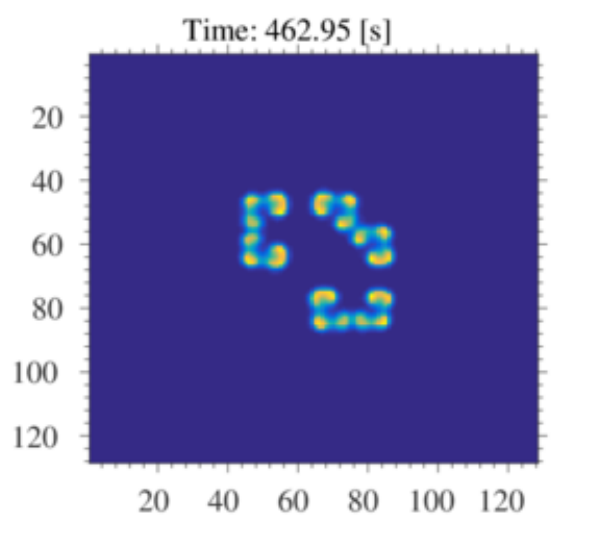}
	\end{subfigure}%

	\begin{subfigure}[t]{0.33\textwidth} 
		\includegraphics[width=0.9\linewidth]{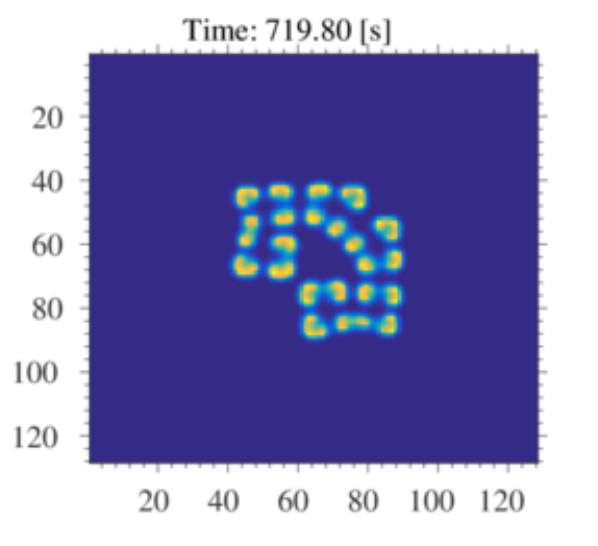}
	\end{subfigure}%
	\hfill  
	\begin{subfigure}[t]{0.33\textwidth} 
		\includegraphics[width=0.9\linewidth]{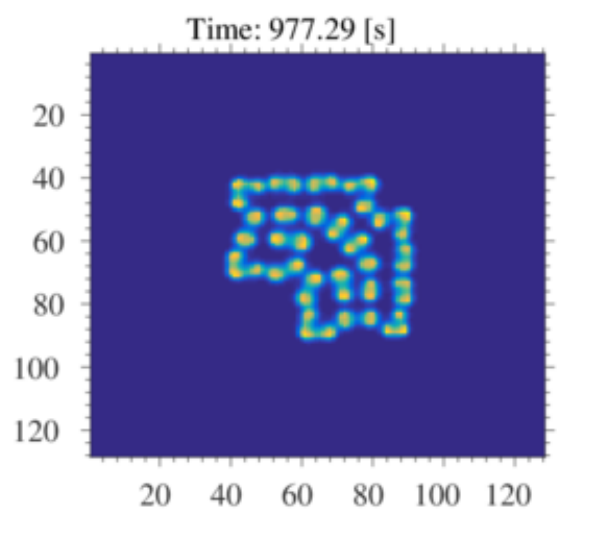}
	\end{subfigure}%
	\hfill
	\begin{subfigure}[t]{0.33\textwidth} 
		\includegraphics[width=0.9\linewidth]{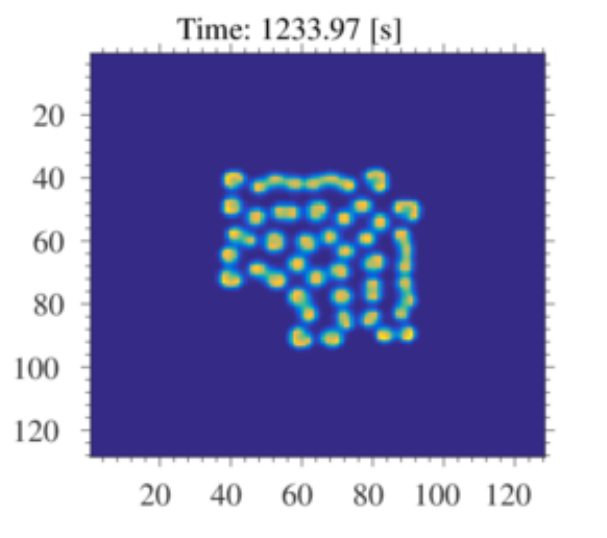}
	\end{subfigure}%
	\hfill
	\caption{Evolution of Gray-Scott model \eqref{eq:Gray-Scott}  solution $u$ in time. A nonlinear diffusion is used, and integration is performed with explicit \mgark methods.}
	\label{fig:Gray-Scott_snapshots}
\end{figure}

\begin{figure}[h]
	\begin{subfigure}[t]{0.5\textwidth} 
		\includegraphics[width=0.95\linewidth]{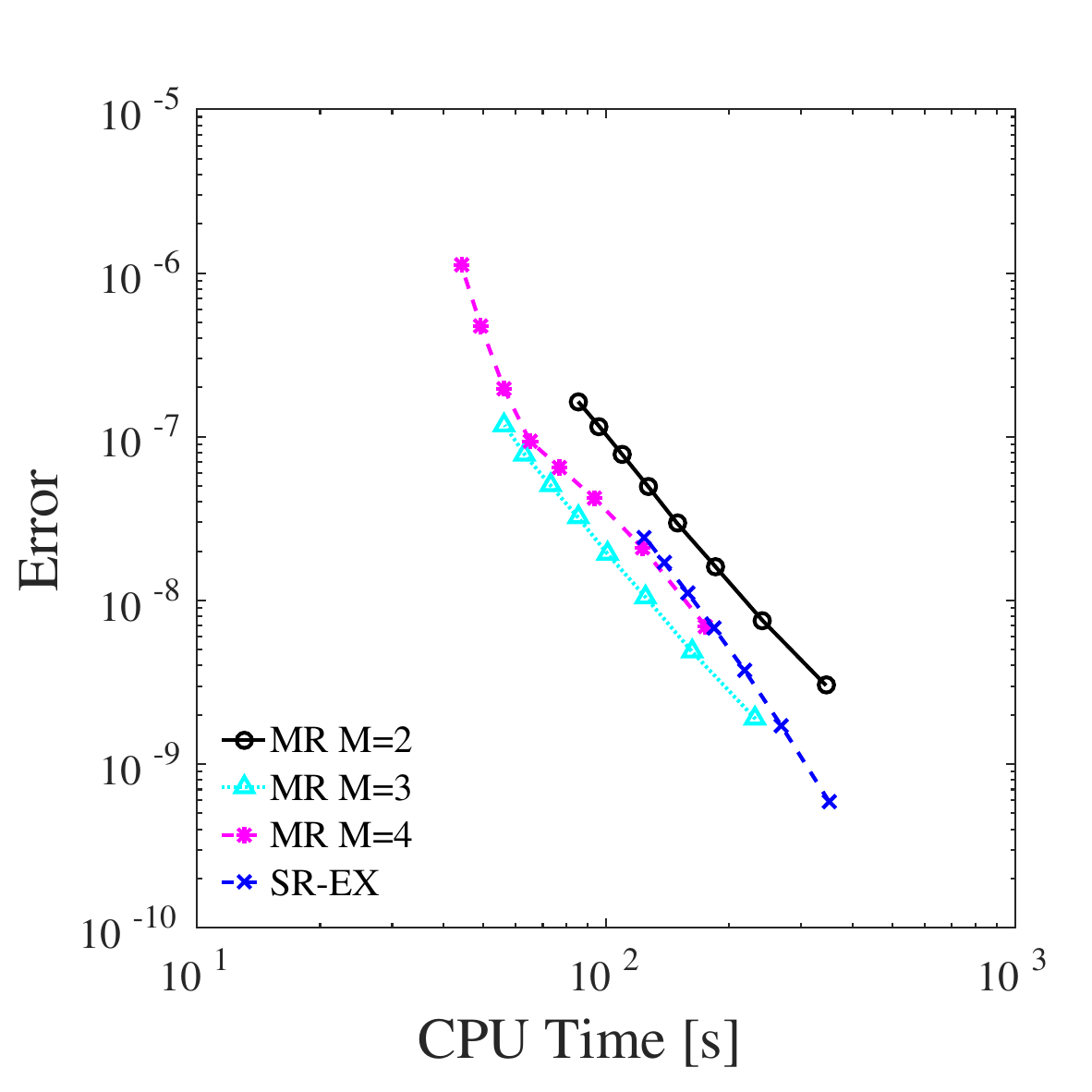}
		\caption{Convergence plot for \method{\explicit}{\explicit}{3}{2}{3}{3}{A}}
		\label{fig:GS_explicit_time}
	\end{subfigure}%
	\hfill
	\begin{subfigure}[t]{0.5\textwidth} 
		\includegraphics[width=0.95\linewidth]{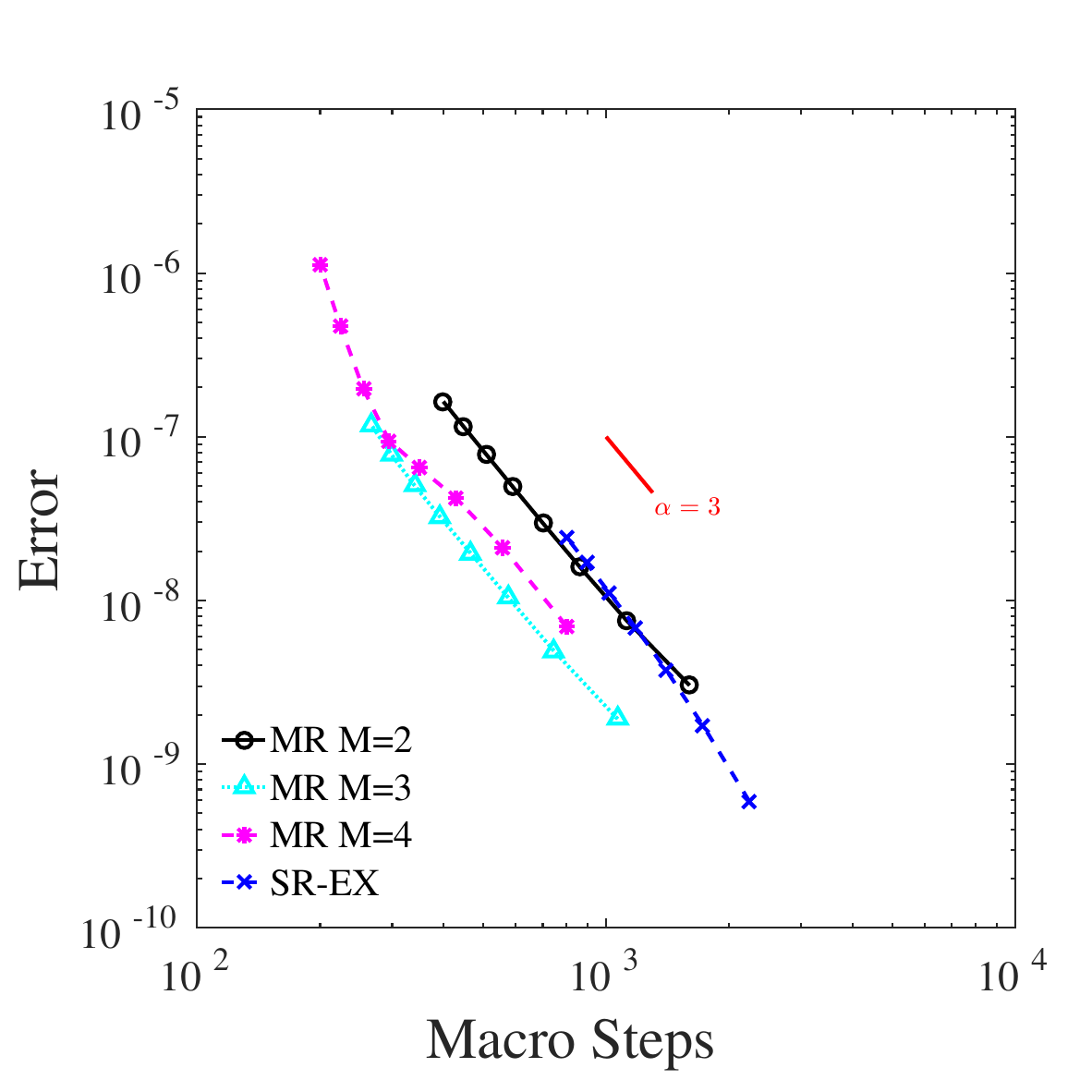}
		\caption{Performance plot for \method{\explicit}{\explicit}{3}{2}{3}{3}{A}}
		\label{fig:GS_explicit_steps}
	\end{subfigure}%

	\begin{subfigure}[t]{0.5\textwidth} 
		\includegraphics[width=0.95\linewidth]{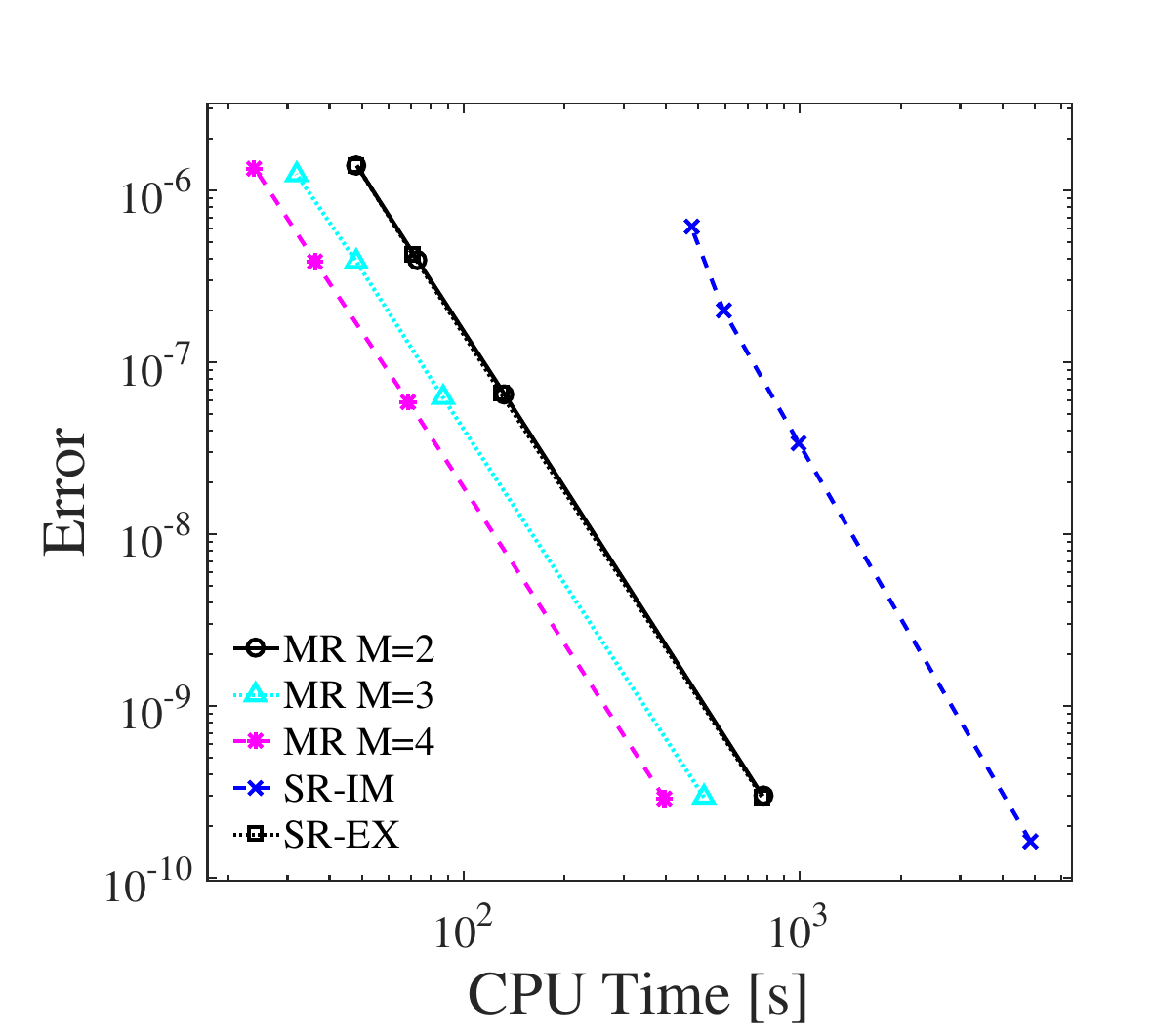}
		\caption{Performance plot for \method{\explicit}{\implicit}{3}{2}{3}{3}{A}}
		\label{fig:GS_exim_time}
	\end{subfigure}%
	\hfill
	\begin{subfigure}[t]{0.5\textwidth} 
		\includegraphics[width=0.95\linewidth]{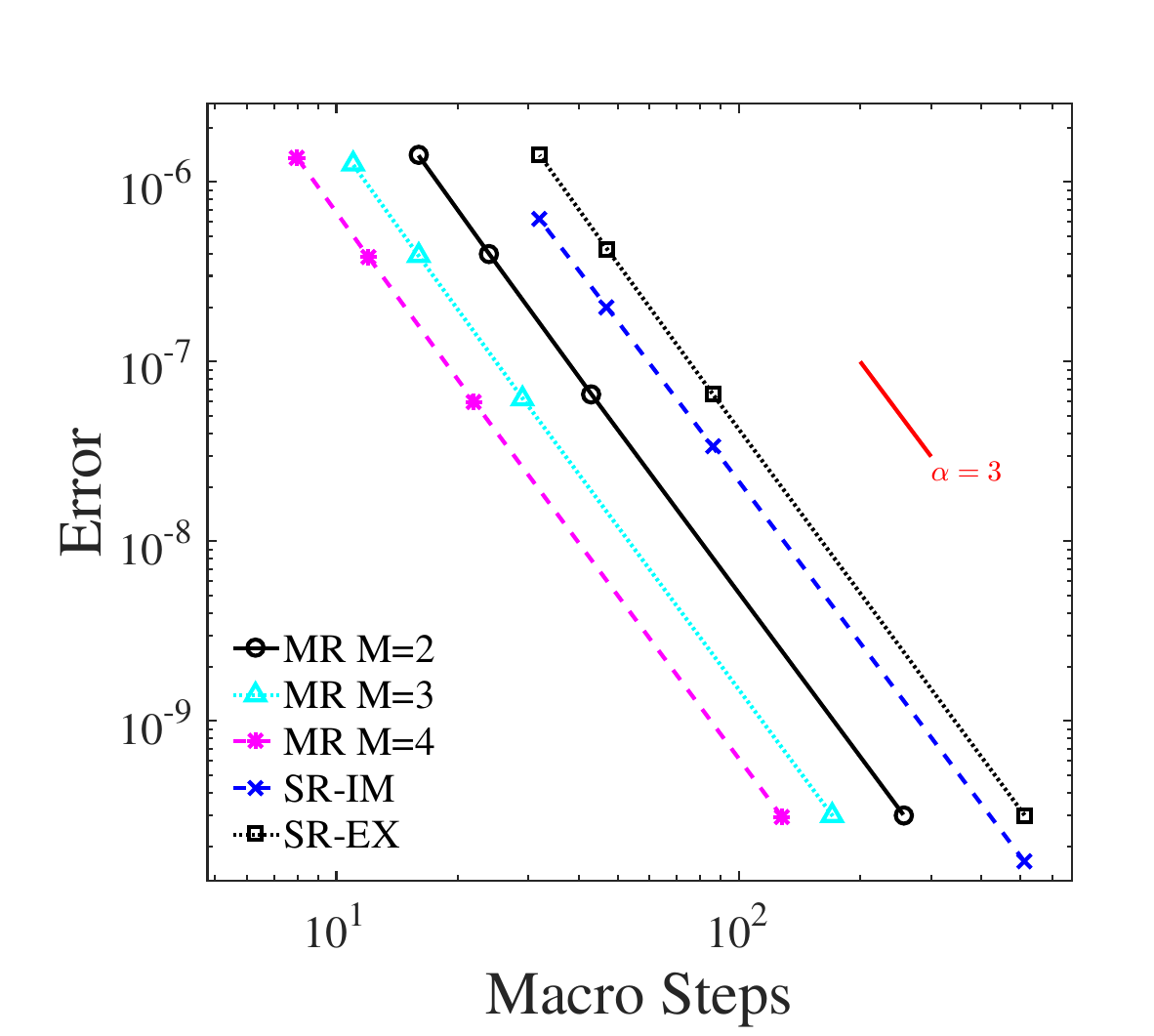}
		\caption{Convergence plot for \method{\explicit}{\implicit}{3}{2}{3}{3}{A}}
		\label{fig:GS_exim_steps}
	\end{subfigure}%
	\hfill
	\caption{Results for \mgark schemes applied to Gray-Scott model  \cref{eq:Gray-Scott}.}
	\label{fig:TimingGS}
\end{figure}

A first version of the model, used to test explicit \mgark methods, includes nonlinear diffusion terms:
\begin{align}
\varepsilon_u = 0.0625\; e^{- \frac{u}{100}} \sin(\pi x) \sin(\pi y), \quad\varepsilon_v = 0.0312 \; e^{- \frac{v}{100}} \sin(\pi x) \sin(\pi y). 
\label{eq:nonlinear_diff_gs_model}
\end{align}
\Cref{fig:Gray-Scott_snapshots} presents the evolution of quantity $u$ over time. \Cref{fig:GS_explicit_time,fig:GS_explicit_steps} present performance diagrams for a fixed-step time integration of Gray-Scott model with nonlinear diffusion. The \mgark method  \method{\explicit}{\explicit}{3}{2}{3}{3}{A} from \cref{subsec:method_EXEX3} with different multirate step ratios is used. Notice that as $M$ increases, it is possible to increase the macro-step size without violating the CFL conditions. The results shown in \cref{fig:GS_explicit_time} indicate that, for a fixed error level, the multirate method shows better performance than single rate one for $M=2$ and $3$.

A second version of the model,  used to test EX-IM \mgark methods, includes linear diffusive terms with parameters $\varepsilon_u = 0.0625$ and $\varepsilon_v = 0.0312$. The diffusion is considered the slow process and is treated implicitly; the constant diffusion operator is leveraged in the solution of linear stage equations.
 \Cref{fig:GS_exim_time,fig:GS_exim_steps} show the performance diagrams for the \mgark method \method{\explicit}{\implicit}{3}{2}{3}{3}{A} from \cref{subsec:method_EXIM3} compared to single rate implicit method of the same order. In \cref{fig:GS_exim_steps} we note the reduced order of convergence for the single rate implicit method. Being a DIRK scheme, this method has a Newton-Krylov iteration for each stage that involves Jacobian-vector products of the full right hand side. Here, as in many practical cases, these are approximated by finite differences. Therefore, the quality of the approximation and the Krylov solver affect the convergence rate and efficiency of the SR-IM method \cite{Sandu_2017_large-CFD}. On the other hand, the stages for multirate methods use only the Jacobian of the linear term, and show full order of convergence and faster computations than both single rate implicit and explicit methods. As $M$ increases the performance of the multirate method improves incurring negligible increase in error.

%

\subsection{Numerical experiments for $H$ and $M$ adaptivity}
We experiment with different adaptivity strategies described in \cref{subsec:adaptivity} using the Gray-Scott model \cref{eq:Gray-Scott} with nonlinear diffusion term \cref{eq:nonlinear_diff_gs_model}. All experiments in this section are performed for a time span of $T=[0,2]$ seconds. 

\begin{figure}[h]
\begin{subfigure}[t]{0.35\textwidth} 
	\includegraphics[width=1\linewidth]{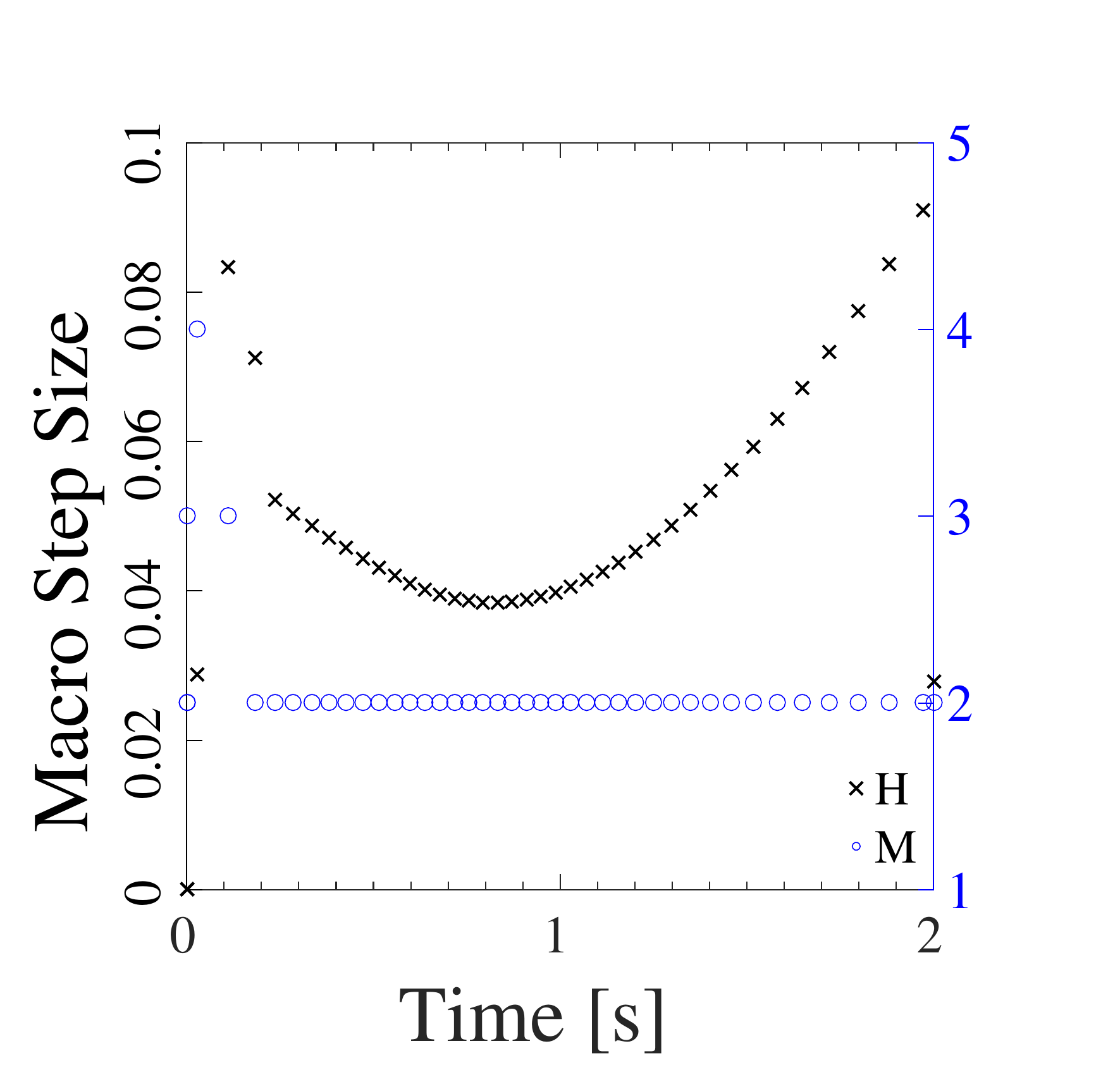}
	\caption{$t^{\{\s\}} / t^{\{\f\}} = 15$.}
	\label{POS_adapt_a}
\end{subfigure}%
\hspace*{-5mm}
\begin{subfigure}[t]{0.35 \textwidth} 
	\includegraphics[width= 1 \textwidth]{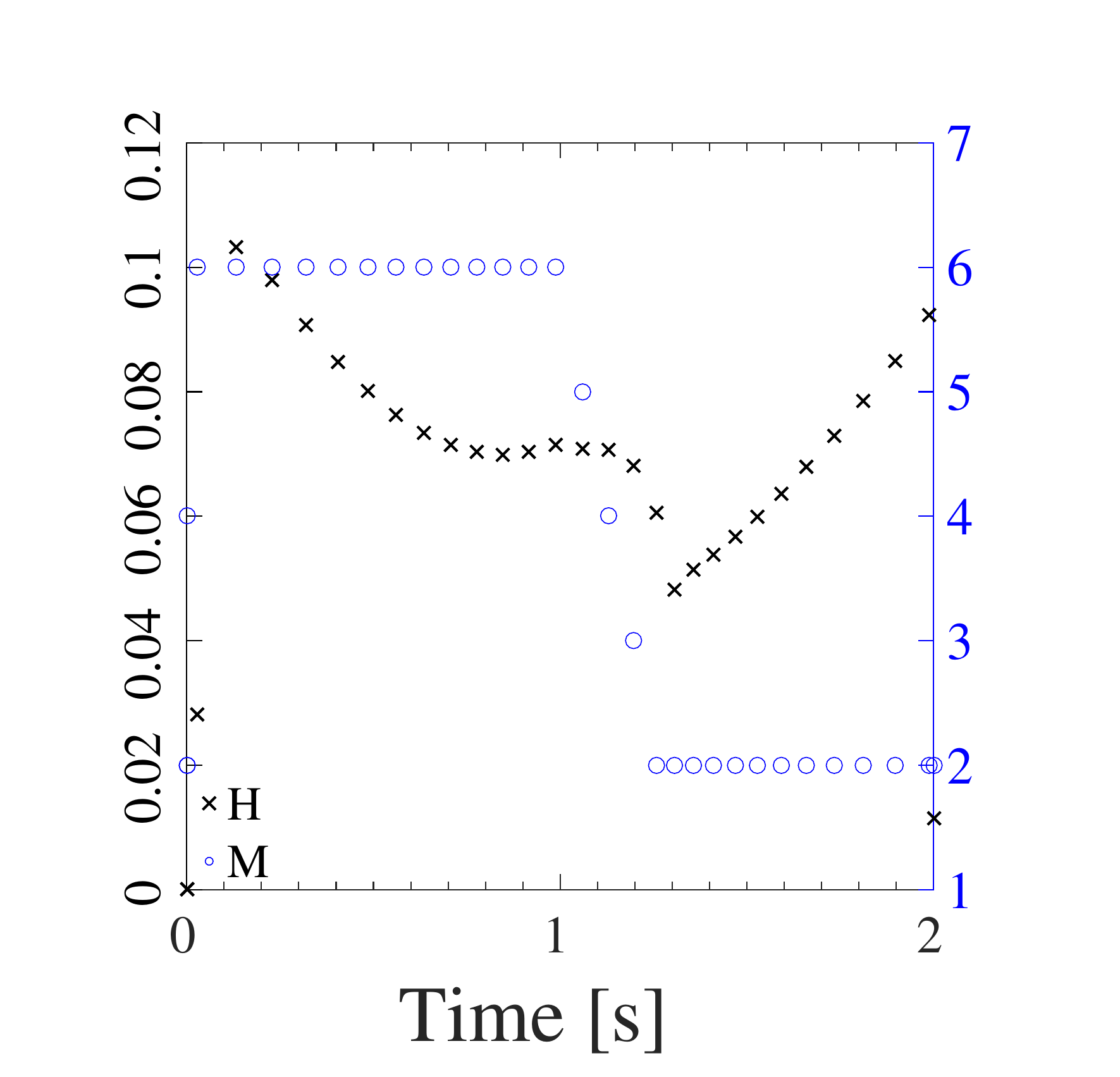}
	\caption{$t^{\{\s\}} / t^{\{\f\}} = 20$.}
	\label{POS_adapt_b}
\end{subfigure}
\hspace*{-5mm}
\begin{subfigure}[t]{0.35 \textwidth} 
	\includegraphics[width= 1\linewidth]{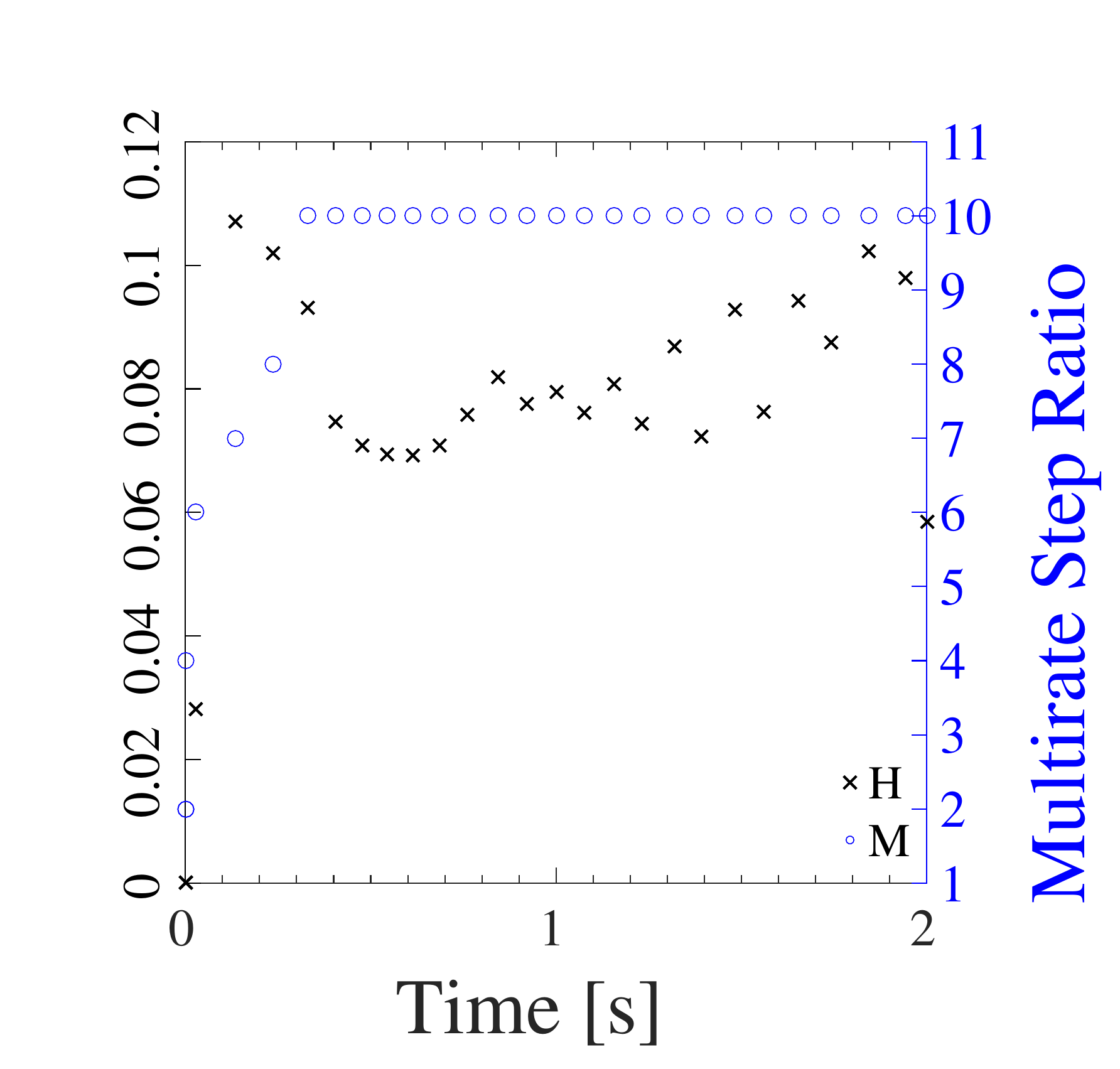}
	\caption{$t^{\{\s\}} / t^{\{\f\}} = 25$.}
	\label{POS_adapt_c}
\end{subfigure}%
\caption{Automatically selected macro-step size and multirate step ratio for the Gray-Scott model \cref{eq:Gray-Scott} integrated over the time span $T = [0, 2]$ using \method{\explicit}{\explicit}{3}{2}{4}{4}[A] method. The efficiency optimization strategy (\cref{subsec:POS_adapt}) is used with $\text{AbsTol}=\text{RelTol}=10^{-4}$.}
\label{fig:POS_adapt}
\end{figure}

\Cref{fig:POS_adapt} shows the automatically selected macro-step size ($H$) and multirate step ratio ($M$) for the naturally adaptive method \method{\explicit}{\explicit}{3}{2}{4}{4}{A}. The efficiency optimization strategy (\cref{subsec:POS_adapt}) is used to select both $H$ and $M$. As the relative cost of evaluating the fast and slow right hand side functions changes in \cref{POS_adapt_a,POS_adapt_b,POS_adapt_c}, the choice of $M$ varies to keep the overall efficiency high. The costs of slow and fast system are computed by timing their respective right hand side evaluations. The one-dimensional integer optimization in \cref{eq:EOS} is then approximated by limiting the possible increase or decrease of $M_{\text{new}}$ to
\begin{align*}
	\max (1,M-1) \leq {M_{\text{new}}} \leq M+2,
\end{align*}
which reduces the computational cost of the optimization and increases the robustness of the algorithm. 

\begin{figure}[h]
	\centering
	\begin{subfigure}[t]{0.35\textwidth} 
		\includegraphics[width=1\linewidth]{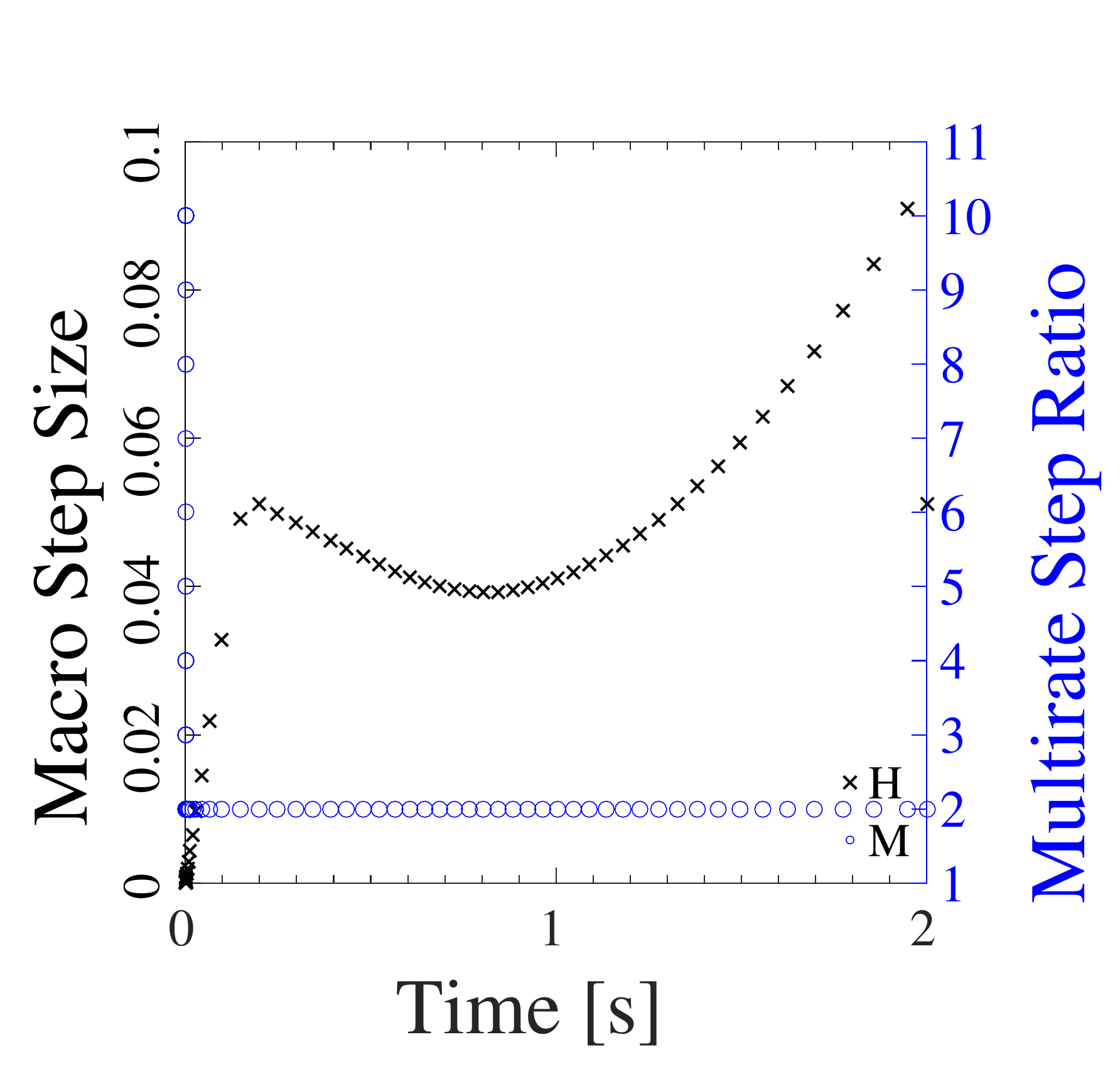}
		\caption{Gray-Scott model with reaction as the fast system}
		\label{fig:BES_adapt_a}
	\end{subfigure}%
	\hspace{5mm}
	\begin{subfigure}[t]{0.35\textwidth} 
		\includegraphics[width=1\linewidth]{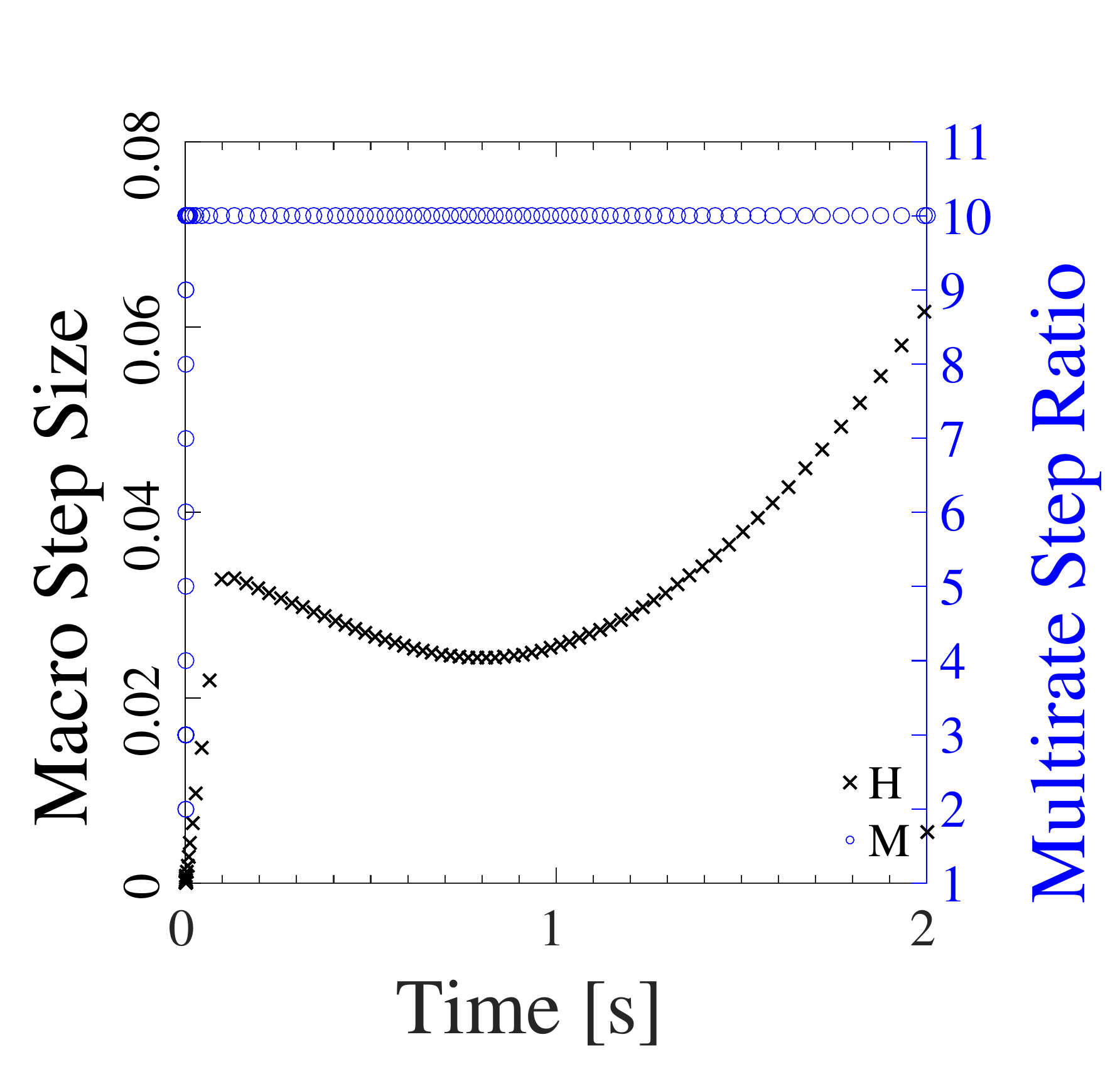}
		\caption{Gray-Scott model with diffusion as the fast system}
		\label{fig:BES_adapt_b}
	\end{subfigure}%
	\hfill
	\caption{Automatically selected macro-step size and multirate step ratio for the Gray-Scott model \cref{eq:Gray-Scott} integrated over the time span $T = [0, 2]$ seconds using \method{\explicit}{\explicit}{3}{2}{4}{4}{A} method. The strategy of balancing the slow and fast errors (\cref{subsec:BES_adapt}) is used with $\text{AbsTol}=\text{RelTol}=10^{-2}$.}
\label{fig:BES_adapt}
\end{figure}

In \cref{fig:BES_adapt} the adaptivity strategy based on balancing the fast and slow errors (\cref{subsec:BES_adapt}) is tested. By interchanging the roles of the fast and slow partitions the choice of multirate step ratio $M$ changes from the maximum allowed value of ten to its minimum allowed value of two, in an attempt to keep the fast and slow errors balanced. The macro-step size selection by the traditional error controller is the same in both cases.

\begin{figure}[h]
	\centering
	\begin{subfigure}[t]{0.45\textwidth} 
		\includegraphics[width=1\linewidth]{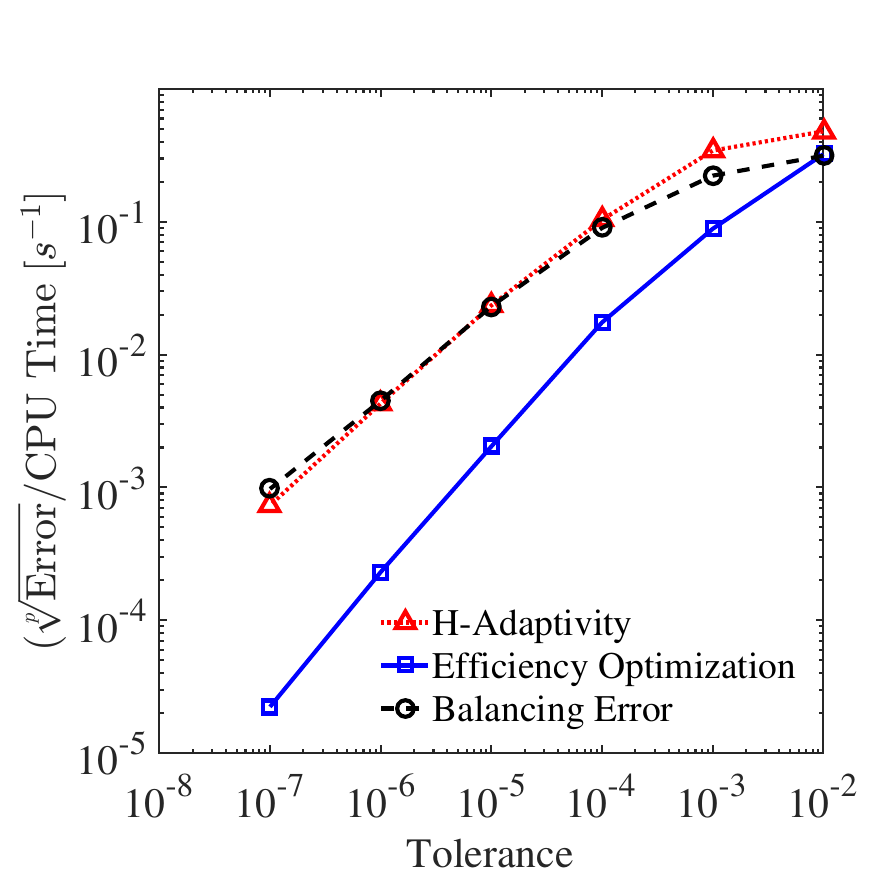}
		\caption{\method{EX}{EX}{2}{1}{2}{2}{A} method.}
		\label{fig:compare_adapt_a}
	\end{subfigure}%
	\hspace{5mm}
	\begin{subfigure}[t]{0.45\textwidth} 
		\includegraphics[width=1\linewidth]{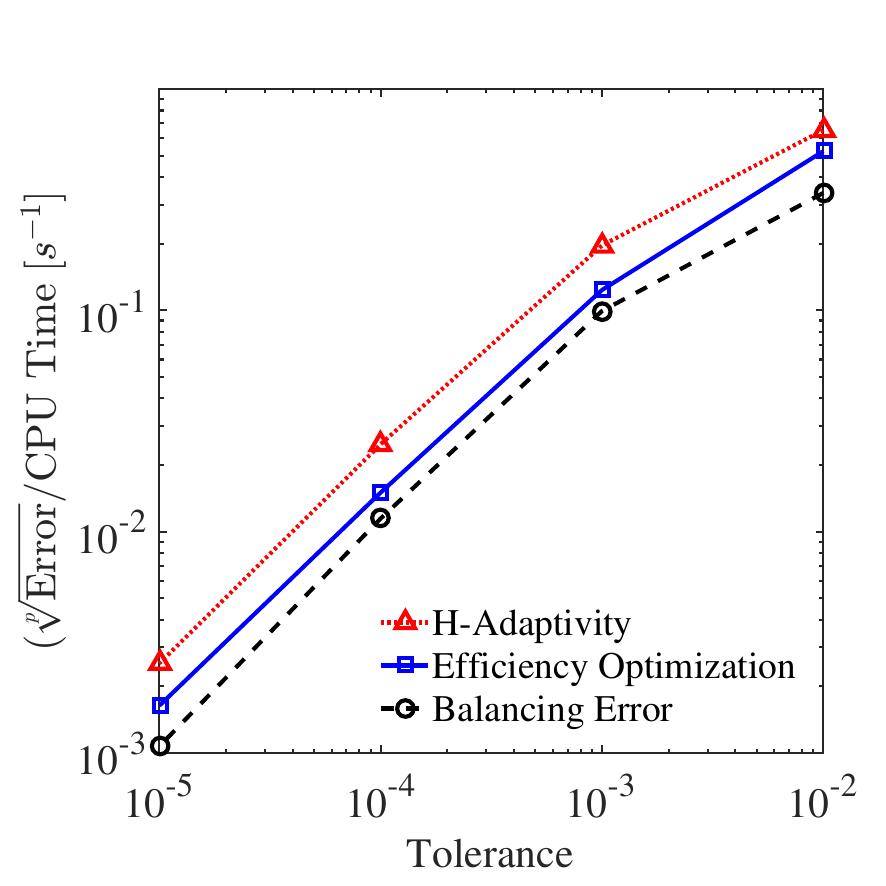}
		\caption{\method{EX}{EX}{4}{3}{5}{5}{A} method. }
		\label{fig:compare_adapt_b}
	\end{subfigure}%
	\hfill
	\caption{Efficiency of different adaptivity strategies for explicit-explicit methods applied to Gray-Scott model with a compute time ratio $t^{\{\s\}} / t^{\{\f\}} = 2.$}
	\label{fig:performance-adapt-strategies}
\end{figure}

Finally, it is instructive to see how different $H$-$M$  adaptivity strategies developed in \cref{subsec:adaptivity} compare against the conventional $H$ adaptivity. We carry out this experiment using the Gray-Scott model. \Cref{fig:compare_adapt_a} shows the error in the final solution scaled by the computation time when the \method{EX}{EX}{2}{1}{2}{2}{A}  method from \cref{subsec:method_EXEX2} is used.   \Cref{fig:compare_adapt_a} repeats the experiment using the \method{EX}{EX}{4}{3}{5}{5}{A} method from  \cref{subsec:method_EXEX4}. In these experiments the efficient single rate base method is employed when the adaptive strategy selects $M=1$. The parameters of the adaptivity algorithm such as macro-step rejection factor are optimized for maximum efficiency. The results indicate that the $H$-$M$ adaptivity strategies are perform better than the classical $H$ error control. 
%
\subsection{Component partitioning experiment} 
\begin{figure}[h]
	\includegraphics[width=\linewidth]{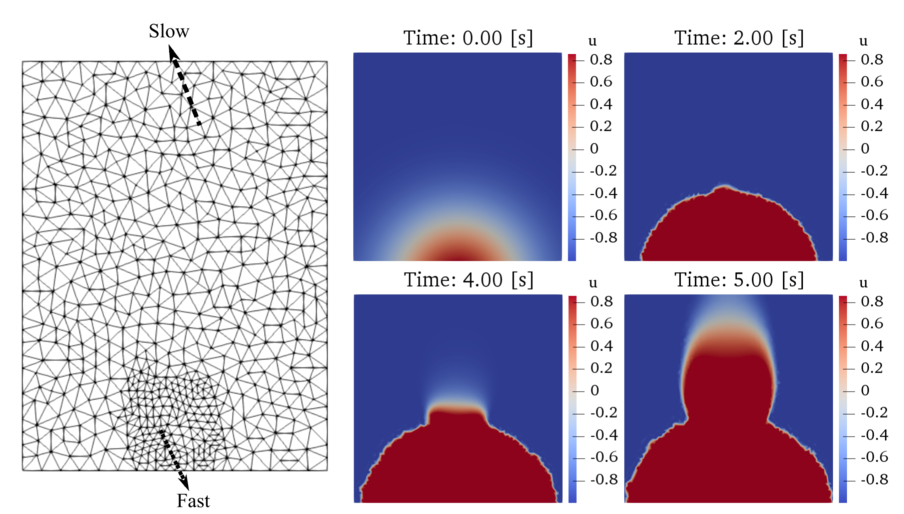}
	\caption{Evolution of the BSVD problem \eqref{eq:BSVD} solution in time.}
	\label{fig:BSVD_snapshots}
\end{figure}

\begin{figure}[h]
	\begin{subfigure}[t]{0.5\textwidth} 
		\includegraphics[width=0.95 \linewidth]{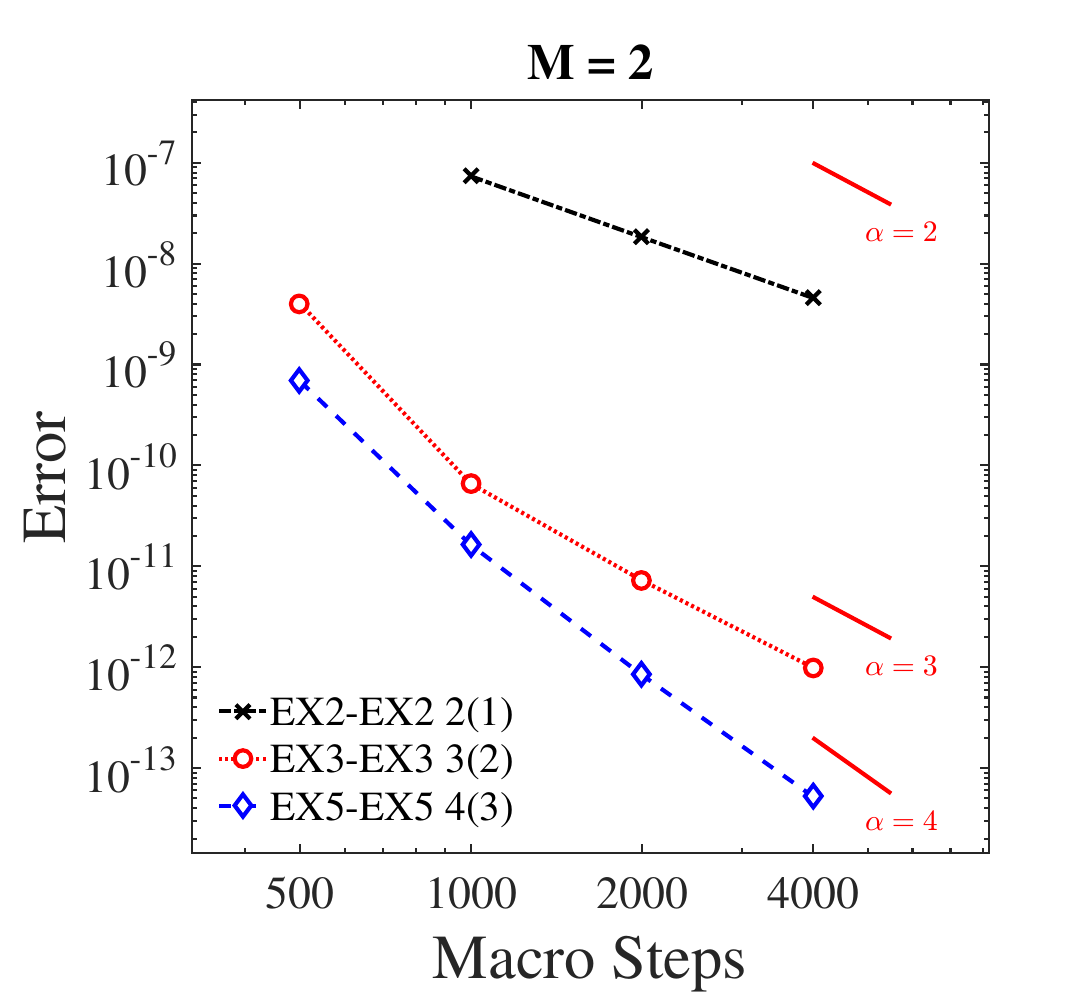}
	\end{subfigure}%
	\hfill
	\begin{subfigure}[t]{0.5\textwidth} 
		\includegraphics[width=0.99 \linewidth]{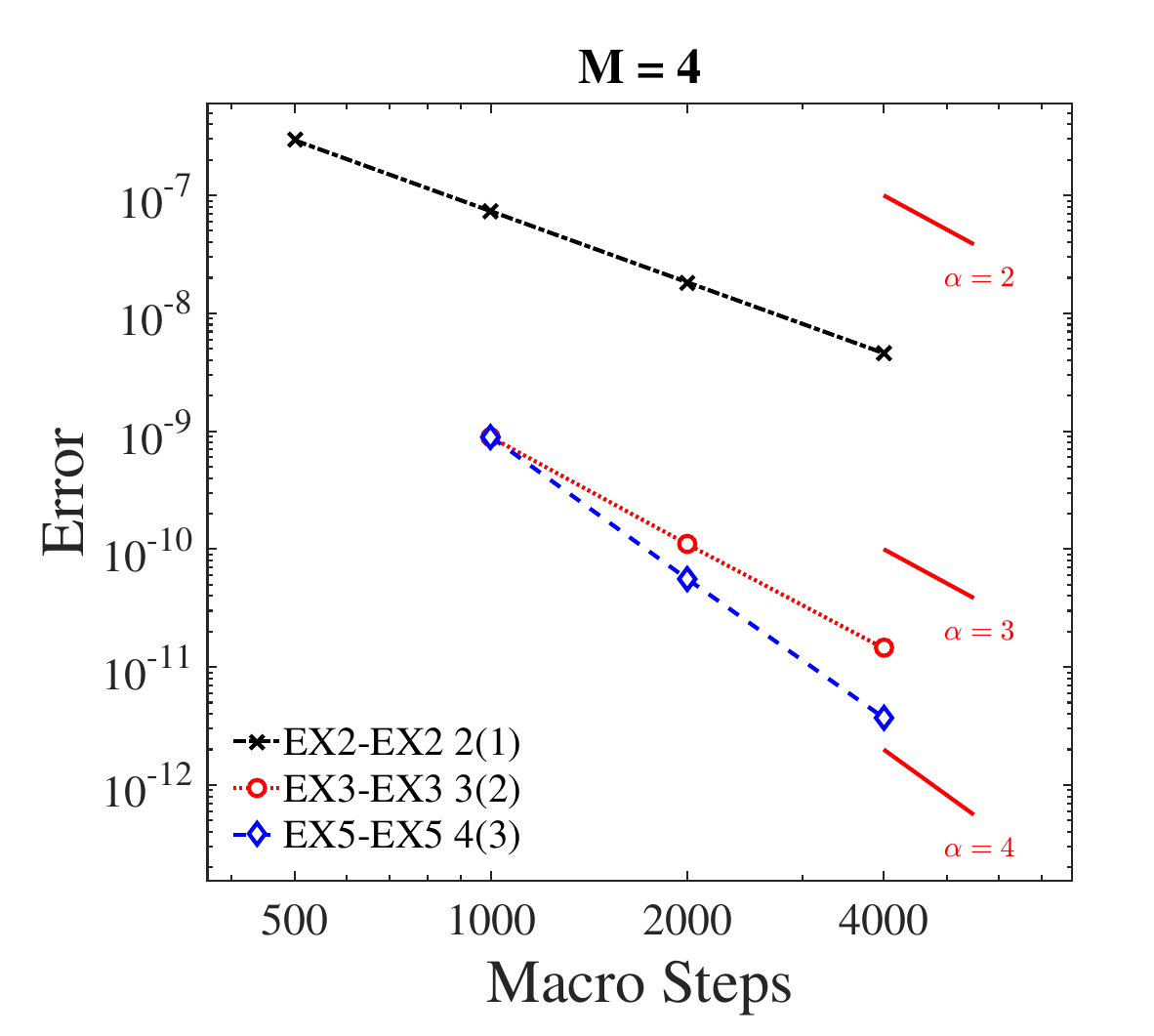}
	\end{subfigure}%
	\hfill
	\begin{subfigure}[t]{0.5\textwidth} 
	\end{subfigure}%
	\caption{Convergence results for BSVD test \cref{eq:BSVD} over the time span $T = [0, 0.2]$ seconds. A fixed macro-step time integration is carried out with varying multirate step ratios $M$ using explicit \mgark type A methods.}
	\label{fig:BSVD_convergence}
\end{figure}

We consider the BSVD reaction-diffusion problem on the unit square $\Omega$ \cite{heineken2006partitioning} with boundary $\partial \Omega$ and outward boundary normal vector $\vec{n}$: 
\begin{subequations}
	\begin{alignat}{2}
	&u_t = \nabla \cdot \left( D(x,y) \nabla u \right)
	+10(1-u^2)(u + 0.6),&& \\
	&u(x,y, 0) = 2 \exp{ \left( -10(x-0.5)^2 -10 (y+0.1)^2 \right)}, \qquad &x,y &\in \Omega, \\
	&D(x,y) = 0.1 \sum_{i=1}^{3} e^{-100(x-0.5)^2 + (y -y_i)^2},&&\\
	&D(x,y) \nabla u \cdot \vec{n} = 0, \qquad & x,y&\in \partial \Omega ,\; t \in [0,t_F].
	\end{alignat}
	\label{eq:BSVD}
\end{subequations}
The parameters values are $y_1 = 0.6, y_2 = 0.75, y_3 = 0.9$ as prescribed in \cite{heineken2006partitioning}. The domain is partitioned into two sub-domains:
 \begin{subequations}
 	\begin{align*}
		\Omega_1 = \{  (x,y)\in \Omega \; : \; \abs{x-0.5}\leq 0.125, ~ y \leq 0.125 \},\qquad \Omega_2 = \Omega \backslash \Omega_1.
 	\end{align*}
 \end{subequations}
A continuous finite element semi-discretization in space with Lagrange polynomial basis of order four leads to a two-way partitioned system of ODEs. We discretized this model using the FEniCS package \cite{AlnaesBlechta2015a} and the evolution of its solution over the time span of $T=[0,5]$ seconds is shown in \cref{fig:BSVD_snapshots}. The subsystem with fewer degrees of freedom (DOFs) is designated as the fast one, and the remaining variables are considered to be the slow system.   In this experiment the ratio of slow to fast DOFs is 28. \cref{fig:BSVD_convergence} shows the error for fixed macro-step time integration using explicit-explicit multirate methods of orders 2, 3  and 4. The largest macro-steps are chosen close to the maximum stable step size for the method. Results reported in \cref{fig:BSVD_convergence} show that, although the numerical order of convergence of the methods fluctuate slightly, they follow their theoretical values closely. Note that the errors are reported only for stable macro-step sizes for different methods and multirate step ratios. 
%
%
%
%
\section{Conclusions}
\label{sec:Conclusions}

This work develops a coherent design strategy for high order \mgark, and constructs several particular schemes of explicit-explicit, explicit-implicit and implicit-explicit types up to order four. Adopting a two-way split system for simplicity, we identify coupling structures for the fast and slow components that strike a balance between stability and computational efficiency of \mgark methods. Furthermore, in the process of designing schemes we explore a variety of optimizations such as local truncation error minimization, FSAL, and stiff accuracy of base and overall implicit methods. When possible, our methods are endowed with telescopic properties to facilitate their application to multi-scale, multi-domain problems. Different design criteria lead to different types of methods, e.g., the A-type is optimized for accuracy, and the S-type is optimized for simplicity and stability. A novel concept of $H$-$M$ adaptivity, based on multiple error estimators and the property of natural adaptivity, is presented.

Our numerical experiments demonstrate several applications of these methods with finite element and finite difference discretized PDEs, where different time-scales in system variables (component partitioning), or in various physical processes of the problem (additive partitioning) are treated with different time steps. In all cases, the proposed methods converge at their theoretical orders. In all cases, speedups are obtained from the multirate approach when compared to the single rate integration. Since the multirate performance depends critically on the efficiency of the implementation, and a high quality implementation is not within the scope of this work, we do not report these speedups, as they are likely considerably smaller than what can realistically be achieved. The future work plans of the authors include the development of high quality implementations of \mgark schemes. Moreover, for many of the new schemes derived here the stability region decreases with increasing step size ratio; future work will focus on the development of S-type schemes where the stability region remains large for any $M$.

This contribution focuses on {\it decoupled} multirate methods satisfying \cref{eq:decoupled-methods-equation}, as they are far less computationally demanding than coupled methods. However, the overall stability of the scheme may be affected by this choice. In future work the authors plan to explore coupled implicit-implicit \mgark methods, and to study their stability via an extended, two-dimensional test problem \cite{Sandu_2013_extrapolatedMR}.

\bibliographystyle{siamplain}
\bibliography{Bib/misc,Bib/sandu,Bib/ode_multirate,Bib/ode_general}
\clearpage
\appendix
\section{Decoupled multirate GARK schemes}
\label{sec:new_mr_schemes}

We use the naming convention {\it \method{FAST}{SLOW}{$p$}{$\widehat{p}$}{$f$}{$s$}{type}}, where $p$ is the method order, $\widehat{p}$ is the embedded order, $f$ is the number of stages in the fast base method, and $s$ is the number of stages in the slow base method. Each component method is either explicit or implicit: FAST, SLOW $\in$ \{EX,IM\}. We distinguish between methods of type A (optimized for accuracy and for better step size control) and methods of type S (optimized for simplicity and for stability), therefore  $type \in$ \{A,S\}.

\subsection{\mgark \method{\explicit}{\explicit}{2}{1}{2}{2}{A}}
\label{subsec:method_EXEX2}
  This explicit method uses a base method from \cite{ralston1962runge} and has telescopic \cref{eqn:mr-telescopic-conditions} and naturally adaptive (\cref{def:naturally-adaptive}) properties.
  
\fitbox{3}{
  	$A^{\{\f,\f\}}=
	  	\begin{bmatrix}
	  	0 & 0 \\
	  	\frac{2}{3} & 0 \\
	  	\end{bmatrix},$
  	&
  	$A^{\{\s,\s\}} =
	  	\begin{bmatrix}
	  		0 & 0 \\
	  		\frac{2}{3} & 0 \\
	  	\end{bmatrix},$
  	&
  	$A^{\{\f,\s,1\}}  =
	  	\begin{bmatrix}
	  	0 & 0 \\
	  	\frac{2}{3 M} & 0 \\
	  	\end{bmatrix},$ \\
  	\fullrow{$A^{\{\f,\s,\lambda\}} =
  		\begin{bmatrix}
  		\frac{3 M^3-11 M^2+20 \lambda  M-20 M-20 \lambda +20}{20 (M-1) M} & -\frac{M (3 M-11)}{20 (M-1)} \\
  		\frac{-3 M^3-9 M^2+60 \lambda  M-20 M-60 \lambda +20}{60 (M-1) M} & \frac{M (M+3)}{20 (M-1)} \\
  		\end{bmatrix}, \qquad \lambda = 2, \ldots, M, $}  \\
  	$A^{\{\s,\f,1\}}  =
	  	\begin{bmatrix}
	  	0 & 0 \\
	  	-\frac{1}{3} (M-2) M & \frac{M^2}{3} \\
	  	\end{bmatrix},$
  	&
  	$A^{\{\s,\f,\lambda\}} =
	  	\begin{bmatrix}
	  	0 & 0 \\
	  	0 & 0 \\
	  	\end{bmatrix},$
  	&
  	$\lambda = 2, \ldots, M,$  \\
  	$b^{{\{\f\}}} =   b^{{\{\s\}}} =
	  	\begin{bmatrix}
	  	\frac{1}{4} & \frac{3}{4} \\
	  	\end{bmatrix}\trsym,$
  	&
    $\widehat{b}^{{\{\f\}}}  =   \widehat{b}^{{\{\s\}}} =
	    \begin{bmatrix}
	    1 & 0 \\
	    \end{bmatrix}\trsym.$ 
	&
}

\begin{figure}[ht!]
	\centering
	\begin{subfigure}[b]{0.45\textwidth}
		\centering
		\includegraphics[width=\textwidth]{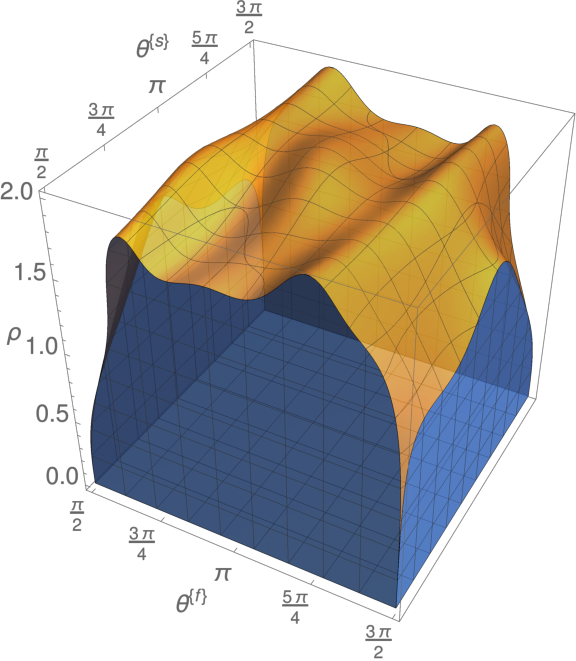}
		\caption{$M=2$.}
	\end{subfigure}
	\hfill
	\begin{subfigure}[b]{0.45\textwidth}
		\centering
		\includegraphics[width=\textwidth]{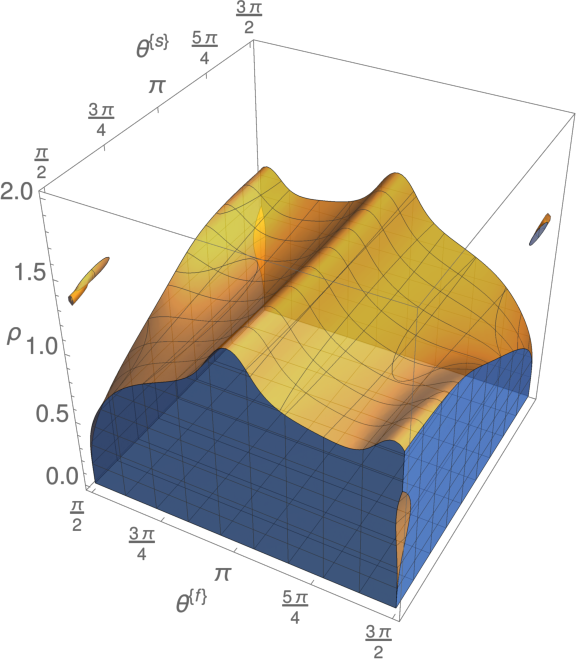}
		\caption{$M=4$.}
	\end{subfigure}
	\caption{\method{\explicit}{\explicit}{2}{1}{2}{2}{A} stability regions.}
\end{figure}

\newpage
\subsection{\mgark \method{\explicit}{\explicit}{2}{1}{2}{2}{S}}
\label{subsec:method_EXEX2s}
  This explicit method uses a two stage base method and has telescopic \cref{eqn:mr-telescopic-conditions} and naturally adaptive (\cref{def:naturally-adaptive}) properties.  The scheme computes the first slow stage, followed by $L_2$ fast steps, then the second slow stage, followed by $M-L_2$ fast steps.  For example, we can take $L_2 = floor(c_2 M)$.

\fitbox{2}{
	$A^{\{\f,\f\}}=
	\begin{bmatrix}
	0 & 0 \\
	c_2 & 0 \\
	\end{bmatrix},$
	&
	$A^{\{\s,\s\}} =
	\begin{bmatrix}
	0 & 0 \\
	c_2 & 0 \\
	\end{bmatrix},$
	\\
	$A^{\{\f,\s,\lambda\}} =
		\begin{bmatrix}
		\frac{\lambda -1}{M} & 0 \\
		\frac{\lambda +c_2-1}{M} & 0 \\
		\end{bmatrix},$
	& $\lambda = 1, \ldots, L_2, $  \\
	$A^{\{\f,\s,\lambda\}} =
	\begin{bmatrix}
	\frac{(\lambda -1) \left(2 c_2-1\right)}{2 M c_2} & \frac{\lambda -1}{2 M c_2} \\
	\frac{M}{3 \left(L_2-M\right)}+\frac{-\lambda +2 c_2 \left(2 \lambda +c_2-2\right)+1}{2 M c_2} & \frac{M}{3 M-3 L_2}+\frac{\lambda -1}{2 M c_2}+\frac{1-\lambda }{M} \\
	\end{bmatrix},$
	& $\lambda = L_2 + 1, \ldots, M, $  \\
	$A^{\{\s,\f,\lambda\}} =
	\begin{bmatrix}
	0 & 0 \\
	\frac{M \left(-2 M+6 c_2+3 L_2-3\right)}{6 L_2} & \frac{M \left(2 M-3 L_2+3\right)}{6 L_2} \\
	\end{bmatrix},$
	&
	$\lambda = 1, \ldots, L_2,$  \\
	$A^{\{\s,\f,\lambda\}} =
	\begin{bmatrix}
	0 & 0 \\
	0 & 0 \\
	\end{bmatrix},$
	&
	$\lambda = L_2 + 1, \ldots, M,$  \\
	$b^{{\{\f\}}} =   b^{{\{\s\}}} =
	\begin{bmatrix}
	\frac{2 c_2-1}{2 c_2} & \frac{1}{2 c_2} \\
	\end{bmatrix}\trsym,$
	&
	$\widehat{b}^{{\{\f\}}}  =   \widehat{b}^{{\{\s\}}} =
	\begin{bmatrix}
	1 & 0 \\
	\end{bmatrix}\trsym.$
}

\begin{figure}[ht!]
	\centering
	\begin{subfigure}[b]{0.45\textwidth}
		\centering
		\includegraphics[width=\textwidth]{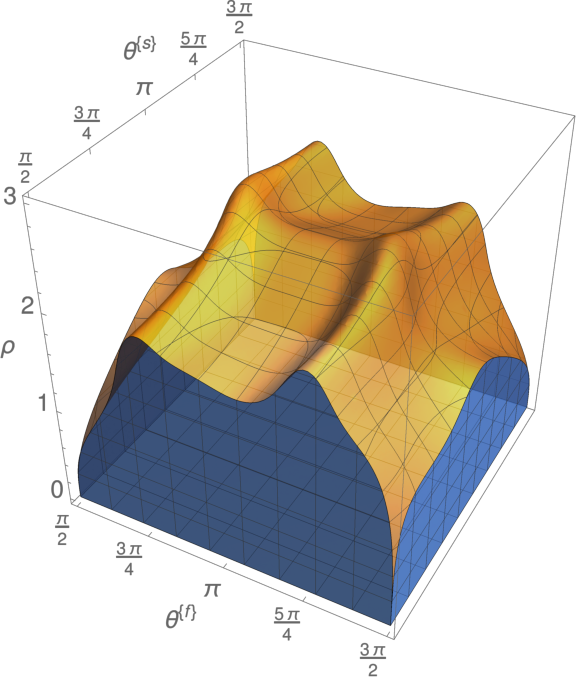}
		\caption{$M=2$.}
	\end{subfigure}
	\hfill
	\begin{subfigure}[b]{0.45\textwidth}
		\centering
		\includegraphics[width=\textwidth]{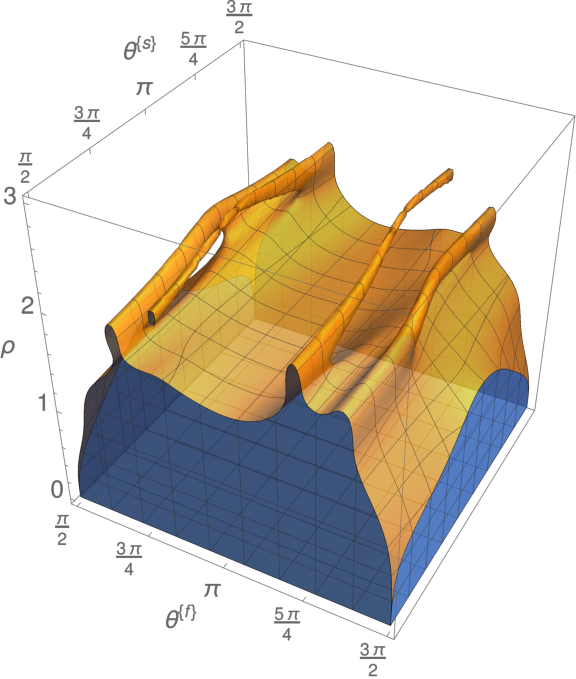}
		\caption{$M=4$.}
	\end{subfigure}
	\caption{\method{\explicit}{\explicit}{2}{1}{2}{2}{S} stability regions with $c_2 = \frac{2}{3}$.}
\end{figure}

\newpage
\subsection{\mgark \method{\explicit}{\implicit}{2}{1}{2}{2}{A}}
\label{subsec:method_EXIM2}
This explicit-implicit method uses the fast method from \cite{ralston1962runge} and slow method from \cite{Alexander_1977_SDIRK}.
The multirate scheme is stiffly accurate \eqref{eqn:stiff-accuracy} in the slow partition and the coupling error is independent of $M$.

\fitbox{3}{
		$A^{\{\f,\f\}}= 
			\begin{bmatrix}
			0 & 0 \\
			\frac{2}{3} & 0 \\
			\end{bmatrix},$
		&
		$A^{\{\s,\s\}} = 
			\begin{bmatrix}
				1-\frac{1}{\sqrt{2}} & 0 \\
				\frac{1}{\sqrt{2}} & 1-\frac{1}{\sqrt{2}} \\
			\end{bmatrix},$
		&
		$A^{\{\s,\f,1\}}=
		\begin{bmatrix}
		M-\frac{M}{\sqrt{2}} & 0 \\
		\frac{1}{4} & \frac{3}{4} \\
		\end{bmatrix},$ 
		\\
		$	A^{\{\f,\s,\lambda\}} =
			\begin{bmatrix}
			\frac{\lambda -1}{M} & 0 \\
			\frac{3 \lambda -1}{3 M} & 0 \\
			\end{bmatrix},$
		&
		$\lambda = 1, \ldots, M,$
  	    &
		\\
		$A^{\{\s,\f,\lambda\}} =
			\begin{bmatrix}
				0 & 0 \\
				\frac{1}{4} & \frac{3}{4} \\
			\end{bmatrix},$
		&
		$\lambda = 2, \ldots, M,$
		&          	
		\\
	   ${b}^{{\{\f\}}}= 
		\begin{bmatrix}
		\frac{1}{4} & \frac{3}{4} \\
		\end{bmatrix}\trsym,$
		&
		$b^{{\{\s\}}} = 
			\begin{bmatrix}
			\frac{1}{\sqrt{2}} & 1-\frac{1}{\sqrt{2}} \\
			\end{bmatrix}\trsym,$
		&
		\\
		$\widehat{b}^{{\{\f\}}}=
			\begin{bmatrix}
				1 & 0 \\
			\end{bmatrix}\trsym,$
			&
			$\widehat{b}^{{\{\s\}}} =
			\begin{bmatrix}
				\frac{3}{5} & \frac{2}{5} \\
			\end{bmatrix}\trsym.$   
		&
}

\begin{figure}[ht!]
\centering
\begin{subfigure}[b]{0.45\textwidth}
	\centering
	\includegraphics[width=\textwidth]{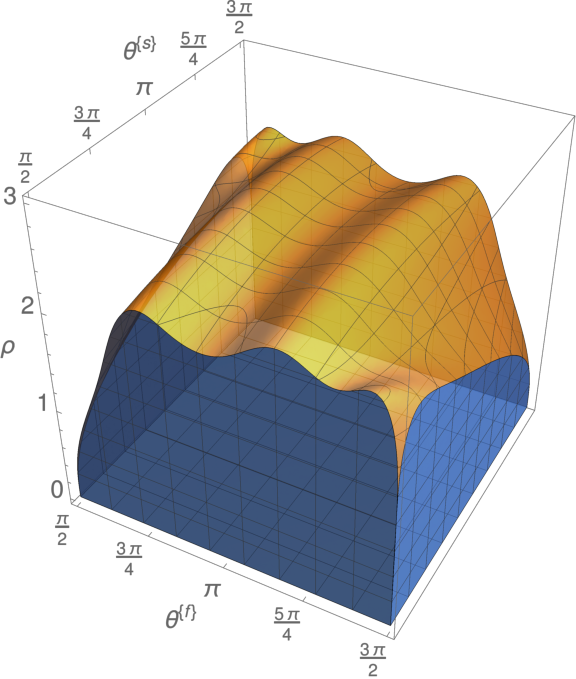}
	\caption{$M=2$.}
\end{subfigure}
\hfill
\begin{subfigure}[b]{0.45\textwidth}
	\centering
	\includegraphics[width=\textwidth]{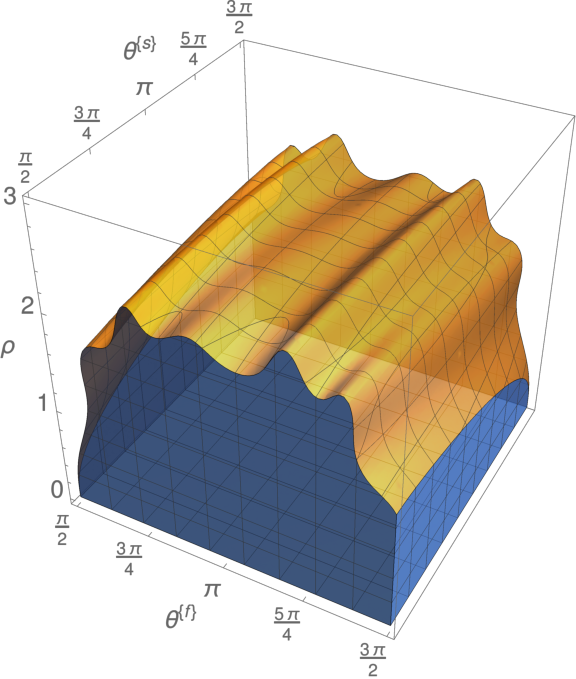}
	\caption{$M=4$.}
\end{subfigure}
\caption{\method{\explicit}{\implicit}{2}{1}{2}{2}{A} stability regions.}
\end{figure}

\newpage
\subsection{\mgark \method{\implicit}{\explicit}{2}{1}{2}{2}{A}}
\label{subsec:method_IMEX2}
This implicit-explicit method uses the fast method from \cite{Alexander_1977_SDIRK} and the slow method from \cite{ralston1962runge}.
The multirate scheme is stiffly accurate \eqref{eqn:stiff-accuracy} in the fast partition.

\fitbox{3}{
		$A^{\{\f,\f\}}= 
		\begin{bmatrix}
		1-\frac{1}{\sqrt{2}} & 0 \\
		\frac{1}{\sqrt{2}} & 1-\frac{1}{\sqrt{2}} \\
		\end{bmatrix},$
		&
		$A^{\{\s,\s\}}=
		\begin{bmatrix}
		0 & 0 \\
		\frac{2}{3} & 0 \\
		\end{bmatrix},$
		&
		$A^{\{\f,\s,M\}}=
		\begin{bmatrix}
		\frac{2 M-\sqrt{2}}{2 M} & 0 \\
		\frac{1}{4} & \frac{3}{4} \\
		\end{bmatrix},$ 
		\\
		$	A^{\{\f,\s,\lambda\}}=
		\begin{bmatrix}
		\frac{2 \lambda -\sqrt{2}}{2 M} & 0 \\
		\frac{\lambda }{M} & 0 \\
		\end{bmatrix},$
		&
		$\lambda = 1, \ldots, M-1,$
		&
		\\
		$A^{\{\s,\f,\lambda\}} = 
		\begin{bmatrix}
		0 & 0 \\
		\frac{2}{3} & 0 \\
		\end{bmatrix},$
		&
		$\lambda = 1, \ldots, M,$
		&          
		\\
		$b^{{\{\f\}}} =
		\begin{bmatrix}
		\frac{1}{\sqrt{2}} & 1-\frac{1}{\sqrt{2}} \\
		\end{bmatrix}\trsym,$
		&
		${b}^{{\{\s\}}}=  
		\begin{bmatrix}
		\frac{1}{4} & \frac{3}{4} \\
		\end{bmatrix}\trsym,$
		&
		\\
		$\widehat{b}^{{\{\f\}}} =
		\begin{bmatrix}
		\frac{3}{5} & \frac{2}{5} \\
		\end{bmatrix}\trsym.$   
		& 
		$\widehat{b}^{{\{\s\}}}= 
		\begin{bmatrix}
		1 & 0 \\
		\end{bmatrix}\trsym,$
		&
}

\begin{figure}[ht!]
\centering
\begin{subfigure}[b]{0.45\textwidth}
	\centering
	\includegraphics[width=\textwidth]{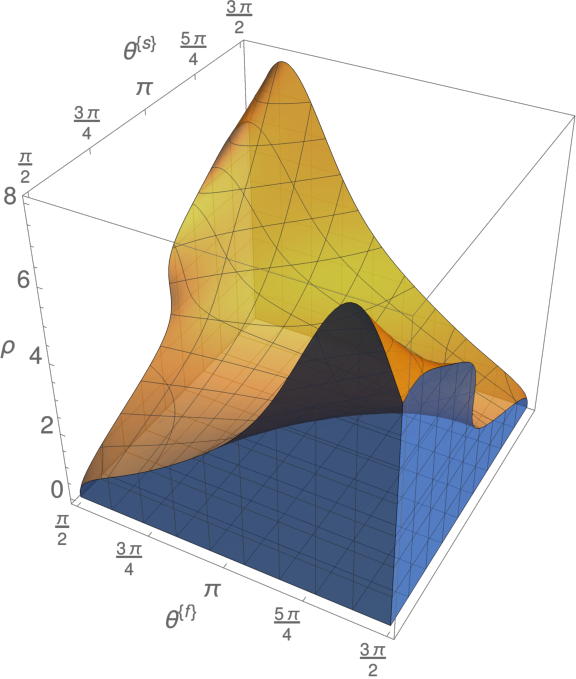}
	\caption{$M=2$.}
\end{subfigure}
\hfill
\begin{subfigure}[b]{0.45\textwidth}
	\centering
	\includegraphics[width=\textwidth]{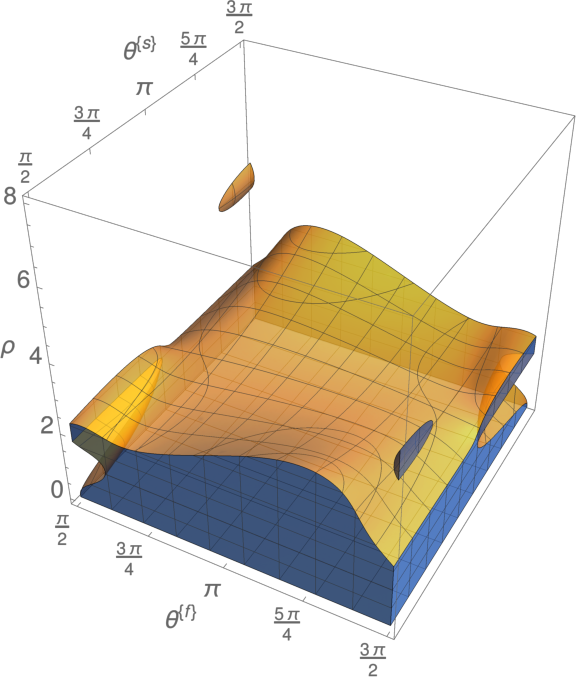}
	\caption{$M=4$.}
\end{subfigure}
\caption{\method{\implicit}{\explicit}{2}{1}{2}{2}{A} stability regions.}
\end{figure}

\newpage
\subsection{\mgark \method{\explicit}{\explicit}{3}{2}{3}{3}{A}}
\label{subsec:method_EXEX3}
This explicit method uses a base method from \cite{ralston1962runge} and has telescopic \cref{eqn:mr-telescopic-conditions} properties. The coupling error is independent of $M$.

\fitbox{3}{
		$A^{\{\f,\f\}}=
		\begin{bmatrix}
		0 & 0 & 0 \\
		\frac{1}{2} & 0 & 0 \\
		0 & \frac{3}{4} & 0 \\
		\end{bmatrix},$
		&
		$A^{\{\s,\s\}} =
		\begin{bmatrix}
		0 & 0 & 0 \\
		\frac{1}{2} & 0 & 0 \\
		0 & \frac{3}{4} & 0 \\
		\end{bmatrix},$
		&
		$A^{\{\f,\s,1\}}  =
		\begin{bmatrix}
		0 & 0 & 0 \\
		\frac{1}{2 M} & 0 & 0 \\
		0 & \frac{3}{4 M} & 0 \\
		\end{bmatrix},$ \\
		\fullrow{$A^{\{\f,\s,\lambda\}} =
			\begin{bmatrix}
			\frac{3 M^3-8 M^2+6 \lambda  M-6 \lambda +6}{6 (M-1) M} & \frac{-3 M^2+8 M-6}{6 (M-1)} & 0 \\
			\frac{-2 M^2+6 \lambda  M-3 M-6 \lambda +3}{6 (M-1) M} & \frac{M}{3 (M-1)} & 0 \\
			\frac{-3 M^3+2 M^2+12 \lambda  M-9 M-12 \lambda +12}{12 (M-1) M} & \frac{3 M^3-2 M^2+6 M-9}{12 (M-1) M} & 0 \\
			\end{bmatrix}, \qquad \lambda = 2, \ldots, M, $}  \\
		\fullrow{$A^{\{\s,\f,1\}} =  
			\begin{bmatrix}
			0 & 0 & 0 \\
			-\frac{1}{66} M (16 M-33) & \frac{8 M^2}{33} & 0 \\
			\frac{1}{264} \left(11 M^4-22 M^3+26 M^2+11 M+44\right) & \frac{1}{88} \left(-11 M^4+22 M^3-16 M^2-11 M+22\right) & \frac{1}{12} \left(M^4-2 M^3+M^2+M+4\right) \\
			\end{bmatrix},$} \\
		\fullrow{$A^{\{\s,\f,\lambda\}} =
			\begin{bmatrix}
			0 & 0 & 0 \\
			0 & 0 & 0 \\
			\frac{-M^4+2 M^3+2 M^2+3 M-4}{24 (M-1)} & \frac{1}{8} \left(M^3-M^2-M+2\right) & \frac{-M^4+2 M^3-M^2+3 M-4}{12 (M-1)} \\
			\end{bmatrix},\qquad \lambda = 2, \ldots, M,$} \\
		$b^{{\{\f\}}} =   b^{{\{\s\}}} =
		\begin{bmatrix}
		\frac{2}{9} & \frac{1}{3} & \frac{4}{9} \\
		\end{bmatrix}\trsym,$
		&
		$\widehat{b}^{{\{\f\}}}  =   \widehat{b}^{{\{\s\}}} =
		\begin{bmatrix}
		\frac{1}{40} & \frac{37}{40} & \frac{1}{20} \\
		\end{bmatrix}\trsym.$ 
		&
}

\begin{figure}[ht!]
\centering
\begin{subfigure}[b]{0.45\textwidth}
	\centering
	\includegraphics[width=\textwidth]{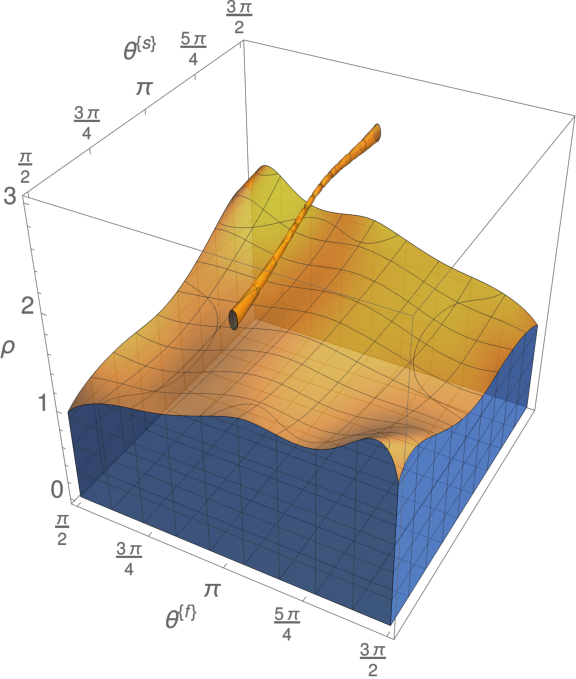}
	\caption{$M=2$.}
\end{subfigure}
\hfill
\begin{subfigure}[b]{0.45\textwidth}
	\centering
	\includegraphics[width=\textwidth]{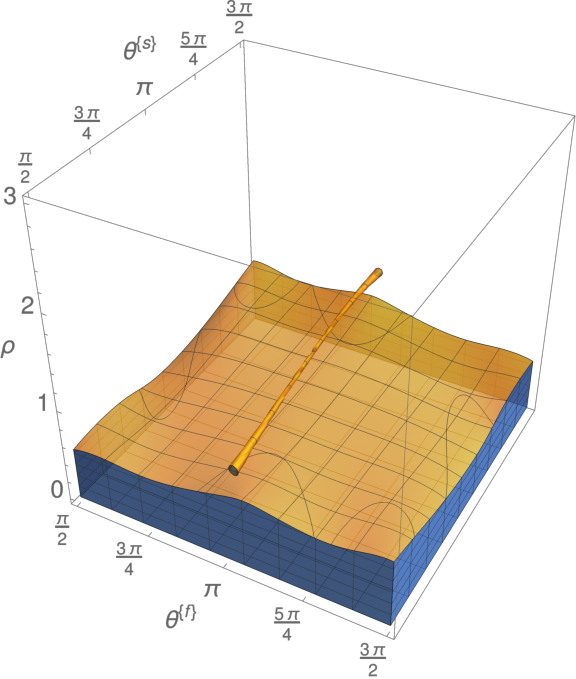}
	\caption{$M=4$.}
\end{subfigure}
\caption{\method{\explicit}{\explicit}{3}{2}{3}{3}{A} stability regions.}
\end{figure}

\newpage
\subsection{\mgark \method{\explicit}{\explicit}{3}{2}{4}{4}{A}}
\label{subsec:method_EXEX3a}
This explicit method is telescopic \cref{eqn:mr-telescopic-conditions} and naturally adaptive (\cref{def:naturally-adaptive}).

\fitbox{3}{
		$A^{\{\f,\f\}}=
		\begin{bmatrix}
		0 & 0 & 0 & 0 \\
		\frac{1}{3} & 0 & 0 & 0 \\
		0 & \frac{5}{9} & 0 & 0 \\
		\frac{833}{7680} & \frac{833}{9216} & \frac{3213}{5120} & 0 \\
		\end{bmatrix},$
		&
		$A^{\{\s,\s\}} =
		\begin{bmatrix}
		0 & 0 & 0 & 0 \\
		\frac{1}{3} & 0 & 0 & 0 \\
		0 & \frac{5}{9} & 0 & 0 \\
		\frac{833}{7680} & \frac{833}{9216} & \frac{3213}{5120} & 0 \\
		\end{bmatrix},$
		&
		$A^{\{\s,\f,\lambda\}} =
		\begin{bmatrix}
		0 & 0 & 0 & 0 \\
		0 & 0 & 0 & 0 \\
		0 & 0 & 0 & 0 \\
		0 & 0 & 0 & 0 \\
		\end{bmatrix},\qquad \lambda = 2, \ldots, M,$
		\\
		\fullrow{$A^{\{\f,\s,1\}} =
			\begin{bmatrix}
 0 & 0 & 0 & 0 \\
\frac{1}{3 M} & 0 & 0 & 0 \\
\frac{5 \left(518 M^3-2140 M^2+2399 M-777\right)}{2331 M (3 M-4)} & -\frac{5 \left(518 M^3-2140 M^2+1622 M+259\right)}{2331 M (3 M-4)} & 0 & 0 \\
\frac{17 \left(141932 M^3-445231 M^2+481160 M-178710\right)}{852480 M (3 M-4)} & -\frac{17 \left(94535 M^3-228442 M^2+142736 M-5180\right)}{340992 M (3 M-4)} &
\frac{3213 M}{5120} & 0 \\
			\end{bmatrix},$}  \\
		\fullrow{$A^{\{\f,\s,\lambda\}} =  
			\begin{bmatrix}
 \frac{\lambda -1}{M} & 0 & 0 & 0 \\
\frac{3 \lambda -2}{3 M} & 0 & 0 & 0 \\
\frac{777 (12 \lambda -7)+(6993 \lambda +12092) M^2-5965 M^3-3 (5439 \lambda +286) M}{2331 M \left(3 M^2-7 M+4\right)} & \frac{5 \left(1193 M^3-3040 M^2+1622
	M+259\right)}{2331 M \left(3 M^2-7 M+4\right)} & 0 & 0 \\
\frac{222 (3072 \lambda -335)+(511488 \lambda +1937719) M^2-867119 M^3-72 (16576 \lambda +13973) M}{170496 M \left(3 M^2-7 M+4\right)} & \frac{17 \left(51007 M^3-119207
	M^2+71368 M-2590\right)}{170496 M \left(3 M^2-7 M+4\right)} & 0 & 0 \\
			\end{bmatrix},\qquad \lambda = 2, \ldots, M,$} \\
		\fullrow{$A^{\{\s,\f,1\}} =
\begin{bmatrix}
	0 & 0 & 0 & 0 \\
	\frac{M}{3}-\frac{34 M^2}{361} & \frac{34 M^2}{361} & 0 & 0 \\
	0 & \frac{5 (1805-981 M) M}{6498} & \frac{5 M (327 M-361)}{2166} & 0 \\
	\frac{M \left(1480461 M^2-3944118 M+3007130\right)}{2772480} & -\frac{119 M \left(3249 M^2-20358 M+18050\right)}{3326976} & -\frac{119 M \left(66063 M^2-78954
		M-18050\right)}{5544960} & (M-1) M^2 \\
\end{bmatrix},$} \\
		$b^{{\{\f\}}} =   b^{{\{\s\}}} =
		\begin{bmatrix}
		\frac{101}{714} & \frac{1}{3} & \frac{1}{6} & \frac{128}{357} \\
		\end{bmatrix}\trsym,$
		&
		$\widehat{b}^{{\{\f\}}}  =   \widehat{b}^{{\{\s\}}} =
		\begin{bmatrix}
		\frac{7}{40} & -\frac{425}{8784} & \frac{100037}{131760} & \frac{188}{1647} \\
		\end{bmatrix}\trsym.$
		&
}

\begin{figure}[ht!]
\centering
\begin{subfigure}[b]{0.45\textwidth}
	\centering
	\includegraphics[width=\textwidth]{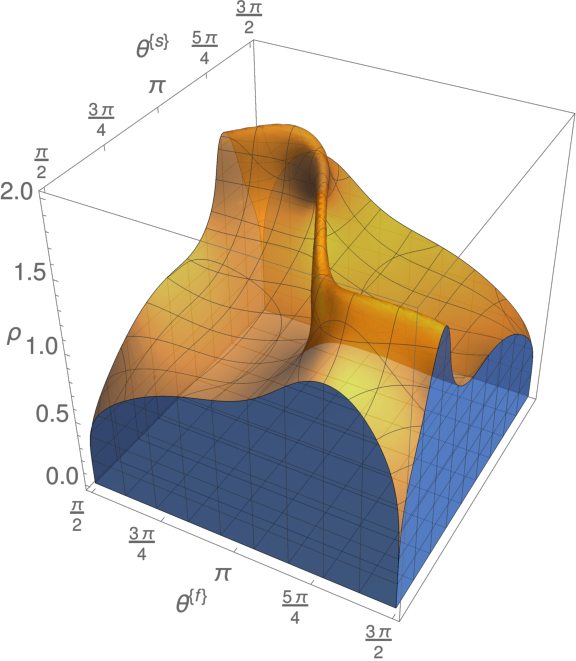}
	\caption{$M=2$.}
\end{subfigure}
\hfill
\begin{subfigure}[b]{0.45\textwidth}
	\centering
	\includegraphics[width=\textwidth]{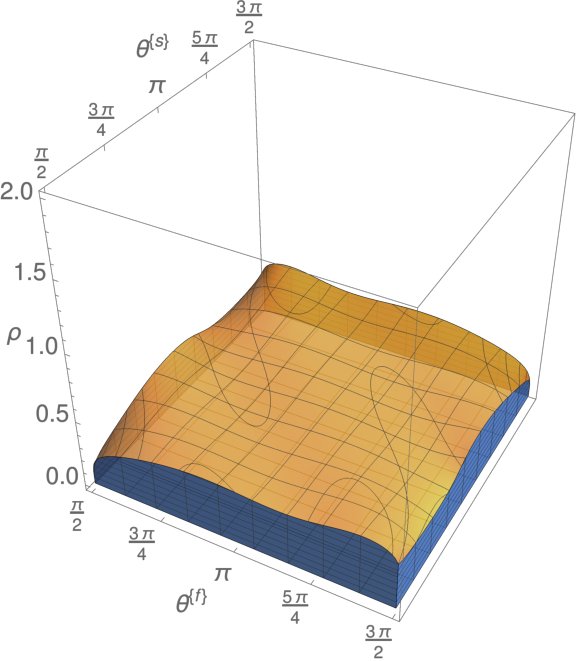}
	\caption{$M=4$.}
\end{subfigure}
\caption{\method{\explicit}{\explicit}{3}{2}{4}{4}{A} stability regions.}
\end{figure}
\newpage
\subsection{\mgark \method{\explicit}{\explicit}{3}{2}{3}{3}{S}}
\label{subsec:method_EXEX3s}
This explicit method uses a three stage base method and has telescopic \cref{eqn:mr-telescopic-conditions} property. Once again, we can take $L_2 = floor(c_2 M)$.

\fitbox{2}{
	$A^{\{\f,\f\}}=A^{\{\s,\s\}} =
	\begin{bmatrix}
	0 & 0 & 0 \\
	c_2 & 0 & 0 \\
	\frac{ \left(3 c_2^2-3 c_2+1\right)}{c_2 (3
		c_2-2)} & \frac{ (c_2-1)}{c_2 (3 c_2-2)} & 0
	\end{bmatrix}$,
	&
	\\
	$A^{\{\f,\s,\lambda\}} =
	\begin{bmatrix}
	\frac{\lambda-1}{M} & 0 & 0\\
	\frac{c_2 + \lambda-1}{M} & 0 & 0\\
	\frac{\lambda}{M} & 0 & 0
	\end{bmatrix},$
	& $\lambda = 1, \ldots, L_2, $  \\
	$A^{\{\f,\s,\lambda\}} =
	\begin{bmatrix}
	\frac{2 \lambda -1}{12 c_2 \left(L_2-M\right)}+\frac{1-2 \lambda }{12 c_2 \left(L_2+M\right)}+\frac{2 \lambda -1}{2 M} & \frac{2 \lambda -1}{12 c_2 \left(L_2+M\right)}+\frac{1-2 \lambda }{12 c_2 \left(L_2-M\right)}-\frac{1}{2 M}  & 0\\
	\frac{2 \lambda -1}{12 c_2 \left(L_2-M\right)}+\frac{1-2 \lambda }{12 c_2 \left(L_2+M\right)}+\frac{2 \lambda -1}{2 M} & \frac{2 \lambda -1}{12 c_2 \left(L_2+M\right)}+\frac{1-2 \lambda }{12 c_2 \left(L_2-M\right)}+\frac{2 c_2-1}{2 M} & 0 \\
	\frac{2 \lambda -1}{12 c_2 \left(L_2-M\right)}+\frac{1-2 \lambda }{12 c_2 \left(L_2+M\right)}+\frac{2 \lambda -1}{2 M} & \frac{2 \lambda -1}{12 c_2 \left(L_2+M\right)}+\frac{1-2 \lambda }{12 c_2 \left(L_2-M\right)}+\frac{1}{2 M}& 0
	\end{bmatrix},$
	& $\lambda = L_2 + 1, \ldots, M, $  \\
	$A^{\{\s,\f,\lambda\}} =
	\begin{bmatrix}
	0 & 0 & 0 \\
	c_2 \left(\frac{2 \lambda  \left(c_2 \left(4 L_2-3\right)-3 L_2+3\right)}{\left(c_2-1\right) \left(3 c_2^2+4 c_2+1\right) \left(L_2+1\right)}+\frac{M}{L_2}\right) & -\frac{c_2 \lambda  \left(c_2 \left(4 L_2-3\right)-3 L_2+3\right)}{\left(c_2-1\right) \left(3 c_2^2+4 c_2+1\right) \left(L_2+1\right)}  & -\frac{c_2 \lambda  \left(c_2 \left(4 L_2-3\right)-3 L_2+3\right)}{\left(c_2-1\right) \left(3 c_2^2+4 c_2+1\right) \left(L_2+1\right)}  \\
	\frac{2 c_2 \lambda  \left(c_2 \left(4 L_2-3\right)-3 L_2+3\right)}{\left(c_2-1\right) \left(3 c_2^2+4 c_2+1\right) \left(L_2+1\right)} &  -\frac{c_2 \lambda  \left(c_2 \left(4 L_2-3\right)-3 L_2+3\right)}{\left(c_2-1\right) \left(3 c_2^2+4 c_2+1\right) \left(L_2+1\right)} & -\frac{c_2 \lambda  \left(c_2 \left(4 L_2-3\right)-3 L_2+3\right)}{\left(c_2-1\right) \left(3 c_2^2+4 c_2+1\right) \left(L_2+1\right)}  
	\end{bmatrix},$
	&
	$\lambda = 1, \ldots, L_2,$  \\
	$A^{\{\s,\f,\lambda\}} =
	\begin{bmatrix}
	0 & 0 & 0\\
	0 & 0 & 0\\
	\frac{\lambda }{3 c_2-2}+\frac{c_2 \left(3 L_2-4\right)-3 L_2+3}{6 c_2-4}+\frac{M}{M-L_2} & \frac{\lambda }{2-3 c_2} & \frac{c_2 \left(4-3 L_2\right)+3 \left(L_2-1\right)}{6 c_2-4}
	\end{bmatrix},$
	&
	$\lambda = L_2 + 1, \ldots, M,$  \\
	$b^{{\{\f\}}} =   b^{{\{\s\}}} =
	\begin{bmatrix}
	\frac{3 c_2-1}{6 c_2} & -\frac{1}{6 \left(c_2-1\right) c_2} & \frac{3 c_2-2}{6 \left(c_2-1\right)} \\
	\end{bmatrix}\trsym,$ \quad
	$\widehat{b}^{{\{\f\}}}  =   \widehat{b}^{{\{\s\}}} =
	\begin{bmatrix}
	\hat{b}_2 \left(c_2-1\right)+\frac{1}{2} &
	\hat{b}_2 &
	\frac{1-2 \hat{b}_2 c_2}{2} \\
	\end{bmatrix}\trsym.$
	&
}
\begin{figure}[ht!]
	\centering
	\begin{subfigure}[b]{0.45\textwidth}
		\centering
		\includegraphics[width=\textwidth]{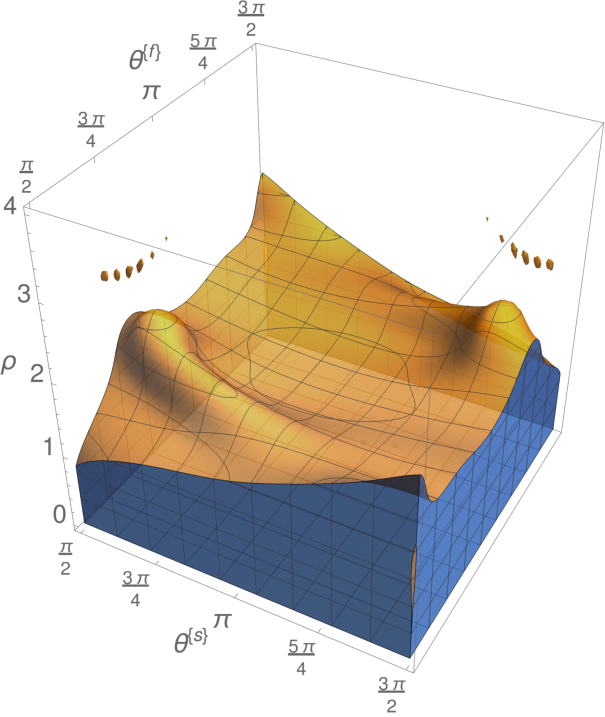}
		\caption{$M=2$.}
	\end{subfigure}
	\hfill
	\begin{subfigure}[b]{0.45\textwidth}
		\centering
		\includegraphics[width=\textwidth]{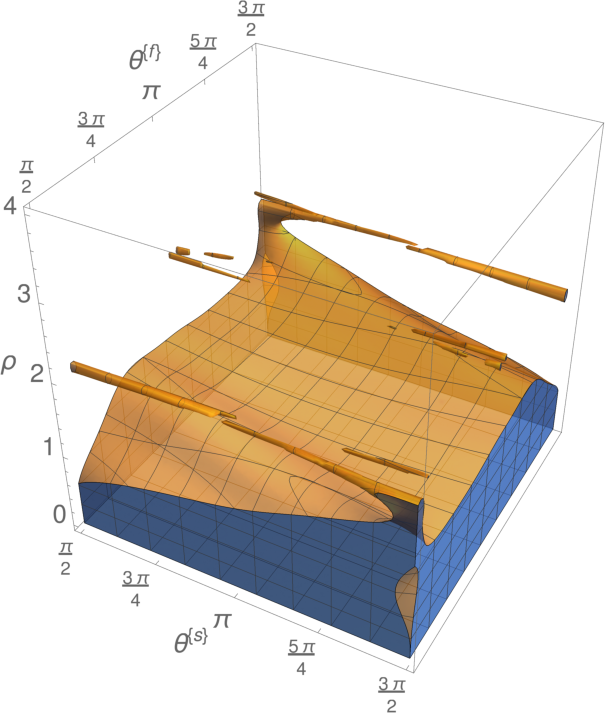}
		\caption{$M=4$.}
	\end{subfigure}
	\hfill
	\caption{Stability plots for \mgark \method{\explicit}{\explicit}{3}{2}{3}{3}{S} method with $c_2 = \frac{1}{2}$ . }
\end{figure}
\newpage
\subsection{\mgark \method{\explicit}{\implicit}{3}{2}{3}{3}{A}}
\label{subsec:method_EXIM3}
This explcit method uses a fast method from \cite{ralston1962runge} and a slow method from \cite{Alexander_1977_SDIRK} and is stiffly accurate \eqref{eqn:stiff-accuracy} in the slow partition.

\fitbox{2}{
		$A^{\{\f,\f\}}= 
		\begin{bmatrix}
		0 & 0 & 0 \\
		\frac{1}{2} & 0 & 0 \\
		0 & \frac{3}{4} & 0 \\
		\end{bmatrix},$
		&
		$A^{\{\s,\s\}} = 
		\begin{bmatrix}
		\gamma  & 0 & 0 \\
		-\frac{2 \left(3 \gamma ^3-9 \gamma ^2+6 \gamma -1\right)}{3 \left(2 \gamma ^2-4 \gamma +1\right)} & \gamma  & 0 \\
		\frac{4 \gamma -1}{4 \left(3 \gamma ^3-9 \gamma ^2+6 \gamma -1\right)} & -\frac{3 \left(2 \gamma ^2-4 \gamma +1\right)^2}{4 \left(3 \gamma ^3-9 \gamma ^2+6 \gamma -1\right)} & \gamma  \\
		\end{bmatrix},$ \\
		\fullrow{$A^{\{\f,\s,\lambda\}} =
		\begin{bmatrix}
		\frac{\lambda -1}{M} & 0 & 0 \\
		\frac{2 \lambda -1}{2 M} & 0 & 0 \\
		\frac{-60 \lambda  \gamma ^3+42 \gamma ^3+18 M \gamma ^2+72 \lambda  \gamma ^2-72 \gamma ^2-36 M \gamma +42 \lambda  \gamma +3 \gamma +9 M-16 \lambda +4}{16 M \left(3 \gamma ^3-9 \gamma ^2+6 \gamma -1\right)} & -\frac{9 \left(2 \gamma ^2-4 \gamma +1\right) (M+3 \gamma -6 \gamma  \lambda )}{16 M \left(3 \gamma ^3-9 \gamma ^2+6 \gamma -1\right)} & 0 \\
		\end{bmatrix}, \qquad \lambda = 1, \ldots, M,$} \\
		\fullrow{$A^{\{\s,\f,1\}}=
		\begin{bmatrix}
		M \gamma  & 0 & 0 \\
		-\frac{M \left(36 M \gamma ^4-36 \gamma ^4-120 M \gamma ^3+126 \gamma ^3+108 M \gamma ^2-138 \gamma ^2-36 M \gamma +51 \gamma +4 M-6\right)}{9 \left(2 \gamma ^2-4 \gamma +1\right)^2} & \frac{4 M^2 \left(9 \gamma ^4-30 \gamma ^3+27 \gamma ^2-9 \gamma +1\right)}{9 \left(2 \gamma ^2-4 \gamma +1\right)^2} & 0 \\
		\frac{2}{9} & \frac{1}{3} & \frac{4}{9} \\
		\end{bmatrix},$} \\
		\fullrow{$A^{\{\s,\f,\lambda\}} =
		\begin{bmatrix}
		0 & 0 & 0 \\
		0 & 0 & 0 \\
		\frac{2}{9} & \frac{1}{3} & \frac{4}{9} \\
		\end{bmatrix}, \qquad \lambda = 2, \ldots, M,$}\\
		${b}^{{\{\f\}}}= 
		\begin{bmatrix}
		\frac{2}{9} & \frac{1}{3} & \frac{4}{9} \\
		\end{bmatrix}\trsym,$
		&
		$b^{{\{\s\}}} = 
		\begin{bmatrix}
		\frac{4 \gamma -1}{4 \left(3 \gamma ^3-9 \gamma ^2+6 \gamma -1\right)} & -\frac{3 \left(2 \gamma ^2-4 \gamma +1\right)^2}{4 \left(3 \gamma ^3-9 \gamma ^2+6 \gamma -1\right)} & \gamma  \\
		\end{bmatrix}\trsym,$ \\
		$\widehat{b}^{{\{\f\}}}=
		\begin{bmatrix}
		\frac{1}{40} & \frac{37}{40} & \frac{1}{20} \\
		\end{bmatrix}\trsym,$
		&
		$\widehat{b}^{{\{\s\}}} =
		\begin{bmatrix}
		\frac{-6 \gamma ^2+6 \gamma -1}{4 \left(3 \gamma ^3-9 \gamma ^2+6 \gamma -1\right)} & \frac{3 \left(4 \gamma ^3-10 \gamma ^2+6 \gamma -1\right)}{4 \left(3 \gamma ^3-9 \gamma ^2+6 \gamma -1\right)} & 0 \\
		\end{bmatrix}\trsym,$ \\
		\fullrow{$\gamma = \frac{1}{2} \left(2+\sqrt{6} \sin \left(\frac{1}{3} \cot ^{-1}{2 \sqrt{2}}\right)-\sqrt{2} \cos \left(\frac{1}{3} \cot ^{-1}{2 \sqrt{2}}\right)\right) \approx 0.43586652150845899942.$}
}

\begin{figure}[ht!]
\centering
\begin{subfigure}[b]{0.45\textwidth}
	\centering
	\includegraphics[width=\textwidth]{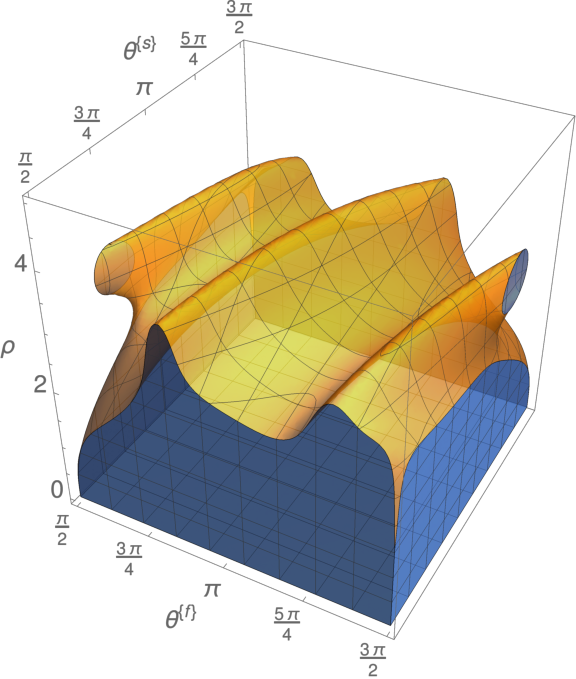}
	\caption{$M=2$.}
\end{subfigure}
\hfill
\begin{subfigure}[b]{0.45\textwidth}
	\centering
	\includegraphics[width=\textwidth]{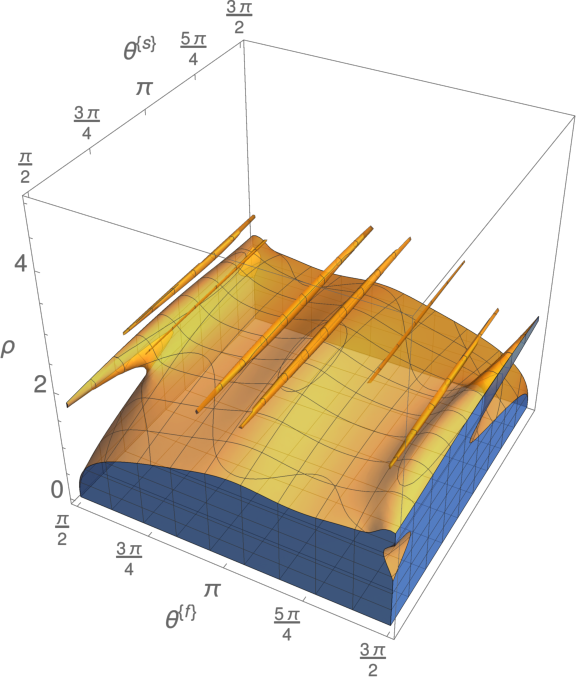}
	\caption{$M=4$.}
\end{subfigure}
\caption{\method{\explicit}{\implicit}{3}{2}{3}{3}{A} stability regions.}
\end{figure}

\newpage
\subsection{\mgark \method{\implicit}{\explicit}{3}{2}{3}{3}{A}}
\label{subsec:method_IMEX3}
This implicit-explicit method uses a fast method from \cite{Alexander_1977_SDIRK} and a slow method from \cite{ralston1962runge} and is stiffly accurate \eqref{eqn:stiff-accuracy} in the fast partition.

\fitbox{2}{
		$A^{\{\f,\f\}}= 
		\begin{bmatrix}
		\gamma  & 0 & 0 \\
		-\frac{2 \left(3 \gamma ^3-9 \gamma ^2+6 \gamma -1\right)}{3 \left(2 \gamma ^2-4 \gamma +1\right)} & \gamma  & 0 \\
		\frac{4 \gamma -1}{4 \left(3 \gamma ^3-9 \gamma ^2+6 \gamma -1\right)} & -\frac{3 \left(2 \gamma ^2-4 \gamma +1\right)^2}{4 \left(3 \gamma ^3-9 \gamma ^2+6 \gamma -1\right)} & \gamma  \\
		\end{bmatrix},$
		&
		$A^{\{\s,\s\}}=
		\begin{bmatrix}
		0 & 0 & 0 \\
		\frac{1}{2} & 0 & 0 \\
		0 & \frac{3}{4} & 0 \\
		\end{bmatrix},$ \\
		\fullrow{$A^{\{\f,\s,\lambda\}}=
		\begin{bmatrix}
		\frac{\gamma +\lambda -1}{M} & 0 & 0 \\
		\frac{6 \lambda  \gamma ^2-12 \lambda  \gamma +3 \gamma +3 \lambda -1}{3 M \left(2 \gamma ^2-4 \gamma +1\right)} & 0 & 0 \\
		\frac{\lambda }{M} & 0 & 0 \\
		\end{bmatrix}, \qquad \lambda = 1, \ldots, M-1,$} \\
		\fullrow{$A^{\{\f,\s,M\}}=
			\begin{bmatrix}
			\frac{M+\gamma -1}{M} & 0 & 0 \\
			\frac{12 M^2 \gamma ^3-36 M \gamma ^3+18 \gamma ^3-36 M^2 \gamma ^2+108 M \gamma ^2-42 \gamma ^2+24 M^2 \gamma -60 M \gamma +21 \gamma -4 M^2+9 M-3}{9 M \left(2 \gamma ^2-4 \gamma +1\right)^2} & -\frac{4 (M-3 \gamma ) \left(3 \gamma ^3-9 \gamma ^2+6 \gamma -1\right)}{9 \left(2 \gamma ^2-4 \gamma +1\right)^2} & 0 \\
			\frac{2}{9} & \frac{1}{3} & \frac{4}{9} \\
			\end{bmatrix},$} \\
		\fullrow{$A^{\{\s,\f,\lambda\}} = 
		\begin{bmatrix}
		0 & 0 & 0 \\
		\frac{1}{2} & 0 & 0 \\
		-\frac{3 \left(12 \gamma ^3+6 M \gamma ^2-18 \gamma ^2-12 M \gamma +6 \gamma +3 M-1\right)}{32 \left(3 \gamma ^3-9 \gamma ^2+6 \gamma -1\right)} & \frac{9 (M+6 \gamma -3) \left(2 \gamma ^2-4 \gamma +1\right)}{32 \left(3 \gamma ^3-9 \gamma ^2+6 \gamma -1\right)} & 0 \\
		\end{bmatrix}, \qquad \lambda = 1, \ldots, M,$} \\
		${b}^{{\{\f\}}}=  
		\begin{bmatrix}
		\frac{4 \gamma -1}{4 \left(3 \gamma ^3-9 \gamma ^2+6 \gamma -1\right)} & -\frac{3 \left(2 \gamma ^2-4 \gamma +1\right)^2}{4 \left(3 \gamma ^3-9 \gamma ^2+6 \gamma -1\right)} & \gamma  \\
		\end{bmatrix}\trsym,$
		&
		$b^{{\{\s\}}} =
		\begin{bmatrix}
		\frac{2}{9} & \frac{1}{3} & \frac{4}{9} \\
		\end{bmatrix}\trsym,$
		\\
		$\widehat{b}^{{\{\f\}}}= 
		\begin{bmatrix}
		\frac{-6 \gamma ^2+6 \gamma -1}{4 \left(3 \gamma ^3-9 \gamma ^2+6 \gamma -1\right)} & \frac{3 \left(4 \gamma ^3-10 \gamma ^2+6 \gamma -1\right)}{4 \left(3 \gamma ^3-9 \gamma ^2+6 \gamma -1\right)} & 0 \\
		\end{bmatrix}\trsym,$
		&
		$\widehat{b}^{{\{\s\}}} =
		\begin{bmatrix}
		\frac{1}{40} & \frac{37}{40} & \frac{1}{20} \\
		\end{bmatrix}\trsym.$ \\
		\fullrow{$\gamma = \frac{1}{2} \left(2+\sqrt{6} \sin \left(\frac{1}{3} \cot ^{-1}{2 \sqrt{2}}\right)-\sqrt{2} \cos \left(\frac{1}{3} \cot ^{-1}{2 \sqrt{2}}\right)\right) \approx 0.43586652150845899942.$}
}

\begin{figure}[ht!]
\centering
\begin{subfigure}[b]{0.45\textwidth}
	\centering
	\includegraphics[width=\textwidth]{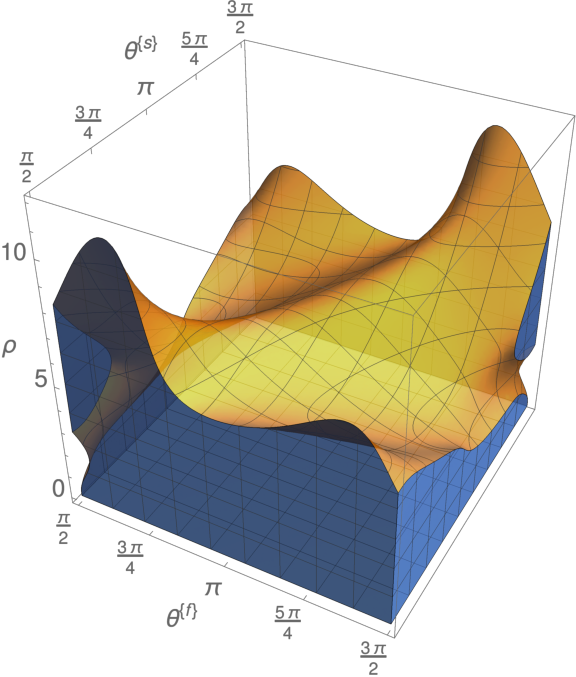}
	\caption{$M=2$.}
\end{subfigure}
\hfill
\begin{subfigure}[b]{0.45\textwidth}
	\centering
	\includegraphics[width=\textwidth]{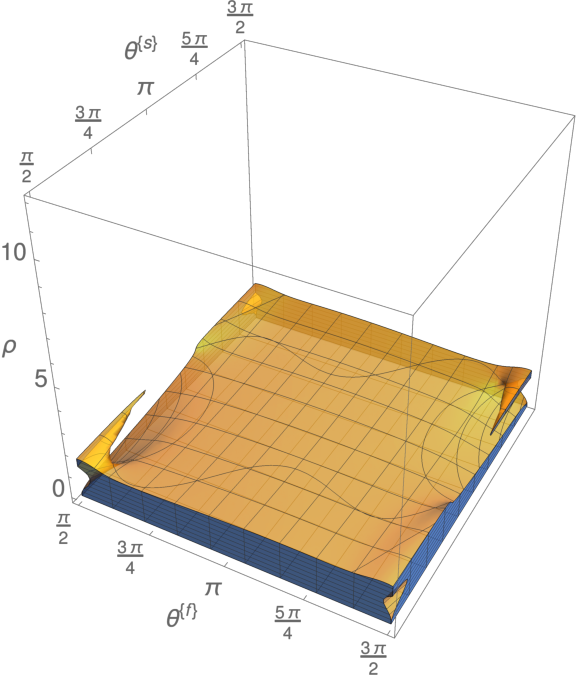}
	\caption{$M=4$.}
\end{subfigure}
\caption{\method{\implicit}{\explicit}{3}{2}{3}{3}{A} stability regions.}
\end{figure}

\newpage
\subsection{\mgark \method{\explicit}{\explicit}{4}{3}{5}{5}{A}}
\label{subsec:method_EXEX4}
This explicit method uses a base method from \cite{sofroniou2004construction}, is telescopic  \cref{eqn:mr-telescopic-conditions} and has the First same as last (FSAL) property.

\fittboxrotate{3}{
		$A^{\{\f,\f\}}=
		\begin{bmatrix}
		0 & 0 & 0 & 0 & 0 \\
		\frac{2}{5} & 0 & 0 & 0 & 0 \\
		-\frac{3}{20} & \frac{3}{4} & 0 & 0 & 0 \\
		\frac{19}{44} & -\frac{15}{44} & \frac{10}{11} & 0 & 0 \\
		\frac{11}{72} & \frac{25}{72} & \frac{25}{72} & \frac{11}{72} & 0 \\
		\end{bmatrix},$
		&
		$A^{\{\s,\s\}} =
		\begin{bmatrix}
		0 & 0 & 0 & 0 & 0 \\
		\frac{2}{5} & 0 & 0 & 0 & 0 \\
		-\frac{3}{20} & \frac{3}{4} & 0 & 0 & 0 \\
		\frac{19}{44} & -\frac{15}{44} & \frac{10}{11} & 0 & 0 \\
		\frac{11}{72} & \frac{25}{72} & \frac{25}{72} & \frac{11}{72} & 0 \\
		\end{bmatrix},$
		&
		$A^{\{\s,\f,\lambda\}} =
		\begin{bmatrix}
		0 & 0 & 0 & 0 & 0 \\
		0 & 0 & 0 & 0 & 0 \\
		0 & 0 & 0 & 0 & 0 \\
		0 & 0 & 0 & 0 & 0 \\
		\frac{11}{72} & \frac{25}{72} & \frac{25}{72} & \frac{11}{72} & 0 \\
		\end{bmatrix},\qquad \lambda = 2, \ldots, M,$
		\\
		\fullrow{$A^{\{\f,\s,1\}} =
			\begin{bmatrix}
			0 & 0 & 0 & 0 & 0 \\
			\frac{2}{5 M} & 0 & 0 & 0 & 0 \\
			\frac{3 \left(10 M^3-30 M^2+22 M-1\right)}{20 M (3 M-4)} & -\frac{3 \left(2 M^3-6 M^2+2 M+3\right)}{4 M (3 M-4)} & 0 & 0 & 0 \\
			0 & \frac{3 \left(10 M^3-50 M^2+116 M-83\right)}{22 M (3 M-4)} & \frac{-30 M^3+150 M^2-282 M+161}{22 M (3 M-4)} & 0 & 0 \\
			\frac{11}{72 M} & \frac{25}{72 M} & \frac{25}{72 M} & \frac{11}{72 M} & 0 \\
			\end{bmatrix},$}  \\
		\fullrow{$A^{\{\s,\f,\lambda\}} =  
			\begin{bmatrix}
			\frac{11 (\lambda -1)}{72 M} & \frac{25 (\lambda -1)}{72 M} & \frac{25 (\lambda -1)}{72 M} & \frac{11 (\lambda -1)}{72 M} & 0 \\
			\frac{-450 M^2+956 \lambda  M-497 M-956 \lambda +776}{450 (M-1) M} & \frac{450 M^2-506 \lambda  M+227 M+506 \lambda -506}{450 (M-1) M} & 0 & 0 & 0 \\
			\frac{-900 M^3+1239 \lambda  M^2+2217 M^2-2891 \lambda  M-97 M+1652 \lambda -1562}{600 M \left(3 M^2-7 M+4\right)} & \frac{900 M^3+561 \lambda  M^2-2937 M^2-1309 \lambda  M+1777 M+748 \lambda +602}{600 M \left(3 M^2-7 M+4\right)} & 0 & 0 & 0 \\
			0 & \frac{-90 M^3+99 \lambda  M^2+197 M^2-231 \lambda  M-205 M+132 \lambda +117}{22 M \left(3 M^2-7 M+4\right)} & \frac{3240 M^3-825 \lambda  M^2-7455 M^2+1925 \lambda  M+8227 M-1100 \lambda -4696}{792 M \left(3 M^2-7 M+4\right)} & -\frac{11 (\lambda -1)}{72 M} & 0 \\
			\frac{11 \lambda }{72 M} & \frac{25 \lambda }{72 M} & \frac{25 \lambda }{72 M} & \frac{11 \lambda }{72 M} & 0 \\
			\end{bmatrix},\qquad \lambda = 2, \ldots, M,$} \\
		\fullrow{$A^{\{\s,\f,1\}} =
			\begin{bmatrix}
			0 & 0 & 0 & 0 & 0 \\
			\frac{2 M}{5} & 0 & 0 & 0 & 0 \\
			-\frac{3}{20} M (5 M-4) & \frac{3 M^2}{4} & 0 & 0 & 0 \\
			\frac{1}{44} M \left(56 M^2-81 M+44\right) & -\frac{5}{44} M^2 (16 M-13) & -\frac{5}{11} (M-3) M^2 & (M-1) M^2 & 0 \\
			\frac{11}{72} & \frac{25}{72} & \frac{25}{72} & \frac{11}{72} & 0 \\
			\end{bmatrix},$} \\
		$b^{{\{\f\}}} =   b^{{\{\s\}}} =
		\begin{bmatrix}
		\frac{11}{72} & \frac{25}{72} & \frac{25}{72} & \frac{11}{72} & 0 \\
		\end{bmatrix}\trsym,$
		&
		$\widehat{b}^{{\{\f\}}}  =   \widehat{b}^{{\{\s\}}} =
		\begin{bmatrix}
		\frac{1251515}{8970912} & \frac{3710105}{8970912} & \frac{2519695}{8970912} & \frac{61105}{8970912} &
		\frac{119041}{747576} \\
		\end{bmatrix}\trsym.$
		&
}

\begin{figure}[ht!]
\centering
\begin{subfigure}[b]{0.45\textwidth}
	\centering
	\includegraphics[width=\textwidth]{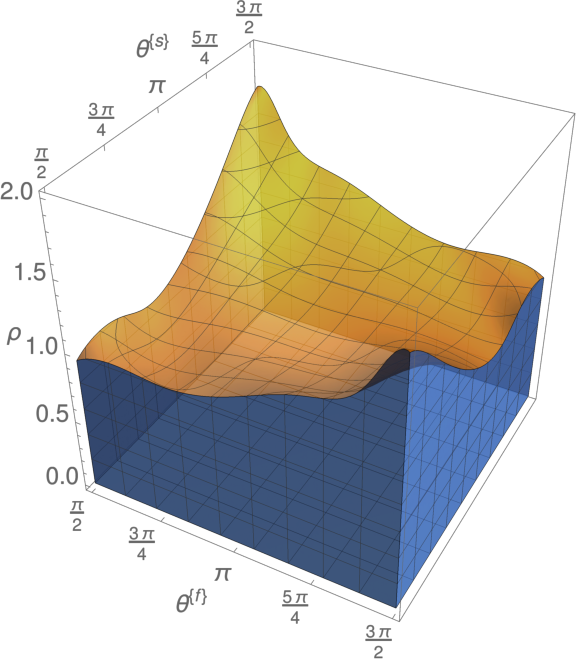}
	\caption{$M=2$.}
\end{subfigure}
\hfill
\begin{subfigure}[b]{0.45\textwidth}
	\centering
	\includegraphics[width=\textwidth]{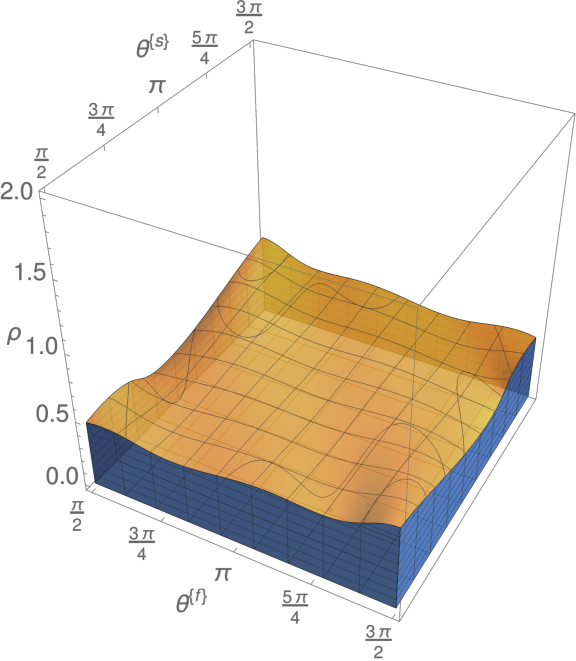}
	\caption{$M=4$.}
\end{subfigure}
\caption{\method{\explicit}{\explicit}{4}{3}{5}{5}{A} stability regions.}
\end{figure}

\newpage
\subsection{\mgark \method{\explicit}{\implicit}{4}{3}{6}{5}{A}}
\label{subsec:method_EXIM4}
This explicit-implicit method uses a fast method from \cite{Fehlberg_1969_construction} and a slow method from \cite{Kennedy_2016_SDIRK-review}. The multirate scheme is stiffly accurate \eqref{eqn:stiff-accuracy} in the slow partition.

\fittboxrotate{3}{
		$A^{\{\f,\f\}}= 
		\begin{bmatrix}
		0 & 0 & 0 & 0 & 0 & 0 \\
		\frac{1}{4} & 0 & 0 & 0 & 0 & 0 \\
		\frac{3}{32} & \frac{9}{32} & 0 & 0 & 0 & 0 \\
		\frac{1932}{2197} & -\frac{7200}{2197} & \frac{7296}{2197} & 0 & 0 & 0 \\
		\frac{439}{216} & -8 & \frac{3680}{513} & -\frac{845}{4104} & 0 & 0 \\
		-\frac{8}{27} & 2 & -\frac{3544}{2565} & \frac{1859}{4104} & -\frac{11}{40} & 0 \\
		\end{bmatrix},$
		&
		$A^{\{\s,\s\}} = 
		\begin{bmatrix}
		\frac{1}{4} & 0 & 0 & 0 & 0 \\
		\frac{13}{20} & \frac{1}{4} & 0 & 0 & 0 \\
		\frac{580}{1287} & -\frac{175}{5148} & \frac{1}{4} & 0 & 0 \\
		\frac{12698}{37375} & -\frac{201}{2990} & \frac{891}{11500} & \frac{1}{4} & 0 \\
		\frac{944}{1365} & -\frac{400}{819} & \frac{99}{35} & -\frac{575}{252} & \frac{1}{4} \\
		\end{bmatrix},$
		&
		$A^{\{\s,\f,1\}}=
		\begin{bmatrix}
		\frac{M}{4} & 0 & 0 & 0 & 0 & 0 \\
		-\frac{1}{100} M (169 M-90) & \frac{169 M^2}{100} & 0 & 0 & 0 & 0 \\
		-\frac{1}{198} M (155 M-132) & \frac{155 M^2}{198} & 0 & 0 & 0 & 0 \\
		-\frac{1}{920} M (497 M-552) & \frac{14 M^2}{23} & -\frac{896 M^2}{10925} & \frac{1183 M^2}{87400} & 0 & 0 \\
		\frac{25}{216} & 0 & \frac{1408}{2565} & \frac{2197}{4104} & -\frac{1}{5} & 0 \\
		\end{bmatrix},$ 
		\\
		\fullrow{$A^{\{\f,\s,\lambda\}} =
		\begin{bmatrix}
		\frac{\lambda -1}{M} & 0 & 0 & 0 & 0 \\
		\frac{4 \lambda -3}{4 M} & 0 & 0 & 0 & 0 \\
		\frac{45 M^3-90 M^2+551 \lambda  M-335 M+90 \lambda -90}{416 M^2} & -\frac{15 \left(3 M^3-6 M^2+9 \lambda  M-5 M+6 \lambda -6\right)}{416 M^2} & 0 & 0 & 0 \\
		\frac{1440 M^3-2880 M^2+6517 \lambda  M-3709 M-3960 \lambda +3960}{2197 M^2} & -\frac{60 \left(24 M^3-48 M^2+72 \lambda  M-59 M-66 \lambda +66\right)}{2197 M^2} & 0 & 0 & 0 \\
		\frac{560 M^3-362 M^2+386 \lambda  M-529 M-1155 \lambda +1155}{273 M^2} & -\frac{5 \left(672 M^3-2439 M^2+4046 \lambda  M-2590 M-1386 \lambda +1386\right)}{1638 M^2} & -\frac{33 (11 M-14 \lambda +7)}{28 M} & \frac{575 (3 M-2 \lambda +1)}{252 M} & 0 \\
		0 & 0 & \frac{160 M^3-109 M^2-300 M+165}{32 M^2} & \frac{-160 M^3+109 M^2+32 \lambda  M+284 M-165}{32 M^2} & 0 \\
		\end{bmatrix}, \qquad \lambda = 1, \ldots, M,$} \\
		\fullrow{$A^{\{\s,\f,\lambda\}} =
		\begin{bmatrix}
		0 & 0 & 0 & 0 & 0 & 0 \\
		0 & 0 & 0 & 0 & 0 & 0 \\
		0 & 0 & 0 & 0 & 0 & 0 \\
		0 & 0 & 0 & 0 & 0 & 0 \\
		\frac{25}{216} & 0 & \frac{1408}{2565} & \frac{2197}{4104} & -\frac{1}{5} & 0 \\
		\end{bmatrix}, \qquad \lambda = 2, \ldots, M,$}	\\
		${b}^{{\{\f\}}}= 
		\begin{bmatrix}
		\frac{25}{216} & 0 & \frac{1408}{2565} & \frac{2197}{4104} & -\frac{1}{5} & 0 \\
		\end{bmatrix}\trsym,$
		&
		$b^{{\{\s\}}} = 
		\begin{bmatrix}
		\frac{944}{1365} & -\frac{400}{819} & \frac{99}{35} & -\frac{575}{252} & \frac{1}{4} \\
		\end{bmatrix}\trsym,$
		&
		\\
		$\widehat{b}^{{\{\f\}}}=
		\begin{bmatrix}
		\frac{16}{135} & 0 & \frac{6656}{12825} & \frac{28561}{56430} & -\frac{9}{50} & \frac{2}{55} \\
		\end{bmatrix}\trsym,$
		&
		$\widehat{b}^{{\{\s\}}} =
		\begin{bmatrix}
		\frac{41911}{60060} & -\frac{83975}{144144} & \frac{3393}{1120} & -\frac{27025}{11088} & \frac{103}{352} \\
		\end{bmatrix}\trsym.$   
		&
}

\begin{figure}[ht!]
\centering
\begin{subfigure}[b]{0.45\textwidth}
	\centering
	\includegraphics[width=\textwidth]{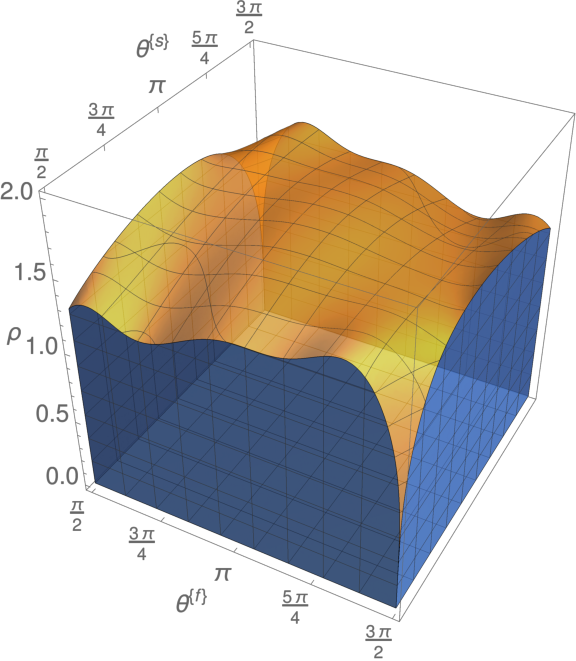}
	\caption{$M=2$.}
\end{subfigure}
\hfill
\begin{subfigure}[b]{0.45\textwidth}
	\centering
	\includegraphics[width=\textwidth]{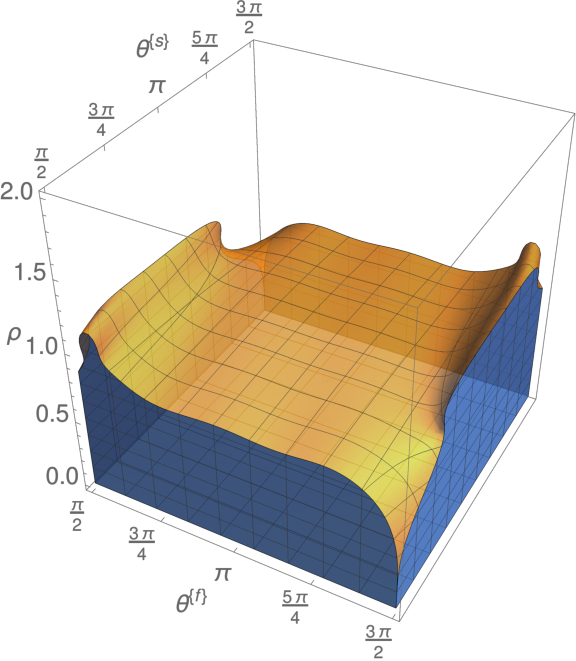}
	\caption{$M=4$.}
\end{subfigure}
\caption{\method{\explicit}{\implicit}{4}{3}{6}{5}{A} stability regions.}
\end{figure}

\newpage
\subsection{\mgark \method{\implicit}{\explicit}{4}{2}{6}{4}{A}}
\label{subsec:method_IMEX4}
This implicit-explicit method uses a slow method from \cite{sofroniou2004construction} and is stiffly accurate \eqref{eqn:stiff-accuracy} in the fast partition.

\fittboxrotate{2}{
		$A^{\{\f,\f\}}= 
		\begin{bmatrix}
		\frac{191}{1000} & 0 & 0 & 0 & 0 & 0 \\
		\frac{209}{1000} & \frac{191}{1000} & 0 & 0 & 0 & 0 \\
		\frac{8466728223}{12920014250} & -\frac{12729769579}{51680057000} & \frac{191}{1000} & 0 & 0 & 0 \\
		\frac{102093693512533448034070599559771}{222819131395744425631166002057000} & -\frac{17248151203963882893894684614}{68098756539041694875050734125} & \frac{783289327941232988291717301}{1938400447113914098574736860} & \frac{191}{1000} & 0 & 0 \\
		\frac{1837041228720545025825201951582239534326}{2195453146940870392428577778808091404375} & -\frac{12181532573386077454382848427541846123}{17628427274186874088578235825358750000} & \frac{4528991149246665992465958589958885289}{8624433917634487442857055565277187500} & \frac{83750160542686187}{606298988321250000} & \frac{191}{1000} & 0 \\
		\frac{2288000}{4732539} & -\frac{2203}{14250} & \frac{247273}{613500} & \frac{30767}{152250} & -\frac{1}{8} & \frac{191}{1000} \\
		\end{bmatrix},$
		&
		$A^{\{\s,\s\}}=
		\begin{bmatrix}
		0 & 0 & 0 & 0 \\
		\frac{2}{5} & 0 & 0 & 0 \\
		-\frac{3}{20} & \frac{3}{4} & 0 & 0 \\
		\frac{19}{44} & -\frac{15}{44} & \frac{10}{11} & 0 \\
		\end{bmatrix},$ \\
		\fullrow{$A^{\{\f,\s,\lambda\}}=
			\begin{bmatrix}
			\frac{1000 \lambda -809}{1000 M} & 0 & 0 & 0 \\
			\frac{5 \lambda -3}{5 M} & 0 & 0 & 0 \\
			\frac{5 \lambda -2}{5 M} & 0 & 0 & 0 \\
			\frac{5 \lambda -1}{5 M} & 0 & 0 & 0 \\
			\frac{\lambda }{M} & 0 & 0 & 0 \\
			\frac{\lambda }{M} & 0 & 0 & 0 \\
			\end{bmatrix}, \qquad \lambda = 1, \ldots, M-1,$} \\
		\fullrow{$A^{\{\f,\s,M\}}=
			\begin{bmatrix}
			\frac{1000 M-809}{1000 M} & 0 & 0 & 0 \\
			\frac{11380195070453 M^3-18408895671188 M^2+14477055081282 M-5016867120000}{8361445200000 M} & -\frac{209 \left(54450694117 M^2-88080840532 M+29261291298\right)}{8361445200000} & 0 & 0 \\
			\frac{-45728475609635251 M^3-421177045491040004 M^2+701106234145018506 M-206755963893960000}{516889909734900000 M} & \frac{409 \left(111805563837739 M^2+1029772727361956 M-450406661149434\right)}{516889909734900000} & 0 & 0 \\
			\frac{313252304037186017 M^3-457232580001772932 M^2+265208779590977398 M-51451316893680000}{257256584468400000 M} & -\frac{203 \left(2350256923212739 M^2-2188535491388044 M-533759817986934\right)}{257256584468400000} & \frac{203 \left(885000 M^2+70000 M-628199\right)}{282071856} & 0 \\
			0 & \frac{-885000 M^2+885000 M+831041}{859500} & \frac{590000 M^2-590000 M-114791}{573000} & \frac{2101}{9000} \\
			\frac{11}{72} & \frac{25}{72} & \frac{25}{72} & \frac{11}{72} \\
			\end{bmatrix},$} \\
		\fullrow{$A^{\{\s,\f,\lambda\}} = 
			\begin{bmatrix}
			0 & 0 & 0 & 0 & 0 & 0 \\
			\frac{2}{5} & 0 & 0 & 0 & 0 & 0 \\
			\frac{3 (500 M-409)}{1045} & -\frac{6 (250 M-309)}{1045} & 0 & 0 & 0 & 0 \\
			\frac{547008637842659863-386281780255161444 M}{152025995207353729} & \frac{3 (2856036493343421 M-3906629787402737)}{2288803570033750} & -\frac{741819 (1303514627 M-4489322119)}{2239523110391875} & -\frac{8391 (18565412379 M-24937027363)}{202099662773750} & 0 & 0 \\
			\end{bmatrix}, \qquad \lambda = 1, \ldots, M,$} \\
		${b}^{{\{\f\}}}=  
		\begin{bmatrix}
		\frac{2288000}{4732539} & -\frac{2203}{14250} & \frac{247273}{613500} & \frac{30767}{152250} & -\frac{1}{8} & \frac{191}{1000} \\
		\end{bmatrix}\trsym,$
		&
		$b^{{\{\s\}}} =
		\begin{bmatrix}
		\frac{11}{72} & \frac{25}{72} & \frac{25}{72} & \frac{11}{72} \\
		\end{bmatrix}\trsym,$
		\\
		$\widehat{b}^{{\{\f\}}}= 
		\begin{bmatrix}
		\frac{1837041228720545025825201951582239534326}{2195453146940870392428577778808091404375} & -\frac{12181532573386077454382848427541846123}{17628427274186874088578235825358750000} & \frac{4528991149246665992465958589958885289}{8624433917634487442857055565277187500} & \frac{83750160542686187}{606298988321250000} & \frac{191}{1000} & 0 \\
		\end{bmatrix}\trsym,$
		&
		$\widehat{b}^{{\{\s\}}} =
		\begin{bmatrix}
		\frac{1}{5} & \frac{1}{4} & \frac{3}{8} & \frac{7}{40} \\
		\end{bmatrix}\trsym.$
}

\begin{figure}[ht!]
	\centering
	\begin{subfigure}[b]{0.45\textwidth}
		\centering
		\includegraphics[width=\textwidth]{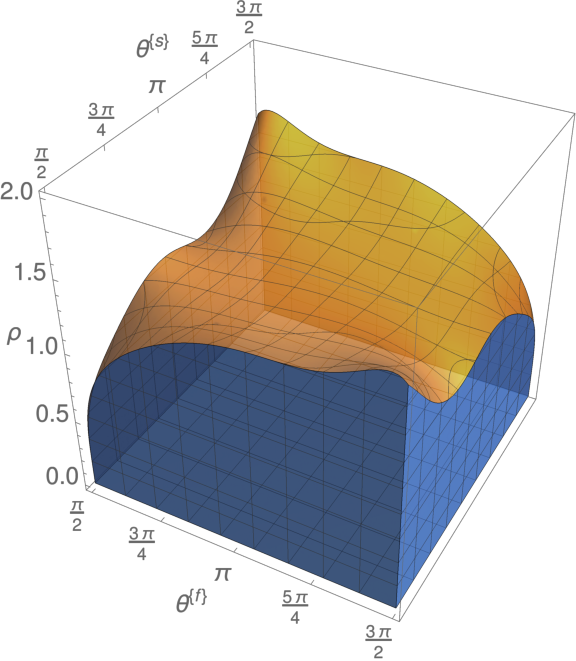}
		\caption{$M=2$.}
	\end{subfigure}
	\hfill
	\begin{subfigure}[b]{0.45\textwidth}
		\centering
		\includegraphics[width=\textwidth]{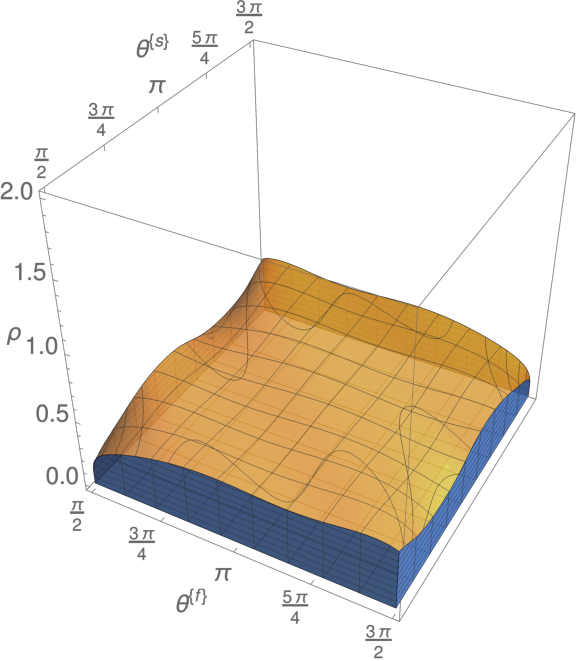}
		\caption{$M=4$.}
	\end{subfigure}
	\caption{\method{\implicit}{\explicit}{4}{2}{6}{4}{A} stability regions.}
\end{figure}

\end{document}